\documentclass[reqno]{amsart}

\usepackage{setspace}

\usepackage{tikz}

\usepackage{hyperref,todonotes,comment}
\usepackage[color,matrix,arrow]{xy}
\usepackage{graphics,amssymb}
\usepackage{pagecolor,lipsum}
\usepackage{latexsym}
\usepackage{mathrsfs}
\xyoption{curve}
\usepackage{hyperref}
\hypersetup{ 
colorlinks,
citecolor=cyan,
filecolor=cyan,
linkcolor=cyan,
urlcolor=cyan
}

\newtheorem{lemma}{Lemma}[section]
\newtheorem{proposition}[lemma]{Proposition}
\newtheorem{theorem}[lemma]{Theorem}
\newtheorem{corollary}[lemma]{Corollary}

\newtheorem{conjecture}[lemma]{Conjecture}
\newtheorem*{theoremA}{Theorem}

\theoremstyle{definition}
\newtheorem{example}[lemma]{Example}
\newtheorem{definition}[lemma]{Definition}
\newtheorem{remark}[lemma]{Remark}

\newcommand{\mfk}[1]{\mathfrak{#1}}
\newcommand{\mbb}[1]{\mathbb{#1}}
\newcommand{\mcl}[1]{\mathcal{#1}}

\newcommand{\msc}[1]{\mathscr{#1}}
\newcommand{\mbf}[1]{\mathbf{#1}}

\DeclareMathOperator{\Hom}{Hom}
\DeclareMathOperator{\End}{End}
\DeclareMathOperator{\RHom}{RHom}

\DeclareMathOperator{\Ext}{Ext}
\DeclareMathOperator{\Aut}{Aut}

\DeclareMathOperator{\rep}{rep}
\DeclareMathOperator{\Rep}{Rep}
\DeclareMathOperator{\Spec}{Spec}

\DeclareMathOperator{\Enh}{Enh}

\DeclareMathOperator{\Fr}{Fr}
\DeclareMathOperator{\Vect}{Vect}

\DeclareMathOperator{\ord}{ord}
\DeclareMathOperator{\Sym}{Sym}
\DeclareMathOperator{\FK}{FK}
\DeclareMathOperator{\res}{res}

\DeclareMathOperator{\Ind}{Ind}
\DeclareMathOperator{\coh}{coh}
\DeclareMathOperator{\Coh}{Coh}
\DeclareMathOperator{\IndCoh}{IndCoh}
\DeclareMathOperator{\QCoh}{QCoh}

\newcommand{\Xq}{X_q}
\newcommand{\dg}{\mathrm{dg}}
\newcommand{\de}{\operatorname{d\acute{e}}}

\newcommand{\sHom}{\mathscr{H}\! om}
\newcommand{\sExt}{\mathscr{E}xt}
\newcommand{\sRHom}{\mathscr{R\!H}\!om}

\newcommand{\ot}{\otimes}
\newcommand{\tN}{\tilde{\mathcal{N}}}

\newcommand{\uqG}{u(G_q)}
\newcommand{\uqB}{u(B_q)}
\newcommand{\su}{u(B_q)}
\newcommand{\dG}{\check{G}}
\newcommand{\dB}{\check{B}}
\renewcommand{\1}{\mathbf{1}}
\renewcommand{\O}{\mathscr{O}}
\renewcommand{\tilde}{\widetilde}

\definecolor{page_color}{HTML}{E1C16E}
\definecolor{text_color}{HTML}{fffff0}
\definecolor{myblue}{RGB}{255,0,255}

\title[]{The half-quantum flag variety and representations for small quantum groups}
\date{December 22, 2022}

\author{Cris Negron}
\address{Department of Mathematics, University of Southern California, Los Angeles, CA 90007}
\email{cnegron@usc.edu}

\author{Julia Pevtsova}
\address{Department of Mathematics, University of Washington, Seattle, WA 98195}
\email{julia@math.washington.edu}

\begin{document}

\maketitle

\begin{abstract}
Consider an almost-simple algebraic group $G$ and a choice of complex root of unity $q$.  We study the category of quasi-coherent sheaves $\msc{X}_q$ on the half-quantum flag variety, which itself forms a sheaf of tensor categories over the classical flag variety $G/B$.  We prove that the category of small quantum group representations for $G$ at $q$ embeds fully faithfully into the global sections of $\msc{X}_q$, and that the fibers of $\msc{X}_q$ over $G/B$ recover the tensor categories of representations for the small quantum Borels.  These relationships hold both at an abelian and derived level.  Subsequently, reduction arguments, from the small quantum group to its Borels, appear algebrogeometrically as ``fiber checking" arguments over $\msc{X}_q$.  We conjecture that $\msc{X}_q$ also contains the category of dg sheaves over the Springer resolution as a full monoidal subcategory, at the derived level, and hence provides a monoidal correspondence between the Springer resolution and the small quantum group.  We relate this conjecture to a known equivalence between dg sheaves on the Springer resolution and the principal block in the derived category of quantum group representations \cite{arkhipovbezrukavnikovginzburg04,bezrukavnikovlachowska07}.
\end{abstract}

\tableofcontents

\section{Introduction}

Let $G$ be an almost-simple algebraic group over the complex numbers and let $q\in \mbb{C}$ be a root of unity of order $>3$.  
We consider representations $\Rep G_q$ of Lusztig's divided power quantum algebra, and the associated small quantum group $\uqG$ for $G$ 
at $q$.   As we recall in Section \ref{sect:Gq}, at such a finite order parameter we have an associated dual group 
$\dG$ and central tensor functor $\Fr:\Rep \dG\to \Rep G_q$ provided by Lusztig's quantum Frobenius \cite{lusztig90II,lusztig93}.\footnote{The appearance of the 
dual group $\dG$ can be safely ignored for the purposes of this introduction}

This paper concerns a certain ``sheaf of tensor categories" over the flag variety
\begin{equation}\label{eq:1}
(\msc{X}_q=)\QCoh(X_q)
\end{equation}
and its relationship with some of the usual quantum suspects; the small (and big) quantum group and its quantum Borels, as well as the Springer resolution and nilpotent cone for $\dG$.  Formally, $\QCoh(X_q)$ is a monoidal category of relative Hopf modules for the quantum function algebra.  Informally, $X_q$ is the noncommutative quotient stack $X_q=\dG/B_q$ and sheaves on $X_q$ are $B_q$-equivariant sheaves on $\dG$.  Based on the latter interpretation, we refer to the category \eqref{eq:1} as the category of sheaves on the \emph{half-quantum flag variety} $X_q$.
\par

Our interest in the category $\QCoh(X_q)$ is based on the leverage it provides in analyses of the small quantum group, via the Borels, and for the (conjectural) role it plays in linking the tensor category of small quantum group representations to the tensor category of sheaves on the Springer resolution.  This first point is borne out explicitly in our subsequent text \cite[Part II]{negronpevtsova5}, where we apply the geometry of $\QCoh(X_q)$ to describe the tensor triangular geometry \cite{balmer10II} of the derived category of small quantum group representations.

We note that this basic reduction principle, which asserts that studies of $G$ or $G_q$ can be reduced to studies of the Borels, has been present in geometric representation theory at least since the middle of the 20th century.  The idea manifests itself already in the classical Beilinson-Bernstein localization theorem \cite{beilinsonbernstein81}, for example, as well as its versions in modular representation theory and for quantum groups (see, in particular, \cite{bezrukavnikovmirkovicrumynin08,arkhipovbezrukavnikovginzburg04,bezrukavnikovlachowska07} and others). The main difference between the goals 
of our analysis and the established localization theorems in the literature is that we aim to preserve the tensor structure of $\Rep \uqG$ when restricting to our family of Borels.

Now, at this point we have alluded to the small quantum Borels for $G$ at $q$ a number of times, but the reader may understand that only the positive and negative Borels for $\uqG$ have been constructed in the literature.  One of the first substantive claims of this text is that the (representations for the) positive and negative Borel in $\uqG$ expand naturally into a \emph{family} of tensor categories which are parametrized by points in the flag variety.  The members of this family are our \emph{small quantum Borels} $\msc{B}_\lambda$ for $G$ at $q$.  We show that the small quantum Borel over a given point $\lambda$ in $\dG/\dB$ is recovered as the fiber of $\QCoh(X_q)$ over $\lambda$.  In this way $\QCoh(X_q)$ can alternatively be understood as a category of representations for ``the universal small quantum Borel" associated to $G$ at $q$.
\par

One should compare here with the classical situation over $\mathbb{C}$ or $\bar{\mathbb{F}}_p$, where the flag variety parametrizes Borel subgroups $B^\lambda\subset G$, or subalgebras $\mfk{b}^\lambda\subset \mfk{g}$, and these subgroups/subalgebras naturally vary as a family over $G/B$.

The remainder of the introduction is dedicated to a more detailed description of the half-quantum flag variety, its fibers, and its connections to the Springer resolution.  We also provide a more precise accounting of the main results of the text.

\subsection{The half-quantum flag variety and small quantum Borels}
\label{sect:intro1}

An explicit description of the half-quantum flag variety can be found in Section \ref{sect:G/Bq_intro}.  However, one can think as follows: quantum Frobenius provides a translation action of the quantum Borel on the quantum group and we take
\[
\QCoh(X_q)=\QCoh(\dG)^{B_q}=\text{the $\ot$-category of $B_q$-equivariant sheaves on }\dG.
\]
The relationship between $\QCoh(X_q)$ and sheaves on the classical flag variety is then provided by a tensor embedding $\zeta^\ast:\QCoh(\dG/\dB)\to \QCoh(X_q)$, and we have a similar embedding from the category of quantum group representations,
\begin{equation}\label{eq:196}
\xymatrix{
 & \QCoh(X_q) \\
\Rep\uqG\ar[ur]^{\kappa^\ast} & & \QCoh(\dG/\dB).\ar[ul]_{\zeta^\ast}
}
\end{equation}
Informally, $\Rep \uqG$ is identified with the category of sheaves on the quotient stack $\dG/G_q$ \cite{arkhipovgaitsgory03}, at which point the maps in \eqref{eq:196} become pullback functors along the naturally occurring projections
\[
\kappa:X_q=\dG/B_q\to \dG/G_q\text{\ \ \ and\ \ \ }\zeta:X_q=\dG/B_q\to \dG/\dB.
\]
The above amalgamation of tensor categories forms the foundation for our analysis, and is described in detail in Sections \ref{sect:G/Bq} and \ref{sect:Kempf}.

In Section \ref{sect:borels} we also introduce a family of small quantum Borels parametrized by points in the flag variety $\lambda:\Spec(K)\to \dG/\dB$.  These small quantum Borels take the forms of tensor categories $\msc{B}_\lambda$ which come equipped with central tensor functors $\res_\lambda:\Rep(u(G_q))\to \msc{B}_\lambda$.  At the identity, we have an identification $\msc{B}_1\cong \Rep(u(B_q))$ of tensor categories over $\Rep(u(G_q))$, and at arbitrary $\lambda$ we have a $\dB_K$-torsor of tensor equivalences $\msc{B}_\lambda \cong \Rep(u(B_q))$.  In Section \ref{sect:fibers_G/B} we realize these quantum Borels as the {\it categorical fibers} of $\QCoh(\Xq)$ over the flag variety, as formulated in Section \ref{sect:flatfamily}. 

\begin{theoremA}[\ref{thm:calc_fib}/\ref{thm:derived_fibI}/\ref{thm:derived_fibII}]
Taking the fiber of the family $\QCoh(\Xq)$ at any geometric point $\lambda:\Spec(K)\to \dG/\dB$ recovers the corresponding small quantum Borel
\[
fib_\lambda:\QCoh(\Xq)|_\lambda\overset{\sim}\to \msc{B}_\lambda.
\]
Furthermore, the composite $fib_\lambda\circ \kappa^\ast:\Rep\uqG\to \msc{B}_\lambda$ recovers the restriction functor $\res_\lambda$.  The analogous calculations also hold at the level of derived categories.
\end{theoremA}

The above calculation allows us to think of $\QCoh(X_q)$ as the total space for the ``smoothly varying" family of categories $\msc{B}_\lambda$, and of $\kappa^\ast$ as a universal restriction functor.  The next theorem summarizes the fundamental properties of the functor $\kappa^\ast$.

\begin{theoremA}[\ref{thm:Kempf}]
The universal restriction functor $\kappa^\ast:\Rep(u(G_q))\to \QCoh(\Xq)$ is an exact, central, fully faithful monoidal embedding.  Furthermore the induced functor on unbounded derived categories
\[
\kappa^\ast:D(u(G_q))\to D(\Xq)
\]
remains fully faithful.
\end{theoremA}

Theorems \ref{thm:calc_fib} and \ref{thm:Kempf} provide an explicit articulation of the notion that the small quantum Borels collectively ``know everything" about the small quantum group.  From this perspective it is not surprising that one can deduce results for the small quantum group via coherent analyses of the small quantum Borels.  This is true of non-tensorial, block-by-block studies of the quantum group representations, as well as studies which emphasize the tensor structure of $\Rep\uqG$ (see e.g.\ \cite{bezrukavnikovmirkovicrumynin08,arkhipovbezrukavnikovginzburg04,bezrukavnikovlachowska07} \cite[Part II]{negronpevtsova5}).

\subsection{The Springer resolution and the half-quantum flag variety}
\label{sect:intro2}

In the previous subsection we discussed $\QCoh(X_q)$ and its relationship, via functors, to the flag variety, small quantum group, and small quantum Borels.  In this subsection we focus on a related enhancement for $\QCoh(X_q)$.  Specifically, the inner morphisms for the action of $\QCoh(\dG/\dB)$ on $\QCoh(X_q)$ produce sheaves $\sHom_{X_q}(M,N)$ 
on $\dG/\dB$ which vary naturally in each factor and admit an identification on the global sections
\[
\Gamma(\dG/\dB,\sHom_{X_q}(M,N))=\Hom_{X_q}(M,N).
\]
\par

Furthermore, the fact that $\zeta^\ast$ is monoidal implies that the localized morphisms $\sHom_{X_q}$ are compatible with the monoidal structure on $\QCoh(X_q)$.  So we obtain a monoidal enhancement $\QCoh^{\Enh}(X_q)$ of $\QCoh(X_q)$ over the flag variety.  This internal sheafy structure naturally derives, so that we have derived sheaf-maps $\sRHom_{X_q}$ and a subsequent enhancement $D^{\Enh}(X_q)$ of the derived category of sheaves on $X_q$.  Again, these derived sheaf-morphisms respect the monoidal structure on $D(X_q)$.  Of particular interest here is the dg algebra
\[
\sRHom_{X_q}(\1,\1)\ \ \in\ \ D(\dG/\dB).
\]

We recall that the Springer resolution $\tN$ is the cotangent bundle for the flag variety $\tN=T^\ast\dG/\dB$, and so we have the affine morphism $p:\tN\to \dG/\dB$.  The Springer resolution is also understood as a distinguished resolution of singularities $\mu:\tN\to \mcl{N}$ for the nilpotent cone in $\operatorname{Lie}(\dG)$.  In the following statement $h$ denotes the Coxeter number for $G$.

\begin{theorem}[\ref{thm:enh_GK}]
Suppose that $q$ is of odd order $\ord(q)>h$, or that $G$ is of type $A_1$.  Then there is an identification of algebras
\begin{equation}\label{eq:261}
H^\ast(\sRHom_{\Xq}(\1,\1))=p_\ast\O_{\tN}.
\end{equation}
\end{theorem}

The calculation of Theorem \ref{thm:enh_GK} can be used immediately to relate the quantum group to the Springer resolution, at least in analyses of tensor triangular geometry for the small quantum group \cite[Part II]{negronpevtsova5}.  However, we view this result as a reflection of a more fundamental relationship between the Springer resolution, half-quantum flag variety, and the (small) quantum group.

\subsection{A formality conjecture}
\label{sect:intro3}

To conclude the introduction, and the paper, we propose a deeper relationship between the half-quantum flag variety and the Springer resolution.

\begin{conjecture}[{Strong formality conjecture \ref{conj:sfc}}]
There is a $\QCoh_{\dg}(\dG/\dB)$-linear, fully faithful, central tensor functor
\[
\eta^\ast:\QCoh_{\dg}(\tN)\to \IndCoh_{\dg}(X_q),
\]
and corresponding fully faithful tensor functor $\QCoh_{\dg}(\mcl{N})\to \Rep_{\dg}(\uqG)$.
\end{conjecture}

In the above statement $\QCoh_{\dg}(Y)$ denotes the derived $\infty$-category of quasi-coherent sheaves on a given scheme $Y$, $\IndCoh_{\dg}(X_q)$ is the formal cocompletion of the $\infty$-category of coherent dg sheaves on $X_q$, and $\Rep_{\dg}(\uqG)$ is the analogous cocompletion of the $\infty$-category of finite-dimensional dg representations for the small quantum group.  Also $\tN$ should be understood here as a dg scheme over $\dG/\dB$, with generators in cohomological degree $2$ (cf.\ \cite{arkhipovbezrukavnikovginzburg04,bezrukavnikovlachowska07}).
\par

This conjecture is informed by the construction of the enhanced derived category for $\Xq$, provided in Sections \ref{sect:D_Enh} and \ref{sect:RGq}, the calculation of Theorem \ref{thm:enh_GK}, and known formality results from Arkhipov, Bezrukavnikov, and Ginzburg \cite{arkhipovbezrukavnikovginzburg04}.

We note that, from Theorem \ref{thm:Kempf}, we have a fully faithful central embedding $\kappa^\ast\colon~\Rep_{\dg}(u(G_q))\to \IndCoh_{\dg}(X_q)$, so that Conjecture \ref{conj:sfc} realizes a correspondence between braided monoidal categories
\begin{equation}\label{eq:266}
\xymatrixcolsep{3mm}
\xymatrixrowsep{3mm}
\xymatrix{
 & \IndCoh_{\dg}(X_q) \\
\QCoh_{\dg}(\tN)\ar[ur]^{\eta^\ast} & & \Rep_{\dg}(u(G_q)).\ar[ul]_{\kappa^\ast}
}
\end{equation}
We furthermore propose the following.

\begin{conjecture}
The push-pull functor along the monoidal correspondence \eqref{eq:266} restricts to an equivalence
\[
\eta_\ast\kappa^\ast:\operatorname{PrinBlock}_{\dg}(u(G_q))\overset{\sim}\longrightarrow \QCoh_{\dg}(\tN).
\]
This equivalence precisely recovers that of Arkhipov-Bezrukavnikov-Ginzburg and Bezrukavnikov-Lachowska \cite{arkhipovbezrukavnikovginzburg04,bezrukavnikovlachowska07}.
\end{conjecture}

Hence the correspondence \eqref{eq:266}, conjecturally, provides a monoidal reconceptualization of the non-monoidal equivalences of \cite{arkhipovbezrukavnikovginzburg04,bezrukavnikovlachowska07}.
\par

Under appropriate centralizing hypotheses \eqref{eq:266} can be seen as a morphism between sheaves on the Springer resolution and the derived category of quantum group representations in the monoidal $4$-category of presentable braided tensor $\infty$-categories, as constructed in \cite{johnsonscheimbauer17}.  Although we won't elaborate on the point here, the latter framing is relevant when one considers the Springer resolution in relation to topological field theories one constructs from the quantum group.  Here one \emph{hopes} to employ the Springer resolution to resolve certain singularities which appear in these TQFTs (cf.\ \cite[\S 19]{negronpevtsova5}).  For additional context we encourage the reader to peruse the texts \cite{costellocreutziggaitto19,schweigertwoike21,gukovetal21,derenzietal,creutzigdimoftegarnergeer}, among others.

\subsection{Acknowledgements}

Thanks to David Ben-Zvi, Roman Bezrukavnikov, Ken Brown, Eric Friedlander, Sergei Gukov, David Jordan, Simon Lentner, Sunghyuk Park, Noah Snyder, and Lukas Woike for useful discussions which influenced the trajectory of this work.  C.\ Negron is supported by NSF grant DMS-2001608 and Simons Collaboration Grant no.\ 999367. J.\ Pevtsova is supported by NSF grants DMS-1901854, 2200832, and the Brian and Tiffinie Pang faculty fellowship. This material is based upon work supported
by the National Science Foundation under Grant No. DMS-1440140, while the first
author was in residence at the Mathematical Sciences Research Institute in Berkeley,
California, and the second author was in digital residence. This research was also funded by the Deutsche Forschungsgemeinschaft (DFG, German Research Foundation) under Germany's Excellence Strategy -- EXC-2047/1 -- 390685813 when the authors visited the Hausdorff Institute for Mathematics in Bonn during the trimester on ``Spectral methods in Algebra, Geometry and Topology".

\tableofcontents


\section{Notation guide and categorical generalities}
\label{sect:notation}

Throughout $k$ is an algebraically closed field of characteristic $0$, and $G$ is an almost-simple algebraic group over $k$.  We let $h$ denote the Coxeter number for $G$, and $B\subset G$ is a fixed choice of Borel, which we recognize as the positive Borel in $G$.
\begin{itemize}
\item $q\in k$ is a root of unity of finite order.  We take $l=\ord(q^2)$ and assume that $\operatorname{ord}(q)$ is not extraordinarily small, as specifically articulated in \cite[Section 35.1.2]{lusztig93}.  When $q$ is of even order we assume that the character lattice for $G$ is strongly admissible, in the sense of \cite[Definition 3.1]{negron21}.\vspace{2mm}

\item $\dG$ is the dual group to $G$ at the given parameter $q$, as defined in Section \ref{sect:Gq}, with corresponding positive Borel $\dB$.  The group $\dG$ is referred to as the metaplectic dual for $G$ in (quantum) geometric Langlands \cite{gaitsgorylysenko18}.\vspace{2mm}

\item $\pi:\dG\to \dG/\dB$ is the quotient map.\vspace{2mm}

\item $\mcl{N}$ is the nilpotent cone in $\operatorname{Lie}(\dG)$, and $\mu:\tN\to \mcl{N}$ is the Springer resolution.  Alternatively, $\tN$ is the cotangent bungle for the dual group $\tN=T^\ast(\dG/\dB)$ and $p:\tN\to \dG/\dB$ is the associated bundle map.\vspace{2mm}

\item A geometric point $x:\Spec(K)\to Y$ in a ($k$-)scheme $Y$ is a map of schemes from the spectrum of an algebraically closed field extension $K/k$.\vspace{2mm}

\item $\Vect$ denotes the category of arbitrary $k$-linear vector spaces and, for any field extension $k\to K$, $\Vect(K)$ denotes the category of arbitrary $K$-linear vector spaces.\vspace{2mm}

\item The symbol $\ot_k$ denotes a vector space tensor product, and the generic symbol $\ot$ denotes the product operation in a given monoidal category $\msc{C}$.  So this product $\ot$ can be a product of sheaves, or a product of group representations, etc.  We also let $\ot_k$ denote the implicit action of $\Vect$ on a linear category.\vspace{2mm}

\item $\1$ is the unit object in a given tensor category $\msc{C}$.
\end{itemize}

We use the term \emph{tensor category} somewhat informally, to indicate a linear monoidal category which is of an algebraic origin (cf.\ \cite[Definition 1.1]{delignemilne82} \cite[Definition 4.1.1]{egno15}).  By a finite tensor category, however, we always mean a finite tensor category in the sense of \cite{etingofostrik04}.

\begin{remark}
The strong admissibility condition, at even order $q$, is inessential.  We only avoid non-admissible even order $q$ because appropriate treatments for the quantum group at such parameters have not appeared in the literature.
\end{remark}

\subsection{Finite-dimensional vs.\ infinite-dimensional representations}
\label{sect:finite}

For our study it has been convenient to work with cocomplete categories, where one generally has enough injectives and can freely use representability theorems.  For this reason we employ categories of arbitrary (possibly infinite-dimensional) representations $\Rep A$ for a given Hopf algebra $A$.  This is, by definition, the category of locally finite representations, or equivalently representations which are the unions of their finite-dimensional subrepresentations.
\par

One recovers the category $\rep A$ of finite-dimensional representations as the subcategory of \emph{dualizable}, or \emph{rigid}, objects in $\Rep A$, i.e.\ objects which admit left and right duals \cite{egno15}.  Since tensor functors preserve dualizable objects \cite[Ex. 2.10.6]{egno15}, restricting to the subcategory of dualizable objects $\rep A\subset \Rep A$ is a natural operation with respect to tensor functors.  In this way one moves freely between the ``small" and ``big" representation categories for $A$.
\par

When we work with derived categories, we take
\[
D(A):=\left\{\begin{array}{c}\text{the unbounded derived category of (generally}\\
\text{infinite-dimensional) $A$-representations}
\end{array}\right\}.
\]
We have the distinguished subcategories $D^b(A)$, $D^+(A)$, etc.\ of appropriately bounded complexes of (generally infinite-dimensional) representations, and a distinguished subcategory of bounded complexes of finite-dimensional representations, which one might write as $D^b(\rep A)$ or $D_{fin}(A)$.  In following the philosophy proposed above, $D_{fin}(A)$ appears as the subcategory of dualizable objects in $D(A)$.

\subsection{Presentable categories}
\label{sect:presentable}

As just stated, we generally work with cocomplete categories in this text.  However, at times we need to be more clear about which \emph{types} of cocomplete categories we consider.  When necessary, we restrict our attention to presentable categories.

We recall that a category $\msc{C}$ is called presentable, or locally presentable, if it is cocomplete and sufficiently compactly generated.  In particular, we assume the existence of a regular cardinal $\tau$ so that $\msc{C}$ is generated by a set of $\tau$-compact object $\{X_i\}_{i\in I}$ under $\tau$-filtered colimits.  This latter condition is called $\tau$-accessibility \cite[Definitions 1.17 \& 2.1]{adamekrosicky94}.
\par

For example, if $\msc{C}$ is cocomplete and compactly generated by an essentially small subcategory of compact objects, then $\msc{C}$ is presentable.  Indeed, in this case $\msc{C}$ is $\tau$-accessible with $\tau=\aleph_0$.  Every (abelian) category we come into direct contact with satisfies such compact generation hypotheses, and is hence presentable.

\subsection{Central functors}

By a central monoidal functor between monoidal categories, we mean a monoidal functor $F:\msc{C}\to \msc{D}$ from a braided monoidal category $\msc{C}$ to a possibly non-braided monoidal category $\msc{D}$ which comes equipped with a natural half-braiding
\[
\gamma_{-,X}:-\ot F(X)\overset{\cong}\to F(X)\ot -
\]
at each object $X$ in $\msc{C}$.  We require these half-braidings to be natural in each coordinate, compatible with the braiding on $\msc{C}$, and associative.  Equivalently, a central monoidal functor is a choice of a monoidal functor $F:\msc{C}\to \msc{D}$ and a choice of a braided monoidal lift of $F$ to the Drinfeld center $\tilde{F}:\msc{C}\to Z(\msc{D})$ \cite[Definition 4.15]{dgno10}.

\subsection{Relative Hopf modules}

For a Hopf algebra $A$, and an $A$-comodule algebra $\O$, we let $_\O\mbf{M}^A$ denote the category of relative $(\O,A)$-Hopf modules, with no finiteness assumptions.  This is the category of simultaneous (left) $\O$-modules and (right) $A$-comodules $M$, for which the coaction $M\to M\ot_k A$ is a map of $\O$-modules \cite[\S 8.5]{montgomery93}.  Here $\O$ acts diagonally on the product $M\ot_k A$, $a\cdot(v\ot b)=a_1v\ot a_2b$.
\par

For a basic example, we consider an affine algebraic group $H$ acting on (the right of) an affine scheme $Y$.  This gives the algebra of functions $\O=\O(Y)$ a comodule structure over $A=\O(H)$.  The fact that the action map $Y\times H\to Y$ is a scheme map says that $\O$ is an $A$-comodule algebra under this comodule structure.  Relative Hopf modules are then identified with equivariant sheaves on $Y$ via the global sections functor,
\[
\Gamma(Y,-):\QCoh(Y)^H\overset{\sim}\to {_{\O(Y)}\mbf{M}^{\O(H)}}.
\]

\subsection{Descent along $H$-torsors}
\label{sect:descent}

Suppose an algebraic group $H$ acts on a scheme $Y$, that the quotient $Y/H$ exists, and that the quotient map $\pi:Y\to Y/H$ is a (faithfully flat) $H$-torsor.  For example, we can consider an algebraic group $G$ and a closed subgroup $H\subset G$ acting on $G$ by translation.  In this case the quotient $G\to G/H$ exists, is faithfully flat, and realizes $G$ as an $H$-torsor over the quotient \cite[Theorem B.37]{milne17}.
\par

For $Y$ as prescribed, pulling back along the quotient map $\pi:Y\to Y/H$ defines an equivalence of categories
\begin{equation}\label{eq:420}
\pi^\ast:\QCoh(Y/H)\overset{\sim}\to \QCoh(Y)^H
\end{equation}
from sheaves on the quotient to $H$-equivariant sheaves on $Y$.  The inverse to this equivalence if provided by faithfully flat descent, or simply ``descent"
\begin{equation}\label{eq:desc}
\operatorname{desc}:\QCoh(Y)^H\to \QCoh(Y/H).
\end{equation}
For details on faithfully flat descent one can see \cite[Expos\'{e} VIII]{sga1} or \cite[\href{https://stacks.math.columbia.edu/tag/03O6}{Tag 03O6}]{stacks}.  One may understand descent simply as the inverse to the equivalence \eqref{eq:420}, which we are claiming exists under the precise conditions outlined above.
\par

Now, for \emph{any} scheme $Y$ with an $H$-action, we have a pair of adjoint functors 
\[
E_{-}:  \Rep H \to \QCoh(Y)^H\ \ \text{and}\ \ -|_{H}: \QCoh(Y)^H \to \Rep H,
\]
where for $V \in \Rep H$ we set $E_V$ to be the vector bundle on $Y$ corresponding to the $H$-equivariant $\O_Y$-module $\O_Y \otimes_k V$, with $\O_Y$ acting on the left factor and $H$ acting diagonally. By construction, $E_V$ is locally free; it is coherent if and only if $V$ is finite-dimensional. The right adjoint $-|_{H}$ is given by taking  global section $\Gamma(Y,-)$ and then forgetting the action of $\O(Y)$.
\par

For an affine algebraic group $Y = G$ and a closed subgroup $H \subset G$, the descent of the equivariant vector bundle $E_V$, defined as above, is the familiar bundle from \cite[I.5.8]{jantzen03}, \cite[3.3]{suslinfriedlanderbendel97b}.

\subsection{Flat families of tensor categories}
\label{sect:flatfamily}

Let $X$ be a (quasi-compact quasi-separated) scheme.  By a \emph{flat family of tensor categories} over $X$ we mean a presentable monoidal category $\msc{C}$ equipped with an exact, cocontinuous, central monoidal functor $w:\QCoh(X)\to \msc{C}$.  We suppose additionally that the product on $\msc{C}$ commutes with small coloimts in each variable, and that $\msc{C}$ is generated by a subcategory of compact dualizable objects.
\par

The fiber of the family $\msc{C}$ along a given map $f:Y\to X$ is the categorical base change
\[
\msc{C}|_f:=\QCoh(Y)\ot_{\QCoh(X)}\msc{C},
\]
as described in Section \ref{sect:base_change}.  Our notion of a flat family of tensor categories is related to a notion of a sheaf of tensor categories, as defined for example in \cite{gaitsgory15}, in the sense that the derived ($\infty$-)category of a flat family of categories over $X$ has the structure of a sheaf of categories over $X$.

\begin{remark}
Our flatness condition, i.e.\ exactness of the structure map $w$, is imposed to create some consistency in calculating fibers at the non-derived level, for families of categories, and the derived level, for sheaves of categories (cf.\ \cite{benzvifrancisnadler10}).
\end{remark}

\subsection{Enriched categories and enhancements}
\label{sect:enrich}

A flat family of tensor categories over a scheme $X$, as in Subsection \ref{sect:flatfamily}, gives a specific example of a category with an enhancement in $\QCoh(X)$. We now describe this (generally weaker) structure in detail. 

By a category $\mcl{T}$ enriched in a given monoidal category $\operatorname{Q}$ we mean a collection of objects $\operatorname{obj}(\mcl{T})$ for $\mcl{T}$ and morphisms $\mcl{H}om_{\mcl{T}}(M,N)$, which are objects in $\operatorname{Q}$, for each pair of objects $M$ and $N$ in $\mcl{T}$.  We suppose, additionally, the existence of associative composition morphisms
\[
\circ:\mcl{H}om_{\mcl{T}}(M,N)\ot\mcl{H}om_{\mcl{T}}(L,M)\to\mcl{H}om_{\mcl{T}}(L,N)
\]
for each triple of objects in $\mcl{T}$.  Basics on enriched categories can be found in \cite[Ch 3]{riehl14}, though we recall some essential points here.
\par

Given any lax monoidal functor $F:\operatorname{Q}\to \operatorname{S}$ we can push forward the morphisms in $\mcl{T}$ along $F$ to get a new category $F\mcl{T}$ which is enriched in $\operatorname{S}$ \cite[Lemma 3.4.3]{riehl14}.  As one expects, the objects in $F\mcl{T}$ are the same as those in $\mcl{T}$, and the morphisms in $F\mcl{T}$ are given as $F\mcl{H}om_{\mcl{T}}(M,N)$.  The composition maps for $F\mcl{T}$ are induced by those of $\mcl{T}$ and the lax monoidal structure on the functor $F$.
\par

Of particular interest is the ``global sections" functor $\Gamma_{\operatorname{Q}}=\Hom_{\operatorname{Q}}(\1,-):\operatorname{Q}\to \operatorname{Set}$, with lax monoidal structure on $\Gamma_{\operatorname{Q}}$ provided by the monoidal and unit structure maps for $\operatorname{Q}$,
\[
\Gamma_{\operatorname{Q}}(A)\times \Gamma_{\operatorname{Q}}(B)\to \Gamma_{\operatorname{Q}}(A\ot B),\ \ (f,g)\mapsto f\ot g.
\]
For any enriched category $\mcl{T}$ over $\operatorname{Q}$ the global sections $\Gamma_{\operatorname{Q}}\mcl{T}$ are then an ordinary category \cite[Definition 3.4.1]{riehl14}.  Note that the category of global sections $\Gamma_{\operatorname{Q}}\mcl{T}$ also acts on $\mcl{T}$, in the sense that any morphism $f\in \Gamma_{\operatorname{Q}}\mcl{H}om_{\mcl{T}}(M,N)$ specifies composition and precomposition maps
\[
f_\ast:\mcl{H}om_{\mcl{T}}(L,M)\to \mcl{H}om_{\mcl{T}}(L,N)\ \ \text{and}\ \ f^\ast:\mcl{H}om_{\mcl{T}}(N,L)\to \mcl{H}om_{\mcl{T}}(M,L)
\]
via the unit structure on $\operatorname{Q}$ and composition in $\mcl{T}$.
\par

Suppose now that $\operatorname{Q}$ is symmetric.  By a \emph{monoidal} category enriched in $\operatorname{Q}$ we mean an enriched category $\mcl{T}$ with a product structure on objects, a unit object, and an associator, where the unit and associator structure maps appear as global isomorphisms for $\mcl{T}$.  We require the existence of tensor maps
\[
tens^{\mcl{T}}:\mcl{H}om_{\mcl{T}}(M,N)\ot \mcl{H}om_{\mcl{T}}(M',N')\to \mcl{H}om_{\mcl{T}}(M\ot M',N\ot N')
\]
which are associative relative to the associators on $\operatorname{Q}$ and $\mcl{T}$, and appropriately compatible with composition.  This compatibility between the tensor and composition morphisms appears as an equality
\[
\begin{array}{l}
(g_1\ot g_2)\circ (f_1\ot f_2)=(g_1\circ f_1)\ot (g_2\circ f_2):\vspace{1mm}\\
\hspace{2cm}\mcl{G}_1\ot \mcl{F}_1\ot \mcl{G}_2\ot\mcl{F}_2\to \mcl{H}om_{\mcl{T}}(L_1\ot L_2,N_1\ot N_2)
\end{array}
\]
for maps $f_i:\mcl{F}_i\to \mcl{H}om_{\mcl{T}}(L_i,M_i)$ and $g_i:\mcl{G}_i\to \mcl{H}om_{\mcl{T}}(M_i,N_i)$.  Such a monoidal structure on $\mcl{T}$ implies a monoidal structure on its category of global sections.
\par

\begin{definition}
An enhancement of a (monoidal) category $T$ is a choice of a $Q$-enriched (monoidal) category $\mcl{T}$ and a (monoidal) equivalence $T\cong \Gamma_{\operatorname{Q}}\mcl{T}$.
\end{definition}

\section{Quantum groups}
\label{sect:quantumgroups}

We recall basic constructions and results for quantum groups.  We also recall a geometric (re)construction of the small quantum group via de-equivariantization along the quantum Frobenius functor.

\subsection{Lusztig's divided power algebra}

We briefly review Lusztig's divided power algebra, leaving the details to the original works \cite{lusztig90II,lusztig93}.  Let $\mfk{g}$ be a semisimple Lie algebra over the complex numbers, and fix some choice of simple roots for $\mfk{g}$.
\par

To begin, one considers the generic quantum universal enveloping algebra
\[
U^{\operatorname{gen}}_v(\mfk{g})=\frac{\mbb{Q}(v)\langle E_\alpha, F_\alpha, K^{\pm 1}_\alpha:\alpha\text{ is a simple root for }\mfk{g}\rangle}{(v\text{-analogs of Serre relations})},
\]
where the $v$-Serre relations are as described in \cite[\S 1.1]{lusztig90II}.  These relations can be written compactly as
\[
\begin{array}{c}
\operatorname{ad}_v(E_\alpha)^{1-\langle \alpha,\beta\rangle}(E_\beta)=0,\ \ \operatorname{ad}_v(F_\alpha)^{1-\langle \alpha,\beta\rangle}(F_\beta)=0,\ \ K_\beta E_\alpha K_\beta^{-1}=v^{(\alpha,\beta)}E_\alpha,\vspace{2mm}\\
K_\beta F_\alpha K_\beta^{-1}=v^{-(\alpha,\beta)}F_\alpha,\ \ [E_\alpha,F_\beta]=\delta_{\alpha\beta}(K_\alpha+K_\alpha^{-1})/(v-v^{-1}).
\end{array}
\]
The algebra $U^{\operatorname{gen}}_v(\mfk{g})$ admits a Hopf algebra structure over $\mbb{Q}(v)$ with coproduct
\begin{equation}\label{eq:424}
\Delta(E_\alpha)=E_\alpha\ot K_\alpha +1\ot E_\alpha,\ \ \Delta(F_\alpha)=F_\alpha\ot 1+K_\alpha^{-1}\ot F_\alpha,\ \ \Delta(K_\alpha)=K_\alpha\ot K_\alpha.
\end{equation}
We consider the distinguished subalgebra of ``Cartan elements"
\[
U^{\operatorname{gen}}_v(\mfk{g})^0=\mbb{Q}(v)\langle K_\alpha,K_\alpha^{-1}:\alpha\text{ simple}\rangle\subset U^{\operatorname{gen}}_v(\mfk{g}),
\]
which we refer to as the toral subalgebra in $U^{\operatorname{gen}}_v(\mfk{g})$.
\par

In $U^{\operatorname{gen}}_v(\mfk{g})$ one has the $\mbb{Z}[v,v^{-1}]$-subalgebra $U_v(\mfk{g})$ generated by the divided powers
\[
E_\alpha^{(i)}=E^i_\alpha/[i]_{d_\alpha}!,\ \ F_\alpha^{(i)}=F^i_\alpha/[i]_{d_\alpha}!,\ \ \text{and }K_\alpha^{\pm 1}.
\]
Here $d_\alpha=|\alpha|^2/|\text{short\ root}|^2$ and $[v]_{d_\alpha}!$ is the $v^{d_i}$-factorial (see \cite[\S 1 and Theorem 6.7]{lusztig90II}).  The subalgebra $U_v(\mfk{g})$ furthermore has a Hopf structure induced by that of $U^{\operatorname{gen}}_v(\mfk{g})$, which is given by the same formulas \eqref{eq:424}.  We define the toral subalgebra as the intersection $U_v(\mfk{g})\cap U^{\operatorname{gen}}_v(\mfk{g})^0$.  This subalgebra not only contains powers of the $K_\alpha$, but also certain divided powers in the $K_\alpha$.
\par

For any value $q\in \mbb{C}^\times$ we have the associated $\mbb{Z}$-algebra map $\phi_q:\mbb{Z}[v,v^{-1}]\to \mbb{C}$ which sends $v$ to $q$, and we change base along $\phi_q$ to obtain Lusztig's divided power algebra
\[
U_q(\mfk{g})=\mbb{C}\ot_{\mbb{Z}[v,v^{-1}]}U_v(\mfk{g})
\]
at parameter $q$.  An important aspect of this algebra is that, at $q$ with $\ord(q^2)=l$, we have
\[
E_\alpha^l=[l]_{d_\alpha}!E_\alpha^{(l)}=0\ \ \text{and}\ \ F_\alpha^l=[l]_{d_\alpha}!F_\alpha^{(l)}=0.
\]
Hence the elements $E_\alpha$ and $F_\alpha$ become nilpotent in $U_q(\mfk{g})$.  One should compare with the distribution algebra associated to an algebraic group in finite characteristic.
\par

In Section \ref{sect:borel} we also consider the divided power algebra $U_q(\mfk{b})$ for the Borel.  This is the subalgebra $U_q(\mfk{b})\subset U_q(\mfk{g})$ generated by the toral elements as well as the divided powers $E_\alpha^{(i)}$ for the positive simple roots.

\subsection{Big quantum groups}
\label{sect:Gq}

We follow \cite{lusztig93} \cite[\S 3]{andersenparadowski95}.  Fix $G$ almost-simple with associated character lattice $X$.  Let $q$ be of finite order and take $\ord(q^2)=l$.  We adopt the specific restrictions on $q$ laid out in Section \ref{sect:notation}, which are almost nothing if one is happy to assume $\ord(q)$ is odd.
\par

We consider the category of representations for the ``big" quantum group associated to $G$,
\[
\rep G_q:=\left\{\begin{array}{c}
\text{finite-dimensional representations for $U_q(\mfk{g})$}\\
\text{which are compatibly graded by the character lattice }X
\end{array}\right\}.
\]
Here when we say a representation $V$ is graded by the character lattice we mean that $V$ decomposes into eigenspaces $V=\oplus_{\lambda\in X} V_\lambda$ for the action of the toral subalgebra in $U_q(\mfk{g})$, with each $V_\lambda$ of the expected eigenvalue.  For example one requires $K_\alpha\cdot v=q^{d_\alpha\langle \alpha,\lambda\rangle}v$ for homogeneous $v\in V_\lambda$ (see \cite[\S 3.6]{andersenparadowski95}).  The category $\rep G_q$ is the same as the category of representations of Lusztig's modified algebra $\dot{U}_q(G)$ \cite{lusztig93} \cite[\S 1.2]{kashiwara94}.
\par

We let $\Rep G_q$ denote the category of generally infinite-dimensional, integrable, $G_q$-representations.  Equivalently, $\Rep G_q$ is the category of $U_q(\mfk{g})$-modules which are the unions of their finite-dimensional $X$-graded submodules.  Note that all objects in $\Rep G_q$ remain $X$-graded.

Let $P$ and $Q$ denote the weight and root lattices for $G$ respectively, $\Phi$ denote the roots in $Q$, and let $(-,-):P\times P\to \mbb{Q}$ denote the normalized Killing form which takes value $2=(\alpha,\alpha)$ on short roots $\alpha$.  For $r=\exp(P/Q)$ we choose an $r$-th root $\sqrt[r]{q}$ of $q$ and define $q^{(\lambda,\mu)}=(\sqrt[r]{q})^{r(\lambda,\mu)}$.  For the sake of specificity, we also assume that $\ord(\sqrt[r]{q})=r\cdot\ord(q)$, and note that such roots exist at arbitrary $q$.  We then have the associated $R$-matrix for $\Rep G_q$ which appears as
\begin{equation}\label{eq:R}
\begin{array}{rl}
R&=(\sum_{n:\Phi^+\to \mbb{Z}_{\geq 0}} \operatorname{coeff}(q,n)E_{\gamma_1}^{n_1}\dots E_{\gamma_r}^{n_r}\ot F_{\gamma_1}^{n_1}\dots F_{\gamma_r}^{n_r})\Omega^{-1}\vspace{1mm}\\
&=\Omega^{-1}+\text{higher order terms}.
\end{array}
\end{equation}
Here $\Omega$ is the global semisimple operator
\[
\Omega=\sum_{\lambda,\mu\in X}q^{(\lambda,\mu)}1_\lambda\ot 1_\mu
\]
with each $1_\lambda$ the natural projection $1_\lambda:V\to V_\lambda$ \cite[Ch.\ 32]{lusztig93} \cite[\S 1]{sawin06}.  We note that the sum in \eqref{eq:R} has only finitely many non-vanishing terms, since the $E_\gamma$ and $F_\gamma$ are nilpotent in $U_q(\mfk{g})$.  The operator $R$ endows $\Rep G_q$ with its standard braiding
\[
\begin{array}{l}
c_{V,W}:V\ot W\to W\ot V,\vspace{1mm}\\
c_{V,W}(v,w)=q^{(\deg(w),\deg(v))}(w\ot v+\sum_{n>0}\operatorname{coeff}(q,n) F_{\gamma_1}^{n_1}\dots F_{\gamma_r}^{n_r}w\ot E_{\gamma_1}^{n_1}\dots E_{\gamma_r}^{n_r}v).
\end{array}
\]

We consider the dual group $\dG$, which is of the same Dynkin type as $G$ at odd order $q$ and of Langlands dual type when $q$ is of even order.  This dual group has character lattice $X^M$, where
\[
X^M:=\{\mu\in X:(\mu,\lambda)\in l\mbb{Z}\ \text{for all }\lambda\in X\},
\]
and form $(-,-)^{\vee}=\frac{1}{l^2}(-,-)$.  The roots for $\dG$ are given by the scalings $\{l_\gamma\gamma:\gamma\in \Phi\}$ \cite[Section 5.1]{negron21} \cite[Section 2.2.5]{lusztig93}.  So, for example, when $G$ is simply-connected and $q$ is of odd order, the dual $\dG$ is just the adjoint form for $G$.
\par

Via Lusztig's quantum Frobenius \cite[Ch.\ 35]{lusztig93} we have a Hopf map $fr^\ast:\dot{U}_q(G)\to \dot{U}(\dG)$,
\[
fr^\ast=\left\{
\begin{array}{l}
E_\gamma,\ F_\gamma\mapsto 0\vspace{1mm}\\
E_\gamma^{(l_\gamma)}\mapsto e_\gamma\vspace{1mm}\\
F_\gamma^{(l_\gamma)}\mapsto f_\gamma\vspace{1mm}\\
1_\lambda\mapsto 1_\lambda\text{ when $\lambda\in X^M$ and $0$ otherwise}
\end{array}\right.
\]
which defines a braided tensor embedding
\[
\Fr:\Rep \dG\to \Rep G_q
\]
whose image is the M\"uger center in $\Rep G_q$ \cite[Theorem 5.3]{negron21}, i.e.\ the full tensor subcategory of all $V$ in $\Rep G_q$ for which $c_{-,V}c_{V,-}:V\ot-\to V\ot-$ is the identity transformation.

\begin{remark}\label{rem:fr}
Our Frobenius map $fr^\ast:\dot{U}_q(G)\to \dot{U}_q(\dG)$ is induced by that of \cite{lusztig93}, but is not precisely the map of \cite{lusztig93}.  Similarly, our dual group $\dG$ is not precisely the dual group $G^\ast$ from \cite{lusztig93}.  Specifically, Lusztig's dual group $G^\ast$ is defined by taking the dual lattice $X^\ast\subset X$ to consist of all $\mu$ with restricted pairings $(\alpha,\mu)\in l\mbb{Z}$ at all simple roots $\alpha$.  This lattice $X^\ast$ contains $X^M$, so that we have an inclusion $\Rep \dG\to \Rep G^\ast$.  We then obtain our quantum Frobenius, functor say, by restricting the more expansive quantum Frobenius $\Rep G^\ast\to \Rep G_q$ from \cite{lusztig93} along the inclusion from $\Rep \dG$.  Very directly, our dual group is dictated by the $R$-matrix while Lusztig's dual group is dictated by representation theoretic considerations.
\end{remark}

\subsection{Quantum function algebras}

For us the quantum function algebra $\O(G_q)$ is a formal device which allows us to articulate certain categorical observations in a ring theoretic language.  The Hopf algebra $\O(G_q)$ is the unique Hopf algebra so that we have an equality
\[
\operatorname{Corep}\O(G_q)=\Rep G_q
\]
of non-full monoidal subcategories in $\Vect$.  Via Tannakian reconstruction \cite{delignemilne82,schauenburg92}, one obtains $\O(G_q)$ (uniquely) as the coendomorphism algebra of the forgetful functor
\[
forget:\Rep G_q\to \Vect,
\]
and the Hopf structure on $\O(G_q)$ is derived from the monoidal structure on $forget$.

\subsection{The small quantum group}
\label{se:small}

For $G$ and $q$ as above the small quantum group $\uqG$ is essentially the small quantum group from Lusztig's original work \cite{lusztig90,lusztig90II}, but with some slight variation in the grouplikes.  We provide a presentation of the small quantum group $\uqG$ at odd order $q$, and refer the reader to \cite{arkhipovgaitsgory03,gainutdinovlentnerohrmann} or \cite{negron21} for details on the even order case.
\par

First, let $A$ denote the character group on the quotient $X/X^M$, $A=(X/X^M)^\vee$.  For any root $\gamma$ let $K_\gamma\in A$ denote the character $K_\gamma:X/X^M\to \mbb{C}^\ast$, $K_\gamma(\bar{\lambda})=q^{(\gamma,\lambda)}$.  We now define
\[
u(G_q):=\frac{k\langle E_\alpha,\ F_\alpha,\ \xi:\alpha\text{ simple roots, }\xi\in A\rangle}
{\left(\begin{array}{c}\text{$q$-Serre relations \cite[(a3)--(a5)]{lusztig90II}, relations from $A$,}\vspace{2mm}\\
\xi\cdot E_\alpha\cdot \xi^{-1}=\xi(\alpha)E_\alpha,\ \xi\cdot F_\alpha\cdot \xi^{-1}=\xi(-\alpha)F_\alpha
\end{array}\right)}.
\]
So, a representation of $\uqG$ is just a representations of Lusztig's usual small quantum group which admits an additional grading by $X/X^M$ for which the $K_\alpha$ act as the appropriate semisimple endomorphisms $K_\alpha \cdot v=q^{(\alpha,\deg(v))}v$.  The algebra $\uqG$ admits the expected Hopf structure, just as in \cite{lusztig90II}, with the $\xi\in A$ grouplike and the $E_\alpha$ and $F_\alpha$ skew primitive.
\par

Now, the simple representations for $\uqG$ are labeled by highest weights $L(\bar{\lambda})$, for $\bar{\lambda}\in X/X^M$.  A given simple $L(\bar{\lambda})$ is dimension $1$, and hence invertible with respect to the tensor product, precisely when $\bar{\lambda}$ satisfies $q^{(\alpha,\lambda)}=1$ at all simple roots $\alpha$.  So if we let $X^\ast\subset X$ denote the sublattice of weights $\lambda$ with $(\lambda,\alpha)\in l\mbb{Z}$ for all simple $\alpha$, then we have $X^M\subset X^\ast$ and the subgroup $X^\ast/X^M\subset X/X^M$ labels all $1$-dimensional simple representations.  These simples form a fusion subcategory in $\Rep \uqG$ so that we have a tensor embedding
\[
\Vect(X^\ast/X^M)\to \Rep \uqG,
\]
where $\Vect(X^\ast/X^M)$ denotes the category of $X^\ast/X^M$ graded vector spaces, and we have the corresponding Hopf quotient $\uqG\to \O(X^\ast/X^M)$.

\begin{example}
When $G$ is of adjoint type $X^\ast=X^M$, so that $\Vect(X^\ast/X^M)$ is trivial.  When $G$ is simply-connected $X^\ast=lP$, $X^M=lQ$, and $\Vect(X^\ast/X^M)$ is isomorphic to representations of the center $Z(G)$.
\end{example}

We note that the $R$-matrix for $G_q$ provides a well-defined global operator on products of $\uqG$-representations, so that we have a braiding on $\Rep \uqG$ given by the same formula
\[
c_{V,W}:V\ot W\to W\ot V,\ \ c_{V,W}(v,w)=R^{21}(w\ot v),
\]
where our notation is as in \cite[\S 8.3]{egno15}.  One can see that this braiding on $\Rep \uqG$ is non-degenerate, in the sense that the M\"uger center vanishes, since the induced form $\bar{\Omega}$ on $X/X^M$ is non-degenerate \cite[Theorem 5.3]{negron21}.  We have the restriction functor
\[
\res:\Rep G_q\to \Rep \uqG
\]
which is braided monoidal.  The following result is essentially covered in works of Andersen and coauthors \cite{andersen03,andersenwen92,andersenparadowski95,andersenpolowen91}.

\begin{proposition}\label{prop:andersen}
\begin{enumerate}
\item For any simple object $L$ in $\Rep \uqG$, there is a simple $G_q$-representation $L'$ so that $L$ is a summand of $\res(L')$.
\item Any projective object in $\Rep G_q$ restricts to a projective in $\Rep \uqG$.
\item For any projective object $P$ in $\Rep \uqG$, there is a projective $G_q$-representation $P'$ so that $P$ is a summand of $\res(P')$.
\end{enumerate}
\end{proposition}

\begin{proof}
For $\lambda\in X^+$ let $\msc{L}(\lambda)$ denote the corresponding simple representation in $\Rep G_q$.  Let $A'=(X/X^\ast)^\ast$ and consider $u'\subset \uqG$ the (quasi-)Hopf subalgebra generated by the $E_\alpha$, $F_\alpha$, and elements in $A'$.  We have the exact sequence of (quasi-)Hopf algebras $1\to u'\to \uqG\to \O(X^\ast/X^M)\to 1$.  From the corresponding spectral sequence on cohomology we find that an object in $\Rep \uqG$ is projective if and only if its restriction to $u'$ is projective, and obviously any object with simple restriction to $u'$ is simple over $\uqG$.
\par

The representation category of $u'$ is the category $\msc{C}_{G_1}$ from \cite{andersenparadowski95}.  So by \cite[Theorem 3.12]{andersenparadowski95} we see that all simple representations for $u'$ are restricted from representations over $G_q$, and by the formula \cite[Theorem 1.10]{andersenwen92} one can see that all simples for $\Rep\uqG$ are summands of a simple from $G_q$.  (The representation ``$\bar{L}_k(\lambda_0)$" from \cite{andersenwen92} will split into $1$-dimensional simples over $\uqG$ and there will be no twists ``$\bar{L}_k(\lambda_i)^{(i)}$" since we assume $\operatorname{char}(k)=0$.)  So we obtain (i).  Statements (i) and (ii) follow from the fact that the Steinberg representation is simple and projective over $G_q$ and restricts to a simple and projective representation over $\uqG$ \cite[Proposition 2.2]{andersenwen92} \cite[Theorem 9.8]{andersenpolowen91} \cite[Corollary 9.7]{andersen03}.
\end{proof}

\subsection{A remark on grouplikes}

In the literature there are basically two choices of grouplikes for the small quantum group which are of interest.  In the first case, we take the small quantum group with grouplikes given by characters on the quotient $X/X^\ast$, where $X^\ast=lP\cap X$.  This is a choice which is relevant for many representations theoretic purposes, and which reproduces Lusztig's original small quantum group \cite{lusztig90,lusztig90II} at simply-connected $G$ and odd order $q$.  In the second case, one proceeds as we have here and considers grouplikes given by characters on the quotient $X/X^M$.  This is a choice relevant for physical applications, as one preserves the $R$-matrix and hence allows for the small quantum group to be employed in constructions and analyses of both topological and conformal field theories, see for example \cite{derenzigeerpatureau18,schweigertwoike21,brochierjordansafronovsnyder,creutziggannon17,feigintipunin,gannonnegron}.
\par

This movement of the grouplikes for the small quantum group corresponds precisely to the choice of dual group to $G_q$ for the quantum Frobenius (discussed above).  One has the maximal choice $G^\ast$, or the choice $\dG$ dictated by the $R$-matrix.

\subsection{De-equivariantization and the small quantum group}
\label{sect:de_equiv}

We have the quantum Frobenius $\Fr:\Rep \dG\to \Rep G_q$ as above, and define the de-equivariantization in the standard way
\[
(\Rep G_q)_{\dG}:=\left\{\begin{array}{c}\text{the category of arbitrary}\\
\text{$\O(\dG)$-modules in }\Rep G_q\end{array}\right\}={_{\O(\dG)}\mbf{M}^{\O(G_q)}}
\]
\cite{dgno10}.  Here we abuse notation to write the image of the object $\O(\dG)$ in $\Rep \dG$ under quantum Frobenius simply as $\O(\dG)\in \Rep G_q$.  The category $(\Rep G_q)_{\dG}$ is monoidal under the product $\ot=\ot_{\O(\dG)}$.  We have the de-equivariantization map $\de:\Rep G_q\to (\Rep G_q)_{\dG}$, which is a free module functor $\de(V):=\O(\dG)\ot_k V$, and the category $(\Rep G_q)_{\dG}$ admits a unique braided monoidal structure so that the de-equivariantization map is a braided monoidal functor \cite[Theorem 1.10]{kirillovostrik02}.  This braiding is given directly by the $R$-matrix for $\Rep G_q$,
\[
c_{M,N}:M\ot N\to N\ot M,\ \ m\ot n\mapsto R^{21}(n\ot m).
\]
\par

To elaborate, objects in $(\Rep G_q)_{\dG}$ are unambiguously $\O(\dG)$-bimodules via an application of the braiding, and one employs this bimodule structure in constructing the product on $(\Rep G_q)_{\dG}$.  The braiding is, again, given by the $R$-matrix.

\begin{remark}
Note that objects in the de-equivariantization are naturally bimodules, but are not necessarily \emph{symmetric} bimodules when $q$ is of even order.  This complicates the situation slightly when working with unrestricted $q$ (cf.\ \cite[\S 7.2]{negron21}).
\end{remark}
\par

We consider for any $M$ in $(\Rep G_q)_{\dG}$ the associated sheaf $M^\sim$ for the right action of $\O(\dG)$, and have the linear functor $(-)^{\sim}:(\Rep G_q)_{\dG}\to \QCoh(\dG)$.  This functor is faithful and so identifies the de-equivariantization with a certain non-full subcategory in $\QCoh(\dG)$, which one might refer to as the category of $G_q$-equivariant sheaves over $\dG$.  We let $\QCoh(\dG)^{G_q}$ denote this category of $G_q$-equivariant sheaves on $\dG$, so that we have an equivalence
\begin{equation}\label{eq:669}
(-)^{\sim}:(\Rep G_q)_{\dG} \overset{\sim}\to \QCoh(\dG)^{G_q}.
\end{equation}
The above equivalence induces a braided monoidal structure on $\QCoh(\dG)^{G_q}$ under which the associated sheaf functor becomes a braided monoidal functor.
\par

The induced monoidal structure on $\QCoh(\dG)^{G_q}$ is the usual one at odd order $q$, i.e.\ the one induced by the ambient category of non-equivariant sheaves, and can be described in purely geometric terms in the even order case as well.  We describe this monoidal structure explicitly in Section \ref{sect:q_mon} below.

\begin{definition}
The quantum Frobenius kernel for $G$ at $q$ is the braided monoidal category of $G_q$-equivariant quasi-coherent sheaves over $\dG$,
\[
\FK{G}_q:=\QCoh(\dG)^{G_q}.
\]
\end{definition}

The compact objects in $\FK{G}_q$ are precisely those equivariant sheaves which are coherent over $\dG$ \cite[Lemma 8.4]{negron21}.  As a consequence of Theorem \ref{thm:ag} below, all coherent sheaves are furthermore dualizable.
\par


We have the following result of Arkhipov and Gaitsgory \cite{arkhipovgaitsgory03} and \cite{negron21}.  

\begin{theorem}[{\cite{arkhipovgaitsgory03} \cite[Proposition 7.3]{negron21}}]\label{thm:ag}
Taking the fiber at the identity in $\dG$ provides an equivalence of braided monoidal categories
\[
1^\ast:\FK{G}_q\overset{\sim}\to \Rep \uqG.
\]
\end{theorem}

\begin{remark}
At odd order $q$ Theorem \ref{thm:ag} is alternately deduced from Takeuchi and Kreimer-Takeuchi \cite[Corollary 1.10]{kreimertakeuchi81} \cite[pg.\ 456 \& Theorem 2]{takeuchi79}.
\end{remark}

The above theorem tells us that the dualizable objects in $\FK{G}_q$ are precisely the compact objects, i.e.\ coherent equivariant sheaves, as claimed above.  This subcategory of coherent sheaves is a finite tensor category which is equivalent to the category of finite-dimensional $\uqG$-representations, via the above equivalence.  One can also show that all objects in $\FK{G}_q$ are flat over $\dG$.  At odd order $q$ this follows by \cite[Theorem 5]{takeuchi79}, and at even order $q$ this can be argued from the materials of Section \ref{sect:q_mon}.
\par


From the above geometric perspective the de-equivariantization map for $(\Rep G_q)_{\dG}$ becomes an equivariant vector bundle map $E_-:\Rep G_q\to \FK{G}_q$, $E_V~=~\O_{\dG}~\ot_k~V$, which we still refer to as the de-equivariantization functor.  One sees directly that the equivalence of Theorem \ref{thm:ag} fits into a diagram
\[
\xymatrix{
 & \Rep G_q\ar[dr]^{\res}\ar[dl]_{E_-}\\
\FK{G}_q\ar[rr]^\sim_{1^\ast} & & \Rep \uqG.
}
\]
From this point on (up to Section \ref{sect:GRT}, that is) we essentially forget about the Hopf algebra $\uqG$, and work strictly with its geometric incarnation $\FK{G}_q$.

\subsection{The $\dG$-action on the quantum Frobenius kernel}

As explained in \cite{arkhipovgaitsgory03,dgno10} we have a translation action of $\dG$ on the category $\FK{G}_q=\QCoh(\dG)^{G_q}$ of $G_q$-equivariant sheaves.  This gives an action of $\dG$ on $\FK{G}_q$ by braided tensor automorphisms.  This action is algebraic, in the precise sense of \cite[Appendix A]{negron21}, and we have the corresponding group map $\dG\to \underline{\Aut}_{\ot}^{br}(\FK{G}_q)$.  In terms of this translation action of $\dG$, the de-equivariantization map from the big quantum group restricts to an equivalence
\[
E_-:\Rep G_q\overset{\sim}\longrightarrow (\FK{G}_q)^{\dG}
\]
onto the monoidal category of $\dG$-equivariant objects in $\FK{G}_q$ \cite[Proposition 4.4]{arkhipovgaitsgory03}.
\par

One can translate much of the analysis in this text from the small quantum group to the big quantum group by restricting to $\dG$-equivariant objects in $\FK{G}_q$.  One can compare, for example, with \cite{arkhipovbezrukavnikovginzburg04}.

\section{The quantum Borel and Kempf's vanishing}
\label{sect:borel}

We give a presentation of the (positive) small quantum Borel which is in line with the presentation of Section \ref{sect:quantumgroups} for the small quantum group.  We subsequently provide a spectral sequence relating cohomology for the big and small quantum Borels, and recall a statement of Kempf's vanishing theorem in the quantum context.

\subsection{Quantum Frobenius for the Borel, and de-equivariantization}

As with the (big) quantum group, we let $\Rep B_q$ denote the category of integrable $U_q(\mfk{b})$-representations which are appropriately graded by the character lattice $X$.  The quantum Frobenius for the quantum group induces a quantum Frobenius for the Borel, $\Fr:\Rep \dB\to \Rep B_q$ \cite{lusztig93}.  This quantum Frobenius identifies $\Rep \dB$ with the full subcategory of $B_q$-representations whose $X$-grading is supported on the sublattice $X^M$ (see Section \ref{sect:Gq}).
\par

The functor $\Fr$ is a fully faithful tensor embedding, in the sense that its image is closed under taking subquotients.  This implies that the corresponding Hopf algebra map $fr:\O(\dB)\to \O(B_q)$, which one can obtain directly by Tannakian reconstruction, is an inclusion \cite[Lemma 2.2.13]{schauenburg92}.
\par

We now consider the restriction functor $\Rep B_q\to \Rep \uqB$, where $\uqB$ is the subalgebra in $\uqG$ generated by the grouplikes and the positive root vectors.  As in the case of the full group, $\uqB$ is naturally a quasi-Hopf algebra at even order $q$ and the restriction functor is given a tensor structure via a choice of ``balancing function", as in \cite[Proposition 4.6]{negron21} \cite[Theorem 5.12.7]{egno15}.  We restrict our attention to the algebra structure on $\uqB$ however to obtain a consistent presentation.
\par

By \cite[Lemma 2.2.13]{schauenburg92}, surjectivity of the above restriction functor implies that the associated coalgebra map $\O(B_q)\to (\uqB)^\ast$ is surjective.  Furthermore, an object in $\Rep B_q$ is in the image of quantum Frobenius if and only if that object has trivial restriction to $\uqB$.  Since $\uqB$ is normal in the big quantum Borel \cite[Proposition 35.3.1]{lusztig93}, it follows that $\O(\dB)$ is identified with the $\uqB$-invariants, or $(\uqB)^\ast$-coinvariants, in the quantum function algebra $\O(\dB)=\O(B_q)^{\uqB}$ via the map $fr:\O(\dB)\to \O(B_q)$.
\par

To rephrase what we have just said; we observe an exact sequence of coalgebras
\[
k\to \O(\dB)\overset{fr}\to \O(B_q)\to (\uqB)^\ast\to k
\]
which corresponds to, and is (re)constructed from, the exact sequence of tensor categories
\[
\Vect\to \Rep \dB\overset{\Fr}\longrightarrow \Rep B_q\to \Rep(\uqB)\to \Vect
\]
\cite[Definition 3.7]{bruguieresnatale11} (cf.\ \cite{etingofgelaki17}).

\begin{proposition}\label{prop:129}
The quantum function algebra $\O(B_q)$ is faithfully flat over $\O(\dB)$, and injective over $\uqB$.
\end{proposition}

\begin{proof}
Faithful flatness of $\O(B_q)$ over $\O(\dB)$ follows by Masuoka \cite[Theorem 1.3]{masuoka91}.  The same result \cite[Theorem 1.3]{masuoka91} tells us that $\O(B_q)$ is coflat over $(\uqB)^\ast$.  Equivalently, $\O(B_q)$ is injectivity over $\uqB$.
\end{proof}

We also need the following centrality lemma.

\begin{lemma}\label{lem:B_center}
The $q$-exponentiated, normalized Killing form $\Omega^{-1}$ (see Section \ref{sect:Gq}) provides the quantum Frobenius functor $\operatorname{Fr}:\Rep \dB\to \Rep B_q$ with a central structure
\[
\operatorname{symm}_{V,W}:V\ot \operatorname{Fr}(W)\to \operatorname{Fr}(W)\ot V,\ \ v\ot w\mapsto \Omega^{-1}(w\ot v).
\]
\end{lemma}

\begin{proof}
The fact that $\Omega$ is defined by a bilinear form $q^{(-,-)}$ implies that the operations $\operatorname{symm}_{V,W}$ are associative, relative to the actions of $\Rep\dB$ and $\Rep B_q$ on the left and right.  The fact that this form vanishes on the dual lattice $X^M\times X^M$ furthermore implies that $\operatorname{symm}_{V,W}$ recovers the usual (trivial) symmetry on $\Rep\dB$ whenever $V$ is in the image of quantum Frobenius, and hence that $\operatorname{symm}_{V,W}$ is compatible with the original braiding on $\Rep \dB$.  All that is left to show is that $\operatorname{symm}_{V,W}$ is a map of $B_q$-representations at all $V$ and $W$.  Equivalently, if we consider the Hopf embedding $fr:\O(\dB)\to \O(B_q)$ defined by quantum Frobenius, we must show that the diagram
\begin{equation}\label{eq:765}
\xymatrix{
\O(B_q)\ot \O(\dB)\ar[r]^{\operatorname{Ad}_{\Omega^{-1}}}\ar[dr]_{\operatorname{swap}} & \O(B_q)\ot \O(\dB)\ar[r] & \O(B_q)\\
 & \O(\dB)\ot\O(B_q)\ar[ur]
}
\end{equation}
commutes (cf.\ \cite[Definition 2.4.4, Theorem 2.4.5]{schauenburg92II}).  Here $\operatorname{swap}$ denotes the trivial vector space symmetry and the two maps to $\O(B_q)$ are given by multiplication.
\par

However, we have the corresponding diagram for the quantum group
\begin{equation}\label{eq:772}
\xymatrix{
\O(G_q)\ot \O(\dG)\ar[r]^{\operatorname{Ad}_R}\ar[dr]_{\operatorname{swap}} & \O(G_q)\ot \O(\dG)\ar[r] & \O(G_q)\\
 & \O(\dG)\ot\O(G_q)\ar[ur] &
}
\end{equation}
\cite[Theorem 2.4.5]{schauenburg92II} and note that the element $R$ acts as $\Omega^{-1}$ when applied to the product $\O(G_q)\ot \O(\dG)$, so that $\operatorname{Ad}_R$ reduces to $\operatorname{Ad}_{\Omega^{-1}}$ in the above diagram.  We then obtain \eqref{eq:765} from \eqref{eq:772} by applying the Hopf surjections $\O(\dG)\to \O(\dB)$ and $\O(G_q)\to \O(B_q)$.
\end{proof}

\begin{remark}
When $q$ is of odd order the form appearing in the operation $\operatorname{symm}_{V,W}$ is identically $1$, and the symmetry of Lemma \ref{lem:B_center} collapses to the standard vector space symmetry.
\end{remark}

As with the quantum group, we consider the monoidal category $(\Rep B_q)_{\dB}={_{\O(\dB)}\mathbf{M}^{\O(B_q)}}$ of $\O(\dB)$-modules in $\Rep B_q$ \cite[Theorem 2.5]{pareigis95}, and define the quantum Frobenius kernel $\FK{B}_q$ for the quantum Borel as the corresponding (monoidal) category of $B_q$-equivariant sheaves on $\dB$,
\[
\FK{B}_q:=\QCoh(\dB)^{B_q}\cong (\Rep B_q)_{\dB}.
\]
The following is an application of Proposition \ref{prop:129} and \cite[Theorem 1]{takeuchi79}.

\begin{corollary}\label{cor:438}
Taking the fiber at the identity $1:\Spec(k)\to \dB$ provides an equivalence of monoidal categories
\[
1^\ast:\FK{B}_q\overset{\sim}\to \Rep\uqB.
\]
\end{corollary}

\begin{proof}
The fact that $1^\ast$ is a linear equivalence follows from the aforementioned results, and the usual monoidal structure on the pullback functor endows $1^\ast$ with a monoidal structure when $q$ is of odd order.  At even order $q$ monoidality of $1^\ast$ is established as in \cite[Proposition 7.3]{negron21}.
\end{proof}

As with the quantum Frobenius kernel, the compact/dualizable objects in $\FK{B}_q$ are precisely those equivariant sheaves which are coherent over $\dB$, and all objects in $\FK{B}_q$ are flat over $\dB$ (cf.\ \cite[Corollary 1.5]{masuoka91}).

\subsection{A spectral sequences for $B_q$-extensions}

\begin{lemma}[{\cite{doi81}}]\label{lem:doi}
An object $V$ in $\Rep B_q$ is injective if and only if $V$ is a summand of some additive power $\oplus_{i\in I} \O(B_q)$.
\end{lemma}

\begin{proof}
Comultiplication provides an injective comodule map $V\to \underline{V}\ot \O(B_q)$, where $\underline{V}$ is the vector space associated to the representation $V$.  Since any cofree comodule is injective \cite{doi81}, this inclusion is split.
\end{proof}

Since the $\uqB$-invariants in $\O(B_q)$ are precisely the classical algebra $\O(\dB)$, we observe the following.

\begin{corollary}\label{cor:143}
If $W$ is injective over $B_q$, and $V$ is a finite-dimensional, then $\Hom_{\uqB}(V,W)$ is an injective $\dB$-representation.
\end{corollary}

\begin{proof}
We have $\Hom_{\uqB}(V,W)=\Hom_{\uqB}(k,{^\ast V}\ot W)$, and ${^\ast V}\ot W$ is injective over $B_q$ in this case.  So it suffices to assume $V=k$, in which case the result follows by Lemma \ref{lem:doi} and the calculation $\O(\dB)=\O(B_q)^{\uqB}$.
\end{proof}

\begin{proposition}\label{prop:151}
Let $V$ and $W$ be in $\Rep(B_q)$, and assume that $V$ is finite-dimensional.  There is a natural isomorphism
\[
\RHom_B(k,\RHom_{\uqB}(V,W))\cong \RHom_{B_q}(V,W),
\]
and subsequent spectral sequence
\[
E_2^{i,j}=\Ext^i_B(k,\Ext^j_{\uqB}(V,W))\ \Rightarrow\ \Ext^{i+j}_{B_q}(V,W).
\]
\end{proposition}

\begin{proof}
The quantum function algebra $\O(B_q)$ is an injective $\uqB$-module, by Proposition \ref{prop:129}.  It follows by Lemma \ref{lem:doi} that any injective $B_q$-representation restricts to an injective $\uqB$-representation.  So the result follows by Corollary \ref{cor:143}.
\end{proof}

\subsection{Kempf vanishing and a transfer theorem}

We recall some essential relations between quantum group representations and representations for the quantum Borel.  The following vanishing result, which first appears in works of Andersen, Polo, and Wen \cite{andersenpolowen91,andersenwen92} with some restrictions on the order of $q$, appears in complete generality in works of Woodock and Ryom-Hasen \cite[Theorem 8.7]{woodock97} \cite[Lemma 4.3, Theorem 5.5]{ryom03}.

\begin{theorem}\label{thm:kempf_vanish}
Let $\mathsf{I}^0$ denote induction from the quantum Borel, $\mathsf{I}^0:\Rep B_q\to \Rep G_q$.
\begin{enumerate}
\item $\mathsf{I}^0(\1)=\1$.
\item The higher derived functors $\mathsf{I}^{>0}(\1)$ vanish.
\end{enumerate}
\end{theorem}

We can now employ the information of Theorem~\ref{thm:kempf_vanish} and follow exactly the proof of~\cite[Theorem 2.1]{clineparshallscottkallen77} to observe the following transfer theorem.

\begin{theorem}[\cite{clineparshallscottkallen77}]\label{thm:cpsk}
For arbitrary $V$ and $W$ in $\Rep G_q$, and $i\geq 0$, the restriction functor $\Rep G_q\to \Rep B_q$ induces an isomorphism on cohomology
\[
\Ext^i_{G_q}(V,W)\overset{\cong}\longrightarrow \Ext^i_{B_q}(V,W).
\]
\end{theorem}

\section{A complete family of small quantum Borels}
\label{sect:borels}

We extend the construction of the positive small quantum Borel, as in Corollary \ref{cor:438} above, to provide a family of small quantum Borels which are parametrized by points in the flag variety $\dG/\dB$.  As far as we understand this construction, and even the assertion that there are small quantum Borels associated to arbitrary Borel subgroups in $\dG$, is new.

\subsection{Small quantum Borels}

For any $k$-point $\lambda:\Spec(k)\to \dG/\dB$ we have the corresponding $\dB$-coset $\iota_\lambda:\dB_\lambda\to \dG$, which is the fiber of $\lambda$ along the quotient map $\pi:\dG\to \dG/\dB$.  This coset is a $\dB$-torsor under the right translation action of $\dB$ and we have the associated central algebra object $\O(\dB_\lambda)$ in $\Rep \dB\subset \Rep B_q$.  We consider the monoidal category $(\Rep B_q)_{\dB_\lambda}$ of $\O(\dB_\lambda)$-modules in $\Rep B_q$ and subsequent category
\[
\msc{B}_\lambda:=\QCoh(\dB_\lambda)^{B_q}\cong (\Rep B_q)_{\dB_\lambda}
\]
of $B_q$-equivariant sheaves on $\dB_\lambda$ \cite[Theorem 2.5]{pareigis95}.
\par

Pulling back along the inclusion $\iota_\lambda:\dB_\lambda\to \dG$ provides a central monoidal functor
\[
\res_\lambda:=\iota_\lambda^\ast:\FK{G}_q\to \msc{B}_\lambda
\]
with central structure induced by the $R$-matrix on global sections
\[
\gamma_{M,N}:M\ot \res_\lambda(N)\to \res_\lambda(N)\ot M,\ \ \gamma_{M,N}(m\ot n)=R^{21}(n\ot m).
\]
Algebraically, $\res_\lambda$ restricts from $\O(\dG)$-modules in $\Rep G_q$ to modules in $\Rep B_q$, then applies base change $\O(\dB_\lambda)\ot_{\O(\dG)}-$.  We have $\msc{B}_1=\FK{B}_q$ and the functor $\res_1:\FK{G}_q\to \FK{B}_q$ is identified with the standard restriction functor for the small quantum group, in the sense that the diagram
\begin{equation}\label{eq:849}
\xymatrix{
\FK{G}_q\ar[rr]^{\res_1}\ar[d]_{1^\ast}^\sim & & \FK{B}_q\ar[d]_{1^\ast}^\sim\\
\Rep \uqG\ar[rr]^{\res} & & \Rep \uqB
}
\end{equation}
commutes.
\par

At a general closed point $\lambda$, any choice of a point $x:\Spec(k)\to \dB_\lambda$ provides a $\dB$-equivariant isomorphism $x\cdot-:\dB\to \dB_\lambda$.  Taking global sections then provides an isomorphism $x:\O(\dB_\lambda)\to \O(\dB)$ of algebra objects in $\Rep \dB$.  So we see that pushing forward along $x$ provides an equivalence of tensor categories $x:\msc{B}_1\to \msc{B}_\lambda$ which fits into a diagram
\begin{equation}\label{eq:502}
\xymatrix{
\FK{G}_q\ar[rr]^x_\sim\ar[d]_{\res} & & \FK{G}_q\ar[d]_{\res_\lambda}\\
\msc{B}_1\ar[rr]^x_\sim & & \msc{B}_\lambda.
}
\end{equation}

Now, let us consider an arbitrary \emph{geometric} point $\lambda:\Spec(K)\to \dG/\dB$.  At $\lambda$ we again have the fiber $\iota_\lambda: \dB_\lambda\to G$, which now has the structure of a $K$-scheme, and which is a torsor over $\dB_K$.  We consider the monoidal category
\[
\msc{B}_\lambda:=\QCoh(\dB_\lambda)^{(B_K)_q}
\]
of equivariant sheaves relative to the base change $(B_K)_q$.  Pulling back along $\iota_\lambda$ again provides a central monoidal functor
\[
\res_\lambda:=\iota_\lambda^\ast:\FK{G}_q\to \msc{B}_\lambda
\]
which factors as a base change map composed with restriction along $\iota_{K,\lambda}$
\[
\res_\lambda=\big(\FK{G}_q\overset{(-)_K}\longrightarrow \FK(G_K)_q\overset{\res_{K,\lambda}}\longrightarrow \msc{B}_\lambda\big).
\]
All of this is to say that, after base change, the construction of $\msc{B}_\lambda$ at a geometric point for $\dG/\dB$ is no different from the construction at a closed point.

\begin{definition}\label{def:B_lambda}
For a given geometric point $\lambda:\Spec(K)\to \dG/\dB$, the $\FK{G}_q$-central, monoidal category $\msc{B}_\lambda$ is called the small quantum Borel at $\lambda$.
\end{definition}

The following Proposition is deduced from the equivalence $x:(\msc{B}_1)_K\overset{\sim}\to \msc{B}_\lambda$ provided by any choice of $K$-point $x:\Spec(K)\to \dB_\lambda$.

\begin{proposition}
\label{prop:proj}
At each geometric point $\lambda:\Spec(K)\to \dG/\dB$ the monoidal category $\msc{B}_\lambda$ has the following properties:
\begin{itemize}
\item $\msc{B}_\lambda$ has enough projectives and injectives, and an object is projective if and only if it is injective (\ \cite{faithwalker67}).
\item $\msc{B}_\lambda$ admits a compact projective generator.
\item The compact objects in $\msc{B}_\lambda$ are precisely those $(B_K)_q$-equivariant sheaves which are coherent over $\dB_\lambda$, and all compact objects are dualizable.
\item Coherent sheaves in $\msc{B}_\lambda$ form a finite tensor subcategory which is of (Frobenius-Perron) dimension $\dim \uqB$.
\item All objects in $\msc{B}_\lambda$ are flat over $\dB_\lambda$.
\item The central tensor functor $\res_{K,\lambda}:\FK(G_K)_q\to \msc{B}_\lambda$ is surjective.
\item $\res_{K,\lambda}:\FK(G_K)_q\to \msc{B}_\lambda$ sends projectives to projectives (\ \cite[Theorem 2.5]{etingofostrik04}).
\end{itemize}
\end{proposition}

\begin{remark}
The analogous construction in finite characteristic $\res_\lambda:\Rep\mathbf{G}_{(1)}\cong\FK\mathbf{G}\to \msc{B}_\lambda$ is explicitly identified, via Tannakian reconstruction \cite[Corollary 2.9, Theorem 3.2]{delignemilne82}, with the closed embedding of the Frobenius kernel $\mathbf{B}^\lambda_{(1)}$ in $\mathbf{G}_{(1)}$, where $\mathbf{B}^\lambda\subset \mathbf{G}_K$ is the Borel subgroup corresponding to $\lambda:\Spec(K)\to \mathbf{G}^{(1)}/\mathbf{B}^{(1)}$.
\end{remark}

\begin{remark}
For an explicit indication that the family $\{\msc{B}_\lambda:\lambda\text{ a geom point for }\dG/\dB\}$ is ``complete" one can see \cite[Theorem 13.1]{negronpevtsova5}.
\end{remark}

\subsection{A notational comment}

As mentioned in Section \ref{sect:quantumgroups}, the small quantum group is almost never referenced in its linear form $\uqG$ in this text.  We will, however, make extensive use of the algebra $\uqB$ throughout.  Furthermore, the algebra $\uqB$ will often appear as a subscript in formulas.  For this reason \emph{we adopt the notation
\[
u:=\uqB
\]
globally throughout this document}.  An unlabeled algebra $u$ which appears anywhere in the text is always the small quantum enveloping algebra $\uqB$ for the positive Borel.

\section{The half-quantum flag variety}
\label{sect:G/Bq}

In this section we introduce the category of sheaves $\QCoh(\Xq)$ on the  
\emph{half-quantum flag variety}.  In light of the construction of the small quantum Borels as equivariant sheaves on the cosets $\dB_\lambda$, one might think of the category $\QCoh(\Xq)$ as a universal small quantum Borel for $G$ at $q$ (see also Section \ref{sect:fibers_G/B} below).
\par

Just as one considers each $\msc{B}_\lambda$ as a $K$-linear monoidal category, one should consider the category of sheaves on the half-quantum flag variety as a $\dG/\dB$-linear monoidal category.  Such linearity can be expressed via an action of the category of sheaves over the classical flag variety on $\QCoh(\Xq)$.

\subsection{The half-quantum flag variety}
\label{sect:G/Bq_intro}

We consider the 
noncommutative space $\Xq:=\dG/B_q$, where the quotient here is interpreted as a stack quotient.  Though this space may not be well-defined, we have a clear presentation of its category of sheaves.

\begin{definition}
The category of quasi-coherent sheaves for the half-quantum flag variety is the abelian monoidal category
\begin{equation}\label{eq:938}
\QCoh(\Xq):=\QCoh(\dG)^{B_q}\cong \{\text{Arbitrary $\O(\dG)$-modules in }\Rep B_q\}.
\end{equation}
We let $\Coh(\Xq)$ denote the full monoidal subcategory of sheaves which are coherent over $\dG$.
\end{definition}

To recall, quantum Frobenius provides a sequence of central tensor functors $\Rep\dG\to \Rep\dB\to \Rep B_q$ through which we view $\O(\dG)$ as a central algebra object in $\Rep B_q$, and consider the corresponding monoidal category of $\O(\dG)$-modules with product $\ot=\ot_{\O(\dG)}$ \cite[Theorem 2.5]{pareigis95} (see also \cite[Theorem 1.5]{kirillovostrik02}).  This is the module category considered at \eqref{eq:938}.

\begin{definition}
An object $M$ in $\QCoh(\Xq)$ is called flat if both operations $-\ot M$ and $M\ot-$ are exact.
\end{definition}

As with the usual flag variety \cite[I.5.8]{jantzen03} \cite[3.3]{suslinfriedlanderbendel97b}, we have an equivariant vector bundle functor
\[
E_-:\Rep B_q\to \QCoh(\Xq),\ \ E_V:=\O_{\dG}\ot_k V.
\]
This functor is exact and monoidal, and has a (non-monoidal) right adjoint provided by the forgetful functor
\[
-|_{B_q}:\QCoh(\Xq)\to \Rep B_q,
\]
which is defined explicitly by applying global sections $M|_{B_q}=\Gamma(\dG,M)$ and forgetting the $\O(\dG)$-action (cf.\ Section \ref{sect:descent}).

\begin{lemma}
The category $\QCoh(\Xq)$ is complete, in the sense that it has all set indexed limits, and has enough injectives.
\end{lemma}

\begin{proof}
The point is that $\QCoh(\Xq)$ is a Grothendieck abelian category.  That is to say, $\QCoh(\Xq)$ is cocomplete, has exact filtered colimits, and admits a generator.  All Grothendieck categories are complete and have enough injectives \cite[Th\'eor\`eme 1.10.1]{grothendieck57}.
\par

The only controversial issue here is the existence of a generator.  However, since $\rep B_q$ is essentially small we can choose a set of representations $\{V_i\}_{i\in I}$ so that each object in $\rep B_q$ admits a non-zero map from some $V_i$.  The object $\oplus_{i\in I}E_{V_i}$ is then a generator for $\QCoh(\Xq)$, where $E_V$ is the vector bundle associated to a given $B_q$-representation $V$.
\end{proof}

\begin{remark}
Consider any lift $W$ of a projective generator for $\Rep u(B_q)$ to $\Rep B_q$.  The category $\QCoh(\Xq)$ is more specifically generated by the orbit of the associated vector bundle $\oplus_{n\geq 0}\mcl{L}_{\rho}^{\ot n}\ot E_{W}$ under the action of the ample line bundle $\mcl{L}_{-\rho}$ for the flag variety \cite[Proposition 4.4]{jantzen03}.
\end{remark}

\subsection{The monoidal structure on $\QCoh(\Xq)$, in geometric terms}
\label{sect:q_mon}

We take a moment to explain how one understands $\QCoh(\Xq)$ explicitly as a monoidal category of \emph{sheaves} on $\dG$.

At even order $q$ there is a left/right distinction for $\O(\dG)$-modules in $\Rep B_q$, and we explicitly realize such modules as sheaves on $\dG$ via the right action.  This gives us our category $\QCoh(\Xq)=\QCoh(\dG)^{B_q}$.  Though $\QCoh(\Xq)$ is not classical, i.e.\ sheaves on an actual stack, its monoidal structure is describable in purely classical terms.
\par

Let $T_q\subset B_q$ denote the quantum torus, which just accounts for the $X$-gradings on objects in $\Rep B_q$.  Note that $T_q$ is an actual group scheme $T_q=\Spec(kX)$, and that our informal translation action of $B_q$ on $\dG$ restricts to an actual translation action of $T_q$ on $\dG$.  We furthermore have the subgroup $A\subset T_q$, $\O(A)=k(X/X^M)$, and restrict this $T_q$-action to obtain an $A$-action on $\dG$.  So we have the sequence of forgetful functors
\[
\QCoh(\Xq)\to \QCoh(\dG)^{T_q}\to \QCoh(\dG)^{A}.
\]
We note that the action of $A$ here is actually trivial, since the $X$-grading on $\O(\dG)$ is concentrated on the subgroup $X^M\subset X$.  Hence $A$-equivariant sheaves on $\dG$ are just sheaves with a grading by the character group $A^\vee=X/X^M$.
\par

For any $\mu\in X$ we consider the algebra automorphism $\mu^{\#}:\O(\dG)\to \O(\dG)$ which is defined on $X$-homogeneous elements as $\mu^{\#}(x)=q^{(\mu,\deg(x))}x$.  These automorphism vanish for $\mu\in X^M$ and hence provide a group map
\begin{equation}\label{eq:1068}
\sigma:(A^\vee)^{op}\to \Aut_{\operatorname{Sch}/k}(\dG)
\end{equation}
from the characters $A^\vee$.
\par

Now, although the objects in $\QCoh(\dG)^A$ are classical, we can introduce some quantum aspects by twisting the tensor product via the automorphisms \eqref{eq:1068}.  In particular, we decompose any object in $\QCoh(\dG)^{A}$ into character sheaves $\msc{F}=\oplus_{\bar{\mu}}\msc{F}_{\bar{\mu}}$, and we define a new product
\begin{equation}\label{eq:1086}
\msc{G}\ot^{\operatorname{new}}\msc{F}:=\bigoplus_{\bar{\mu}\in A^\vee}\sigma(\bar{\mu})^\ast\msc{G}\ot_{\O_{\dG}} \msc{F}_{\bar{\mu}}.
\end{equation}
So we obtain a natural monoidal structure on $\QCoh(\dG)^A$ under which it is interpreted as a monoidal category of quantum equivariant sheaves.  Under this quantum tensor product the forgetful functor
\[
F:\QCoh(\Xq)\to \QCoh(\dG)^A
\]
is naturally monoidal, with structure maps $F(M)\ot^{\operatorname{new}} F(N)=F(M\ot N)$ given by the identity.  One checks this directly by examining how the form $\Omega^{-1}$ translates right $\O(\dG)$-modules to left $\O(\dG)$-modules in $\Rep B_q$.

\begin{proposition}\label{prop:1093}
An object in $\QCoh(\Xq)$ is flat whenever its image in $\QCoh(\dG)$, under the forgetful functor, is flat.
\end{proposition}

\begin{proof}
Since the product on $\QCoh(\Xq)$ is given by the product \eqref{eq:1086} on $\QCoh(\dG)^A$, it suffices to show that objects in $\QCoh(\dG)^A$ are flat for the new product whenever they are flat for the classical product $\ot_{\O_{\dG}}$.  However, this is clear since a sum of sheaves on $\dG$ is flat if and only if each summand is flat, and since the automorphisms $\sigma(\bar{\mu})^\ast$ preserve flatness.
\end{proof}

We note that every object in $\QCoh(\Xq)$ admits a surjection $E\to M$ from a sum of equivariant vector bundles, and such bundles have flat image in $\QCoh(\dG)$.  So $\QCoh(\Xq)$ has enough flat sheaves.

\begin{corollary}
Every object in $\QCoh(\Xq)$ admits a bounded resolution by flat sheaves.
\end{corollary}

\begin{proof}
Follows from Proposition \ref{prop:1093} and the fact that $\dG$ has finite flat dimension.
\end{proof}

Via a completely similar analysis one sees that the quantum product on $\QCoh(\dG)^A$ also realizes the monoidal structure on $\FK{G}_q=\QCoh(\dG)^{G_q}$ in geometric terms.  The corresponding result for the small quantum Borel $\msc{B}_\lambda=\QCoh(\dB_\lambda)^{B_q}$ is obtained by replacing $A$-equivariant sheaves on $\dG$ with the category $\QCoh(\dB_\lambda)^A$ of $A$-equivariant sheaves on the $\dB_\lambda$, endowed with a tensor structure as in \eqref{eq:1086}.

\begin{remark}
The above presentation is equally valid at odd order $q$.  The group map \eqref{eq:1068} just happens to be trivial in this case.
\end{remark}

\subsection{$\QCoh(\Xq)$ as a flat family of categories over $\dG/\dB$}

The quantum Frobenius maps for $G_q$ and $B_q$ fit into a diagram of central tensor functors
\[
\xymatrix{
\Rep \dG\ar[rr]^{\Fr}\ar[d]_{\res} & & \Rep G_q\ar[d]^{\res}\\
\Rep \dB\ar[rr]^{\Fr} & & \Rep B_q.
}
\]
This diagram then implies the existence of a fully faithful monoidal embedding
\begin{equation}\label{eq:zeta}
\QCoh(\dG/\dB)\overset{\pi^\ast}\to \QCoh(\dG)^{\dB}\to \QCoh(\dG)^{B_q}=\QCoh(\Xq).
\end{equation}
We let $\zeta^\ast$ denote this embedding.  We note that $\zeta^\ast:\QCoh(\dG/\dB)\to \QCoh(\Xq)$ inherits a natural central structure which is defined on global sections by the $q$-exponentiated Killing form
\begin{equation}\label{eq:symm}
\operatorname{symm}_{M,\msc{F}}:M\ot \zeta^\ast(\msc{F})\to \zeta^\ast(\msc{F})\ot M,\ \ \operatorname{symm}_{M,\msc{F}}(m\ot s)=\Omega^{-1}(s\ot m).
\end{equation}
In accordance with the language of Section \ref{sect:flatfamily} we have the following.

\begin{proposition}\label{prop:flatfam}
The central tensor functor $\zeta^\ast$ gives $\QCoh(\Xq)$ the structure of a flat family of tensor categories over $\dG/\dB$.
\end{proposition}

\begin{proof}
The functor $\zeta^\ast$ is exact, as its a composite of exact functors, and the central structure is given above.  We must show that $\QCoh(\Xq)$ is generated by a subcategory of compact rigid objects.
\par

We have that $\QCoh(\Xq)$ is generated by the equivariant vector bundles $E_V$, with $V$ a finite-dimensional $B_q$-representation.  Each such vector bundle $E_V$ is dualizable as it is the image of a dualizable object under the monoidal functor $E_-:\Rep B_q\to \QCoh(\Xq)$ \cite[Exercise 2.10.6]{egno15}.  These $E_V$ are additionally compact as the generating representations $V$ are compact in $\Rep B_q$, and the right adjoint $-|_{B_q}$ to the equivariant vector bundle functor commutes with colimits.  So $\QCoh(\Xq)$ is in fact compactly-rigidly generated.
\end{proof}

For $\msc{F}$ in $\QCoh(\dG/\dB)$ we let
\[
\msc{F}\star-:\Coh(\Xq)\to \Coh(\Xq)
\]
denote the corresponding action map, $\msc{F}\star-=\zeta^\ast(\msc{F})\ot-$.  The operation $\msc{F}\star-$ is exact whenever $\msc{F}$ is flat over $\dG/\dB$, and $-\star M$ is exat whenever $M$ is flat over $\Xq$.
\par

Although $\zeta^\ast$ is not precisely the pullback equivalence $\pi^\ast:\QCoh(\dG/\dB)\overset{\sim}\to \QCoh(\dG)^{\dB}$, due to the appearance of quantum Frobenius, we will sometimes abuse notation and write simply $\pi^\ast(\msc{F})$ for the object $\zeta^\ast(\msc{F})$.

\subsection{Sheafy morphisms over $\dG/\dB$}

We have just seen that $\QCoh(\Xq)$ admits a natural module category structure over $\QCoh(\dG/\dB)$, and so becomes a flat family of tensor categories over the flag variety.  Inner morphisms for this module category/sheaf structure provide a sheaf-Hom functor for $\QCoh(\Xq)$. 

\begin{lemma}\label{lem:994}
For any $M$ in $\QCoh(\Xq)$, the operation $-\star M:\QCoh(\dG/\dB)\to \QCoh(\Xq)$ has a right adjoint $\sHom_{\Xq}(M,-)$.  The functors $\sHom_{\Xq}(M,-)$ are furthermore natural in $M$, so that we have a bifunctor
\[
\sHom_{\Xq}:\QCoh(\Xq)^{op}\times \QCoh(\Xq)\to \QCoh(\Xq).
\]
\end{lemma}

\begin{proof}
The functor $-\star M$ commutes with colimits and thus admits a right adjoint.  Naturality in $M$ follows by Yoneda's lemma.
\end{proof}

We note that the functor $\sHom_{\Xq}$ is left exact in both coordinates, since the functor
\[
\Hom_{\dG/\dB}(\msc{F},\sHom_{\Xq}(-,-))\cong \Hom_{\dG/\dB}(\msc{F}\star-,-)
\]
is left exact in each coordinate at arbitrary $\msc{F}$.
\par

In Sections \ref{sect:shom1}--\ref{sect:shom2} we observe that, via general nonsense with adjunctions (cf.\ \cite{ostrik03,etingofostrik04}), the bifunctor $\sHom_{\Xq}$ admits natural composition and monoidal structure maps.  These structure maps localize the monoidal structure on $\QCoh(\Xq)$, in the sense that they recover the composition and tensor structure maps for $\QCoh(\Xq)$ after taking global sections.  In the language of Section \ref{sect:enrich}, we are claiming specifically that the pairing
\begin{equation}\label{eq:991}
(\!\ \operatorname{obj}\QCoh(\Xq),\ \sHom_{\Xq}\!\ )
\end{equation}
provides a monoidal enhancement of $\QCoh(\Xq)$ in the category of sheaves over the flag variety.
\par

Before delving further into these issues, we explain a fundamental relationship between the half-quantum flag variety and the small quantum group.  This relationship motivates our consideration of the aforementioned enhancement, and our consideration of the half-quantum flag variety more generally.

\section{The universal restriction functor $\kappa^\ast:\FK{G}_q\to \QCoh(\Xq)$}
\label{sect:Kempf}

Let us consider again the quantum Frobenius kernel $\FK{G}_q$.  We have the obvious forgetful functor
\begin{equation}\label{eq:funcTOR}
\FK{G}_q=\QCoh(\dG)^{G_q}\to \QCoh(\dG)^{B_q}=\QCoh(\Xq)
\end{equation}
This functor is immediately seen to be exact and monoidal, and the $R$-matrix for the quantum group provides it with a central structure (see Section \ref{sect:Kempf_Z} below).  We refer to the above functor as the universal restriction functor, and denote it
\[
\kappa^\ast:\FK{G}_q\to \QCoh(\Xq).
\]
Our language is informed by certain observations which we articulate in Section \ref{sect:fibers_G/B}.
\par

Below we prove that the universal restriction functor is fully faithful, and induces a fully faithful functor on unbounded derived categories as well.  In the statement of the following theorem, $D(\FK{G}_q)$ denotes the unbounded derived category of complexes in $\FK{G}_q$, and $D(\Xq)$ is the unbounded derived category of complexes in $\QCoh(X_q)$.

\begin{theorem}\label{thm:Kempf}
The central tensor functor $\kappa^\ast:\FK{G}_q\to \QCoh(\Xq)$ is fully faithful, and induces a fully faithful monoidal embedding
\[
\kappa^\ast(=\operatorname{L}\kappa^\ast):D(\FK{G}_q)\to D(\Xq)
\]
on the corresponding unbounded derived categories.
\end{theorem}

\begin{proof}
Follows immediately by Theorem \ref{thm:kempf_ext} below.
\end{proof}

Theorem \ref{thm:Kempf} is essentially a consequence of Kempf vanishing \cite{kempf76}, or rather the quantum analog of Kempf vanishing provided in \cite{andersenpolowen91,woodock97,ryom03}.

\subsection{Compact sheaves in $D(\Xq)$}

We begin with a little lemma, the proof of which is deferred to Section \ref{sect:relative_p}.

\begin{lemma}\label{lem:1090}
Suppose that $V$ is a finite-dimensional $B_q$-representation and that the restriction of $V$ to the small quantum Borel $\uqB$ is projective.  Then the associated vector bundle $E_V$ over $\Xq$ is compact in $D(\Xq)$.  Similarly, any summand of $E_V$ is compact in $D(\Xq)$.
\end{lemma}

Now, by Proposition \ref{prop:andersen}, all coherent projectives in $\FK{G}_q$ are summands of vector bundles $E_V$ with $V$ projective over $G_q$, and hence projective over $\uqG$ and $\uqB$ as well (see Proposition~\ref{prop:proj}).  So Lemma \ref{lem:1090} implies the following.

\begin{corollary}\label{cor:1110}
The functor $\kappa^\ast:D(\FK{G}_q)\to D(\Xq)$ sends compact objects in $D(\FK{G}_q)$ to compact objects in $D(\Xq)$.
\end{corollary}

\subsection{Analysis of universal restriction via extensions}

\begin{theorem}\label{thm:kempf_ext}
The forgetful $\kappa^\ast$ induces an isomorphism on cohomology
\begin{equation}\label{eq:1093}
\Ext^i_{\FK{G}_q}(M,N)\overset{\cong}\longrightarrow \Ext^i_{\Xq}(M,N),
\end{equation}
for all $i$ and all $M$ and $N$ in $D(\FK{G}_q)$.
\end{theorem}

\begin{proof}
The forgetful functor is exact and hence induces a map on unbounded derived categories.  We first claim that the map on extensions is an isomorphism whenever $M$ is in the image of the de-equivariantization/equivariant vector bundle map $E_-:\Rep G_q\to \FK{G}_q$, and $N$ is arbitrary in $\FK{G}_q$.  Recall that this vector bundle functor has an exact right adjoint $-|_{G_q}:\FK{G}_q\to \Rep G_q$ provided by taking global sections and forgetting the $\O(\dG)$-action, and we have the analogous adjunction for $\QCoh(\Xq)$, as discussed in Section \ref{sect:G/Bq_intro}.  (We make no notational distinction between vector bundles in $\FK{G}_q$ and in $\QCoh(\Xq)$, and rely on the context to distinguish the two classes of sheaves.)
\par

Since the equivariant vector bundle functors for $\FK{G}_q$ and $\QCoh(\Xq)$ are exact, the adjoints $-|_{G_q}$ and $-|_{B_q}$ preserve injectives.  So we have isomorphisms
\[
\Ext^\ast_{\FK{G}_q}(E_V,N)\overset{\cong}\to \Ext^\ast_{G_q}(V,N|_{G_q})\ \ \text{and}\ \ \Ext^\ast_{\Xq}(E_V,N)\overset{\cong}\to \Ext^\ast_{B_q}(V,N|_{B_q})
\]
which fit into a diagram
\[
\xymatrix{
\Ext^\ast_{\FK{G}_q}(E_V,N)\ar[rr]^{\kappa^\ast}\ar[d]_\cong & & \Ext^\ast_{\Xq}(E_V,N)\ar[d]^\cong\\
\Ext^\ast_{G_q}(V,N|_{G_q})\ar[rr]^{\operatorname{restrict}}_{\cong} & & \Ext^\ast_{B_q}(V,N|_{B_q}).
}
\]
The bottom morphism here is an isomorphism by Kempf vanishing, or more precisely by the transfer result of Theorem \ref{thm:cpsk}.  So we conclude that the top map is an isomorphism.
\par

We now understand that the map \eqref{eq:1093} is an isomorphism whenever $M=E_V$ for some $V$ in $\Rep G_q$ and $N$ is arbitrary in $\FK{G}_q$.  It follows that \eqref{eq:1093} is an isomorphism whenever $M$ is a summand of a vector bundle $E_V$ and $N$ is arbitrary in $\FK{G}_q$, and hence whenever $M$ is simple or projective in $\FK{G}_q$ by Proposition \ref{prop:andersen}.
\par

Since coherent projectives in $\FK{G}_q$ are compact in $D(\FK{G}_q)$, and the functor $D(\FK{G}_q)\to D(\Xq)$ preserves compact objects by Corollary \ref{cor:1110}, it follows that \eqref{eq:1093} is an isomorphism whenever $M$ is a coherent projective in $\FK{G}_q$ and $N$ is in the localizing subcategory
\begin{equation}\label{eq:1118}
\operatorname{Loc}(\FK{G}_q)=\operatorname{Loc}(D(\FK{G}_q)^{\heartsuit})\subset D(\FK{G}_q)
\end{equation}
generated by the heart of the derived category.  But, this localizing subcategory is all of $D(\FK{G}_q)$ \cite[\S 5.10]{krause10}, so that the map \eqref{eq:1093} is an isomorphism whenever $M$ is coherent projective and $N$ is arbitrary in $D(\FK{G}_q)$.  We use the more precise identification
\[
\operatorname{Loc}(\operatorname{proj}\FK{G}_q)=D(\FK{G}_q)
\]
\cite[\S 5.10]{krause10}, where $\operatorname{proj}\FK{G}_q$ denotes the category of coherent projectives in $\FK{G}_q$, to see now that \eqref{eq:1093} is an isomorphism at arbitrary $M$ and $N$.
\end{proof}

\subsection{A stronger centrality for the embedding $\kappa^\ast$}
\label{sect:Kempf_Z}

We have argued above that one should consider $\QCoh(\Xq)$ as a flat family of tensor categories over the flag variety.  Hence the appropriate ``Drinfeld center" for $\QCoh(\Xq)$ should be the centralizer
\[
Z_{\dG/\dB}(\QCoh(\Xq)):=\text{the centralizer of $\QCoh(\dG/\dB)$ in }Z(\QCoh(\Xq)).
\]
We have the central structure $\tilde{\kappa}^\ast:\FK{G}_q\to Z(\QCoh(\Xq))$ for universal restriction which is provided by the $R$-matrix in the expected way,
\begin{equation}\label{eq:1237}
\gamma_{N,M}:N\ot \kappa^\ast(M)\to \kappa^\ast(M)\ot N,\ \ \gamma_{N,M}(n,m)=R^{21}(m\ot n).
\end{equation}
One sees directly that for $N$ in the image of the embedding $\zeta:\QCoh(\dG/\dB)\to \QCoh(\Xq)$, $\gamma_{N,M}$ reduces to the symmetry provided by the exponentiated Killing form $\operatorname{symm}_{N,M}$.  Since the Killing form is a symmetric, square zero operation on all such products $N\ot M$ it follows that the image of $\tilde{\kappa}^\ast$ does in fact lie in the centralizer of $\QCoh(\dG/\dB)$, so that we restrict to obtain the desired ``$\QCoh(\dG/\dB)$-linear" central structure
\begin{equation}\label{eq:1241}
\tilde{\mu}^\ast:\FK{G}_q\to Z_{\dG/\dB}(\QCoh(\Xq)).
\end{equation}
When we speak of the universal restriction functor as a central tensor functor we mean specifically the functor $\kappa^\ast:\FK{G}_q\to \QCoh(\Xq)$ along with the lift \eqref{eq:1241} provided by the half-braidings \eqref{eq:1237}.

\subsection{The geometry of $\QCoh(\Xq)$ and the quantum Frobenius kernel}

We conclude the section with some remarks on the ``geometry" of the category $\QCoh(\Xq)$, and how this geometry does (not) transfer to $\FK{G}_q$.
\par

As was explained previously, the $\QCoh(\dG/\dB)$-action on $\QCoh(\Xq)$ endows this category with a kind of external geometry.  This external geometry is realized by the construction of sections
\[
\QCoh(\Xq)|_f:=\QCoh(S)\ot_{\QCoh(\dG/\dB)}\QCoh(\Xq)
\]
over arbitrary maps $f:S\to \dG/\dB$.  (See section \ref{sect:fibers_G/B}, and compare with \cite{gaitsgory15}.)  This local structure on $\QCoh(\Xq)$ does not restrict to provide any local structure on the category $\FK{G}_q$.  In particular, $\FK{G}_q$ is not stable under the action of $\QCoh(\dG/\dB)$ on $\QCoh(\Xq)$.  This is clear since $\QCoh(\dG/\dB)$ has infinitely many invertible objects, provided by the collection of $1$-dimensional $\dB$-representations, while $\FK{G}_q$ has only finitely many invertible objects.
\par

However, the action of $\QCoh(\dG/\dB)$ provides an additional ``internal geometry" on the category $\QCoh(\Xq)$.  This internal geometry is realized by the inner-Homs $\sHom_{\Xq}$, as described in Lemma \ref{lem:994}, and these inner-Homs provide an enhancement for $\QCoh(\Xq)$ in the category of sheaves over the flag variety:
\[
\QCoh(\Xq)\ \rightsquigarrow\ \QCoh^{\Enh}(\Xq):=(\operatorname{obj}\QCoh(\Xq),\ \sHom_{\Xq}).
\]
This internal geometry is less stable than the external geometry discussed above, but \emph{does} endow the category of small quantum group representations with a kind of local structure over $\dG/\dB$.  In particular, since the embedding $\kappa^\ast:\FK{G}_q\to \QCoh(\Xq)$ is fully faithful, the aforementioned enhancement for $\QCoh(\Xq)$ \emph{will} restrict to a monoidal enhancement for the small quantum group
\[
\FK{G}_q\ \rightsquigarrow\ (\operatorname{obj}\FK{G}_q,\ \sHom_{\FK{G}_q}),
\]
where we take $\sHom_{\FK{G}_q}:=\sHom_{\Xq}(\kappa^\ast-,\kappa^\ast-)$.  This local structure for the category of quantum group representations is not discussed in detail in this text, but plays a fundamental role in our related study of support theory \cite[Part II]{negronpevtsova5}.
\par

In Section \ref{sect:shom1} and \ref{sect:shom2} we clarify some basic claims about the internal geometry of $\QCoh(\Xq)$, as discussed above.  We subsequently explain how the enhancement $\QCoh^{\Enh}(\Xq)$ for $\QCoh(\Xq)$ readily derives to provide an enhancement for the derived category $D(\Xq)$.  We then turn to a discussion of relationships between $\Xq$ and the quantum Borels $\msc{B}_\lambda$.

\section{Structure of $\sHom$ I: Linearity, composition, and tensoring}
\label{sect:shom1}

In the next two sections we provide an analysis of the inner-Hom, or sheaf-Hom, functor for the action of $\QCoh(\dG/\dB)$ on $\QCoh(\Xq)$.  Here we show that sheaf-Homs admit natural composition and monoidal structure maps, and so provide a monoidal enhancement $\QCoh^{\Enh}(\Xq)$ for the monoidal category of sheaves on the half-quantum flag variety.  Many of the results of this section are completely general and completely formal.  So some proofs are sketched and/or delayed to the appendix.
\par

In the subsequent section, Section \ref{sect:shom2}, we provide an explicit description of the sheafy morphisms $\sHom_{\Xq}$ and describe objects in $\QCoh(\Xq)$ which are projective for this functor.

\begin{remark}
Sections \ref{sect:shom1}--\ref{sect:Enh_derived} are purely technical, as we are only verifying some expected properties for the inner-Hom functor, both at the abelian and derived levels.  So one might skim these sections on a first reading.  The next substantial results come in Sections \ref{sect:fibers_G/B} and \ref{sect:fibers_sHom}, where we calculate the fibers of the category $\QCoh(\Xq)$ and the mapping sheaves $\sRHom_{\Xq}$ over the flag variety, respectively.
\end{remark}

\subsection{$\QCoh(\dG/\dB)$-linearity of sheaf-Hom}

The adjoint to the identity map $id:\sHom_{\Xq}(M,N) \to \sHom_{\Xq}(M,N)$ provides an evaluation morphism
\[
ev\colon\sHom_{\Xq}(M,N)\star M\to N.
\]
The evaluation is just the counit for the $(\star,\sHom)$-adjunction.  For any $\msc{F}$ in $\QCoh(\dG/\dB)$ the map
\[
id\ot ev:\msc{F}\star(\sHom_{\Xq}(M,N)\star M)\to \msc{F}\star N
\]
provides a natural morphism $\msc{F}\star\sHom_{\Xq}(M,N)\to \sHom_{\Xq}(M,\msc{F}\star N)$ in the category of sheaves over the flag variety.  In analyzing this natural morphism it is helpful to consider a notion of projectivity for the sheaf-Hom functor.

\begin{definition}\label{def:loc_inj/proj}
An object $M$ in $\QCoh(\Xq)$ is called relatively projective (resp.\ relatively injective) if the functor $\sHom_{\Xq}(M,-)$ (resp.\ $\sHom_{\Xq}(-,M)$) is exact.
\end{definition}

Discussions of relatively injective and projective sheaves are provided in Lemma \ref{lem:loc_inj} and Section \ref{sect:proof1090} below, respectively.

\begin{lemma}[{cf.\ \cite[Lemma 3.3]{ostrik03}}]\label{lem:sHom-linear}
Consider the structural map
\begin{equation}\label{eq:469}
\msc{F}\star \sHom_{\Xq}(M,N)\to \sHom_{\Xq}(M,\msc{F}\star N)
\end{equation}
at $M$ and $N$ in $\QCoh(\Xq)$, and $\msc{F}$ in $\QCoh(\dG/\dB)$.  The map \eqref{eq:469} is an isomorphism under any of the following hypotheses:
\begin{itemize}
\item $\msc{F}$ is a coherent vector bundle.
\item $M$ is relatively projective and $\msc{F}$ is coherent.
\item $M$ is coherent and relatively projective, and $\msc{F}$ is arbitrary.
\item $M$ admits a presentation $E'\to E\to M$ by coherent, relatively projective sheaves and $\msc{F}$ is flat.
\end{itemize}
\end{lemma}

\begin{proof}
When $\msc{F}$ is a vector bundle it is dualizable, so that we have (explicit) natural isomorphisms \cite[Lemma 2.2]{ostrik03}
\[
\begin{array}{rl}
\Hom_{\dG/\dB}(-,\msc{F}\star \sHom_{\Xq}(M,N))& \cong \Hom_{\Xq}(\msc{F}^\vee\ot-,\sHom_{\Xq}(M,N))\vspace{1mm}\\
&\cong \Hom_{\Xq}((\msc{F}^\vee\ot -)\star M,N)\vspace{1mm}\\
&\cong \Hom_{\Xq}(-\star M,\msc{F}\star N)\vspace{1mm}\\
&\cong \Hom_{\dG/\dB}(-,\sHom_{\Xq}(M,\msc{F}\star N))
\end{array}
\]
and deduce an isomorphism $\msc{F}\star \sHom_{\Xq}(M,N)\cong \sHom_{\Xq}(M,\msc{F}\star N)$ via Yoneda.  One traces the identity map through the above sequence to see that this isomorphism is in fact just the structure map \eqref{eq:469}.  The second statement follows from the first, after we resolve $\msc{F}$ by vector bundles.  The third statement follows from the second and the fact that $\sHom_{\Xq}(M,-)$ commutes with colimits in this case.  The fourth statement follows from the second by resolving $M$ by relatively projective coherent sheaves.
\end{proof}

\begin{remark}
We will see at Corollary \ref{cor:enough_proj} below that all coherent sheaves in $\QCoh(\Xq)$ admit a resolution by relative projectives.  So the fourth point of Lemma \ref{lem:sHom-linear} simply says that $\sHom_{\Xq}(M,-)$ is linear with respect to the action of flat sheaves over $\dG/\dB$, whenever $M$ is coherent.
\end{remark}

\begin{remark}
One should compare Lemma \ref{lem:sHom-linear} with the familiar case of a linear category.  For a linear category $\msc{C}$, i.e.\ a module category over $\Vect$, we have natural maps $V\ot_k\Hom_{\msc{C}}(A,B)\to \Hom_{\msc{C}}(A,V\ot_k B)$ which one generally thinks of as associated to an identification between $V\ot_k-$ and a large coproduct.  This map is an isomorphism provided $V$ is sufficiently finite (dualizable), or $A$ is sufficiently finite (compact).
\end{remark}

\subsection{Enhancing $\QCoh(\Xq)$ via sheaf-Hom}

Evaluation for $\sHom_{\Xq}$ provides natural composition functions
\[
\circ:\sHom_{\Xq}(M,N)\ot\sHom_{\Xq}(L,M)\to \sHom_{\Xq}(L,N)
\]
which are adjoint to the maps
\[
\sHom_{\Xq}(M,N)\ot\sHom_{\Xq}(L,M)\star L\overset{id\ot ev}\longrightarrow \sHom_{\Xq}(M,N)\star M\overset{ev}\to N.
\]
We also have monoidal structure maps
\[
t:\sHom_{\Xq}(M,N)\ot\sHom_{\Xq}(M',N')\to \sHom(M\ot M',N\ot N')
\]
which are adjoint to the composition
\begin{equation}\label{eq:1380}
\begin{array}{l}
\sHom_{\Xq}(M_1,N_1)\ot \sHom_{\Xq}(M_2,N_2)\star (M_1\ot M_2)\vspace{2mm}\\
\hspace{1cm}\overset{symm}\longrightarrow(\sHom_{\Xq}(M_1,N_1)\star M_1)\ot (\sHom_{\Xq}(M_2,N_2)\star M_2)\vspace{2mm}\\
\hspace{2cm}\overset{ev\ot ev}\longrightarrow N_1\ot N_2.
\end{array}
\end{equation}
One can check the following basic claim, for which we provide a proof in Appendix \ref{sect:A} (cf.\ \cite{ostrik03,etingofostrik04}).

\begin{proposition}\label{prop:enriched}
The composition and monoidal structure maps for $\sHom_{\Xq}$ are associative, and are compatible in the sense of Section \ref{sect:enrich}.
\end{proposition}

This result says that the pairing of objects from $\QCoh(\Xq)$, along with the sheaf-morphisms $\sHom_{\Xq}$, constitutes a monoidal category enriched in the symmetric monoidal category of quasi-coherent sheaves on $\dG/\dB$.

\begin{definition}
We let $\QCoh^{\Enh}(\Xq)$ denote the enriched monoidal category
\[
\QCoh^{\Enh}(\Xq):=\big(\operatorname{obj}\QCoh(\Xq),\ \sHom_{\Xq}\big),
\]
with composition and tensor structure maps as described above.
\end{definition}

As the notation suggests, the category $\QCoh^{\Enh}(\Xq)$ does in fact provide an enhancement for the category of sheaves on the half-quantum flag variety.

\begin{theorem}\label{thm:enhA}
The adjunction isomorphism
\[
\Gamma(\dG/\dB,\sHom_{\Xq}(-,-))\overset{\cong}\to \Hom_{\Xq}(-,-)
\]
induces an isomorphism of monoidal categories
\[
\Gamma(\dG/\dB,\QCoh^{\Enh}(\Xq))\overset{\cong}\to \QCoh(\Xq).
\]
\end{theorem}

The proof of Theorem \ref{thm:enhA} is outlined in Appendix \ref{sect:A}, and is essentially the same as \cite[proof of Lemma 3.4.9]{riehl14}, for example.  Let us enumerate the main points here in any case.
\par

Recall that the linear structure on $\QCoh(\Xq)$ corresponds to an action of $\Vect$ on $\QCoh(\Xq)$.  The inner-Homs with respect to this action are the usual vector space of morphisms $\Hom_{\Xq}$ with expected evaluation
\[
ev:\Hom_{\Xq}(M,N)\ot_k M\to N,\ \ f\ot m\mapsto f(m).
\]
This evaluation map specifies a unique binatural morphism
\begin{equation}\label{eq:978}
\Hom_{\Xq}(M,N)\ot_k \O_{\dG/\dB}\to \sHom_{\Xq}(M,N)
\end{equation}
which is compatible with evaluation, in the sense that the diagram
\begin{equation}\label{eq:982}
\xymatrix{
\Hom_{\Xq}(M,N)\ot_k M\ar[d] \ar[r] & N\\
\sHom_{\Xq}(M,N)\star M\ar[ur] & 
}
\end{equation}
commutes.  If we consider $\Hom_{\Xq}(M,N)$ as a constant sheaf of vector spaces on $\dG/\dB$, the map \eqref{eq:978} is specified by a morphism of sheaves $\Hom_{\Xq}(M,N)\to \sHom_{\Xq}(M,N)$, which is in turn specified by its value on global sections.
\par

One can check that the global sections of the map $\Hom_{\Xq}(M,N)\to \sHom_{\Xq}(M,N)$ of \eqref{eq:978} recovers the adjunction isomorphism referenced in Theorem \ref{thm:enhA}.  One then uses compatibility with evaluation \eqref{eq:982} to see that the adjunction isomorphism is compatible with composition and the monoidal structure maps, and hence that $\QCoh^{\Enh}(\Xq)$ provides an enhancement of $\QCoh(\Xq)$ over the flag variety, as claimed.

\subsection{Implications for the quantum Frobenius kernel $\FK{G}_q$}

Recall that the functor $\kappa^\ast:\FK{G}_q=\Coh(\dG)^{G_q}\to \QCoh(\Xq)$ is a monoidal embedding.  By restricting along this embedding the enhancement $\QCoh^{\operatorname{Enh}}(\Xq)$ for $\QCoh(\Xq)$ restricts to a monoidal enhancement for the quantum Frobenius kernel.

\begin{theorem}
Let $\sHom_{\FK{G}_q}$ denote the restriction of the inner-$\operatorname{Hom}$s $\sHom_{\Xq}$ to the full monoidal subcategory $\FK{G}_q$ via the functor $\kappa^\ast$,
\[
\sHom_{\FK{G}_q}(M,N):=\sHom_{\Xq}(\kappa^\ast M,\kappa^\ast N).
\]
Then the pairing $(\operatorname{obj}\FK{G}_q,\ \sHom_{\FK{G}_q})$ provides a monoidal enhancement for the quantum Frobenius kernel $\FK{G}_q$ in the category of quasi-coherent sheaves over the flag variety $\dG/\dB$.
\end{theorem}

\subsection{A local-to-global spectral sequence}

The following lemma says that injectives in $\QCoh(\Xq)$ are relatively injective for the sheaf-Hom functor.

\begin{lemma}\label{lem:loc_inj}
If $M$ is flat over $\Xq$ then the functor $\sHom_{\Xq}(M,-)$ sends injectives in $\QCoh(\Xq)$ to injectives in $\QCoh(\dG/\dB)$.  Additionally, when $I$ is an injective object $\QCoh(\Xq)$, the functor $\sHom_{\Xq}(-,I)$ is exact.
\end{lemma}

\begin{proof}
Suppose that $I$ is injective over $\Xq$ and that $M$ is flat.  Then the functor $\Hom_{\Xq}(-\ot M,I)$ is exact, and hence $\Hom_{\Xq}(-\star M,I)$ is an exact functor from $\QCoh(\dG/\dB)$.  Via adjunction we find that the functor $\Hom_{\Xq}(-,\sHom_{\Xq}(M,I))$ is exact.  So $\sHom_{\Xq}(M,I)$ is injective.
\par

Now, if $I$ is injective then for all vector bundles $\msc{E}$ over $\dG/\dB$ the operation $\Hom_{\dG/\dB}(\msc{E}\star -,I)$ is exact, which implies that each functor
\[
\Hom_{\dG/\dB}(\msc{E},\sHom_{\Xq}(-,I))
\]
is exact.  This is sufficient to ensure that $\sHom_{\Xq}(-,I)$ is exact.
\end{proof}

Recall that we have the natural identification
\[
\begin{array}{rl}
\Hom_{\Xq}(M,N)&\cong \Hom_{\dG/\dB}(\O_{\dG/\dB},\sHom_{\Xq}(M,N))\vspace{1mm}\\
&=\Gamma(\dG/\dB,\sHom_{\Xq}(M,N))
\end{array}
\]
provided by adjunction.  We therefore obtain a natural map
\[
\RHom_{\Xq}(M,N)\to \operatorname{R}\Gamma(\dG/\dB,\sRHom_{\Xq}(M,N)),
\]
where we derive $\sHom_{\Xq}(M,-)$ by taking injective resolutions.  Lemma \ref{lem:loc_inj} implies that this map is a quasi-isomorphism.

\begin{corollary}\label{cor:372}
The natural map $\RHom_{\Xq}(M,N)\to \operatorname{R}\Gamma(\dG/\dB,\sRHom_{\Xq}(M,N))$ is a quasi-isomorphism.  Hence we have a local-to-global spectral sequence
\[
H^\ast(\dG/\dB,\sExt^\ast_{\Xq}(M,N))\ \Rightarrow\ \Ext^\ast_{\Xq}(M,N)
\]
at arbitrary $M$ and $N$ in $\QCoh(\Xq)$.
\end{corollary}

\begin{remark}
Spectral sequences analogous to Corollary \ref{cor:372} can be found in much earlier works of Suslin, Friedlander, and Bendal \cite[Theorem 3.6]{suslinfriedlanderbendel97b}.  So certain pieces of the enhancement $\QCoh^{\Enh}(\Xq)$, and the universal restriction functor, had already been employed in works which appeared as early as the 90's.
\end{remark}

\section{Structure of $\sHom$ II: explicit description of sheaf-Hom}
\label{sect:shom2}

We show that the morphisms $\sHom_{\Xq}(M,N)$ are explicitly the descent
\begin{equation}\label{eq:nmbvz}
\sHom_{\Xq}(M,N)=\text{descent of }\Hom_{\QCoh(\dG)^{\uqB}}(M,N)^\sim
\end{equation}
of morphisms in the category $\QCoh(\dG)^{\uqB}$ of $\uqB$-equivariant sheaves over $\dG$, whenever $M$ is coherent.  Here $\uqB$ is taken to act trivially on $\dG$, so that $\uqB$ acts by $\O_{\dG}$-linear endomorphisms on such sheaves.  Also, implicit in the formula \eqref{eq:nmbvz} is a claim that the above morphisms over $\QCoh(\dG)^{\uqB}$ admit natural, compatible actions of $\O(\dG)$ and $\dB$, so that the associated sheaf over $\dG$ is $\dB$-equivariant.  We then apply descent to produce a corresponding sheaf on the flag variety.  (See Proposition \ref{prop:shom_exp} below.)  We also describe the composition and tensor structure maps for $\sHom_{\Xq}$ in terms of the formula \eqref{eq:nmbvz}.

\subsection{$B_q$-equivariant structure on $\sHom_{\dG}$}
\label{sect:HomG}

Recall that for an algebraic group $H$ acting on a finite-type scheme $Y$, the usual sheaf-Hom functor $\sHom_Y$ provides the inner-Homs for the tensor action of $\Coh(Y)^H$ on itself.  We claim that this is also true when we act via a quantum group.  Specifically, we claim that when $M$ and $N$ are $B_q$-equivariant sheaves on $\dG$, and $M$ is coherent, the sheaf-morphisms $\sHom_{\dG}(M,N)$ admit a natural $B_q$-equivariant structure.  Indeed, the functor $\sHom_{\dG}(M,-)$, with its usual evaluation morphism, provides the right adjoint for the action of $M$ on $\QCoh(\Xq)$.  We describe the equivariant structure on $\sHom_{\dG}(M,N)$ explicitly below.
\par
Let $-|_{B_q}:\QCoh(\Xq)\to \Rep B_q$ denote the global sections functor, which we understand as adjoint to the vector bundle map $E_-:\Rep B_q\to \QCoh(\Xq)$.  The $U_q(\mfk{b})$-actions on $M|_{B_q}$ and $N|_{B_q}$ induce an action of $U_q(\mfk{b})$ on the linear morphisms $\Hom_k(M|_{B_q},N|_{B_q})$, via the usual formula $(x\cdot f)(m)=x_1f(S(x_2)m)$.  The compatibility between the $\O(\dG)$ and $U_q(\mfk{b})$-actions on $M$ and $N$ ensure that the subspace of \emph{right} $\O(\dG)$-linear maps
\[
\Hom_{\O(\dG)}(M|_{B_q},N|_{B_q})\subset \Hom_k(M|_{B_q},N|_{B_q})
\]
forms a $U_q(\mfk{b})$-subrepresentation, and the natural \emph{left} action of $\O(\dG)$ provides the space $\Hom_{\O(\dG)}(M|_{B_q},N|_{B_q})$ with the structure of a $U_q(\mfk{b})$-equivariant $\O(\dG)$-module.
\par

Now, provided $M$ is coherent, the action of $U_q(\mfk{b})$ on the above $\O(\dG)$-linear morphism space integrates to a $B_q$-action, so that $\Hom_{\O(\dG)}(M|_{B_q},N|_{B_q})$ is naturally an object in the category of relative Hopf modules $_{\O(\dG)}\mbf{M}^{\O(B_q)}$.  Indeed, one can observe such integrability by resolving $M$ by vector bundles $E_W\to E_V\to M$.  Since
\[
\Gamma(\dG,\sHom_{\dG}(M,N))=\Hom_{\O(\dG)}(M|_{B_q},N|_{B_q})
\]
we see that $\sHom_{\dG}(M,N)$ is naturally a $B_q$-equivariant sheaf on $\dG$.
\par

We now understand that we have an endofunctor
\[
\sHom_{\dG}(M,-):\QCoh(\Xq)\to \QCoh(\Xq)
\]
provided by usual sheaf-Hom.  The evaluation maps exhibiting $\sHom_{\dG}(M,-)$ as the right adjoint to the functor $-\ot M$ are the expected ones
\[
ev:\sHom_{\dG}(M,N)\ot M\to N,\ \ f\ot m\mapsto f(m).
\]

\subsection{Explicit description of $\sHom_{\Xq}$}

Suppose that $M$ is a coherent sheaf over $\Xq$.  We consider the ``quotient stack" $Y_q=\dG/\uqB$ via the trivial action, and the corresponding category of sheaves $\QCoh(Y_q)=\QCoh(\dG)^{\uqB}$.  These are just sheaves with an action of $\uqB$ by sheaf endomorphisms, and we have the sheafified morphisms
\[
\sHom_{Y_q}=\sHom_{\dG}(-,-)^{\uqB}.
\]
Left exactness of the invariants functor ensures that the sections of $\sHom_{Y_q}$ over opens are as expected,
\[
\begin{array}{rl}
\sHom_{Y_q}(M,N)(U)&=\Hom_{\O(U)}(M(U),N(U))^{\uqB}\\
&=\Hom_{\O(U)\ot \uqB}(M(U),N(U)).
\end{array}
\]

We have the faithful inclusion $\QCoh(X_q)\to \QCoh(Y_q)$, as in Section \ref{sect:q_mon}, and restrict the domain of $\sHom_{Y_q}$ to obtain an operation
\[
\sHom_{Y_q}:\QCoh(\Xq)^{op}\times \QCoh(\Xq)\to \QCoh(\dG)^{\dB},
\]
where the natural $\dB$-action on $\sHom_{Y_q}(M,N)$ is deduced as in Section \ref{sect:HomG}.

\begin{proposition}\label{prop:shom_exp}
Let $M$ be in $\Coh(\Xq)$.  Then, at arbitrary $N$ in $\QCoh(\Xq)$, we have a natural identification
\[
\sHom_{\Xq}(M,N)=\text{\rm descent of the $\dB$-equivariant sheaf }\sHom_{Y_q}(M,N).
\]
Under the subsequent natural isomorphism $\pi^\ast\sHom_{\Xq}(M,-)\cong \sHom_{Y_q}(M,-)$, the evaluation maps for $\sHom_{\Xq}$ are identified with the morphisms
\[
\sHom_{Y_q}(M,N)\ot M\to N,\ \ f\ot m\mapsto f(m).
\]
\end{proposition}

\begin{proof}
By the materials of Section \ref{sect:HomG} we have the adjunction
\[
\Hom_{\Xq}(L\ot M,-)\cong \Hom_{\Xq}(L,\sHom_{\dG}(M,-)).
\]
For any $F$ in $\QCoh(\dG)^{\dB}\subset \QCoh(\Xq)$ we have
\[
\Hom_{\Xq}(F,-)=\Hom_{\QCoh(\dG)^{\dB}}(F,(-)^{\su}).
\]
So for $\msc{F}$ in $\QCoh(\dG/\dB)$ the above two formulae give
\begin{equation}\label{eq:1409}
\begin{array}{rl}
\Hom_{\Xq}(\msc{F}\star M,-)&=\Hom_{\Xq}(\pi^\ast(\msc{F})\ot M,-)\vspace{1mm}\\
&=\Hom_{\Xq}(\pi^\ast(\msc{F}),\sHom_{\dG}(M,-))\vspace{1mm}\\
&=\Hom_{\QCoh(\dG)^{\dB}}(\pi^\ast(\msc{F}),\sHom_{\dG}(M,-)^{\su})\vspace{1mm}\\
&\cong \Hom_{\dG/\dB}(\msc{F},\text{desc of\ }\sHom_{Y_q}(M,-)).
\end{array}
\end{equation}
The above formula demonstrates the descent of the sheaf $\sHom_{Y_q}(M,-)$ is the right adjoint to $-\star M$, and hence identifies $\sHom_{\Xq}(M,-)$ with the descent of $\sHom_{Y_q}(M,-)$.  Tracing the identity map through the sequence \eqref{eq:1409} calculates evaluation for $\sHom_{\Xq}$ as the expected morphism
\[
\pi^\ast\sHom_{\Xq}(M,N)\ot M\cong \sHom_{Y_q}(M,N)\ot M\to N,\ \ f\ot m\mapsto f(m).
\]
\end{proof}

For arbitrary $M$ in $\QCoh(\Xq)$ we may write $M$ as a colimit $M=\varinjlim_\alpha M_\alpha$ of coherent sheaves to obtain
\[
\sHom_{\Xq}(M,N)=\varprojlim_\alpha \sHom_{\Xq}(M_\alpha,N),
\]
where the final limit is the limit in the category of quasi-coherent sheaves over $\dG/\dB$.  So, Proposition \ref{prop:shom_exp} provides a complete description of the inner-Hom functor.

\subsection{Relatively projective sheaves}
\label{sect:relative_p}

For any finite-dimensional $B_q$-representation $V$ we consider the functor
\[
\Hom_{\su}(V,-):\QCoh(\Xq)\to \QCoh(\dG)^{\dB},
\]
where specifically $\Hom_{\su}(V,-)=\Hom_k(V,-)^{\uqB}$.  We have the following description of sheaf-Hom for the equivariant vector bundles which refines the description of Proposition \ref{prop:shom_exp}.

\begin{proposition}
For any finite-dimensional $B_q$-representation $V$, there is a natural identification between $\sHom_{\Xq}(E_V,-)$ and the descent of the functor $\Hom_{\su}(V,-)$.
\end{proposition}

\begin{proof}
For any $\msc{F}$ in $\QCoh(\dG/\dB)$ we calculate
\[
\begin{array}{l}
\Hom_{\Xq}(\msc{F}\star E_V,-)=\Hom_{\Xq}(\pi^\ast(\msc{F})\ot_k V,-)\\
\hspace{1cm}=\Hom_{\Xq}(\pi^\ast(\msc{F}),-\ot_k V^\ast)\\
\hspace{1cm}=\Hom_{\Xq}(\pi^\ast(\msc{F}),\Hom_k(V,-))\\
\hspace{1cm}=\Hom_{\QCoh(\dG)^{\dB}}(\pi^\ast(\msc{F}),\Hom_k(V,-)^{\su})=\Hom_{\dG/\dB}(\msc{F},\operatorname{desc\ of\ }\Hom_{\su}(V,-)).
\end{array}
\]
Thus, by uniqueness of adjoints, we find that $\sHom_{\Xq}(E_V,-)$ is identified with the descent of the functor $\Hom_{\su}(V,-)$.
\end{proof}

As a corollary we observe a natural class of relative projective sheaves in $\QCoh(\Xq)$.

\begin{corollary}\label{cor:loc_proj}
Suppose that $V$ is a finite-dimensional $B_q$-representation which is projective over $\uqB$.  Then the functor $\sHom_{\Xq}(E_V,-)$ is exact.
\end{corollary}

\begin{proof}
To establish exactness of $\sHom_{\Xq}(E_V,-)$ it suffices to show that the functor $\Hom_{\su}(V,-)$ is exact.  But this follows by projectivity of $V$ over $\uqB$.
\end{proof}

Note that any coherent sheaf $M$ in $\QCoh(\Xq)$ admits a surjection $E\to M$ from an equivariant vector bundle $E=E_V$, with $V$ finite-dimensional and projective over $\uqB$.  So Corollary \ref{cor:loc_proj} implies that the category of coherent sheaves over $\Xq$ has enough relative projectives.

\begin{corollary}\label{cor:enough_proj}
Any coherent sheaf $M$ in $\QCoh(\Xq)$ admits a resolution $\cdots\to E^{-1}\to E^0\to M$ by coherent, relatively projective sheaves $E^i$.
\end{corollary}

\subsection{A proof of Lemma \ref{lem:1090}}
\label{sect:proof1090}

At Lemma \ref{lem:1090} above, we have claimed that the vector bundle $E_V$ is compact in the unbounded derived category $D(\Xq)$ whenever the given $B_q$-representation $V$ is finite-dimensional and projective over $\uqB$.  We can now prove this result.

\begin{proof}[Proof of Lemma \ref{lem:1090}]
We have the functor $\sHom_{\Xq}(E_V,-)$ which is the descent of the functor $\Hom_{\su}(V,-)$.  This functor is exact, and finiteness of $V$ implies that $\Hom_{\su}(V,-)$ commutes with set indexed sums.  Hence these inner-Homs provide a well-defined operation
\[
\sHom_{\Xq}(E_V,-):D(\Xq)\to D(\dG/\dB)
\]
which commutes with set indexed sums.  We then have
\[
\Ext_{\Xq}^i(E_V,-)=\Ext_{\dG/\dB}^i(\O_{\dG/\dB},\sHom_{\Xq}(E_V,-))
\]
at each integer $i$ by Corollary \ref{cor:372}.  Compactness of $\O_{\dG/\dB}$ over $\dG/\dB$ therefore implies compactness of $E_V$ over $\Xq$.
\end{proof}

\subsection{Composition and tensor structure maps}
\label{sect:comp_ot_exp}

Suppose that $M$ is coherent in $\QCoh(\Xq)$.   From Proposition \ref{prop:shom_exp} we have an identification
\begin{equation}\label{eq:1450}
\pi^\ast\sHom_{\Xq}(M,-)\cong \sHom_{Y_q}(M,-)
\end{equation}
under which the evaluation morphisms for $\sHom_{\Xq}$ pull back to the usual evaluation morphisms
\[
\sHom_{Y_q}(M,N)\ot M\to N,\ \ f\ot m\mapsto f(m)
\]
for $\sHom_{Y_q}$.  It follows that, under the identification \eqref{eq:1450}, the composition and tensor maps
\[
\circ:\sHom_{\Xq}(M,N)\ot\sHom_{\Xq}(L,M)\to \sHom_{\Xq}(L,N)
\]
and
\[
tens:\sHom_{\Xq}(M,N)\ot\sHom_{\Xq}(M',N')\to \sHom_{\Xq}(M\ot M',N\ot N')
\]
pull back to the expected morphisms
\[
\begin{array}{c}
\pi^\ast\circ:\sHom_{Y_q}(M,N)\ot\sHom_{Y_q}(L,M)\to \sHom_{Y_q}(L,N)\\
(f, g)\mapsto f\circ g
\end{array}
\]
and
\begin{equation}\label{eq:1688}
\begin{array}{c}
\pi^\ast tens:\sHom_{Y_q}(M,N)\ot\sHom_{Y_q}(M',N')\to \sHom_{Y_q}(M\ot M',N\ot N')\\
(f,f')\mapsto f\ot f'.
\end{array}
\end{equation}

\begin{remark}
To be clear that, when evaluating a product $f\ot f'$ \eqref{eq:1688} on homogenous sections of $M\ot M'$, we pick up ``Koszul signs"
\[
(f\ot f')(m\ot m')=q^{(\deg(f'),\deg(m))}f(m)\ot f'(m')=\pm f(m)\ot f'(m').
\]
The is due to the presence of the half-braiding in the formula \eqref{eq:1380}.  As usual, these signs vanish at odd order $q$.
\end{remark}

\section{The enhanced derived category}
\label{sect:Enh_derived}

In Section \ref{sect:shom1} we saw that the sheaf-Hom functor $\sHom_{\Xq}$ provides a monoidal enhancement for $\QCoh(\Xq)$ over the classical flag variety.  At this point we want to provide a corresponding enhancement for the (unbounded) derived category $D(\Xq)$ of quasi-coherent sheaves over the half-quantum flag variety.  Given the information we have already collected, this move to the derived setting is a relatively straightforward.  We record some of the details here.

\subsection{$\sHom_{\Xq}$ for complexes}

Let $\operatorname{dgQCoh}(\Xq)$ denote the category of quasi-coherent dg sheaves on $\Xq$.  This is the category of quasi-coherent sheaves $M$ with a grading $M=\oplus_{n\in \mbb{Z}} M^n$, and a degree $1$ square zero map $d_M:M\to M$.  Morphisms in this category are the usual morphisms of complexes.  We similarly define $\operatorname{dgQCoh}(\dG/\dB)$.  Consider the forgetful functor $f:\operatorname{dgQCoh}(\Xq)\to \QCoh(\Xq)$.
\par

Let $M$ and $N$ be in $\operatorname{dgQCoh}(\Xq)$.  For an open $U\subset \dG/\dB$ and $M$ in $\QCoh(\Xq)$ take $M|_U:=\O_U\star M$.  A section
\[
s:\O_U\to \sHom_{\Xq}(fM,fN)|_U
\]
over an open $U\subset \dG/\dB$ is said to be homogenous of degree $n$ if, for all $i\in \mbb{Z}$, the restriction $s|_{M^i}:M^i|_U\to N|_U$ has image in $N^{i+n}|_U$.  Here $s:M|_U\to N|_U$ is specifically the composite
\[
M|_U=\O_U\star M|_U\to \sHom_{\Xq}(fM,fN)|_U\star M|_U\overset{ev|_U}\to N|_U.
\]
The collection of degree $n$ maps in $\sHom_{\Xq}(fM,fN)$ forms a subsheaf which we denote
\[
\sHom_{\Xq}^n(M,N)\subset \sHom_{\Xq}(fM,fN).
\]
The sum of all homogeneous morphisms also provides a subsheaf
\[
\oplus_{n\in \mbb{Z}}\sHom^n_{\Xq}(M,N)\subset \sHom_{\Xq}(fM,fN).
\]

\begin{definition}
For $M$ and $N$ in $\operatorname{dgQCoh}(\Xq)$ we define the inner-$\Hom$ complex to be the dg sheaf consisting of all homogenous inner morphisms
\[
\sHom_{\Xq}(M,N)=\oplus_{n\in \mbb{Z}}\sHom^n_{\Xq}(M,N)
\]
equipped with the usual differential $d_{\sHom(M,N)}=(d_N)_\ast-(d_M)^\ast$.
\end{definition}

We note that, when $M$ is bounded above and $N$ is bounded below, this complex is just the expected one
\[
\sHom_{\Xq}(M,N)=\oplus_{n\in \mbb{Z}}(\oplus_{i}\sHom_{\Xq}(M^i,N^{i+n}))\ \ \text{along with } d_{\sHom}.
\]
Evaluation for $\sHom_{\Xq}(fM,fN)$ restricts to an evaluation map for the Hom complex $ev:\sHom_{\Xq}(M,N)\star M\to N$.  This evaluation map induces an adjunction
\[
\Hom_{\operatorname{dgQCoh}(\Xq)}(\msc{F}\star M,N)\overset{\cong}\to \Hom_{\operatorname{dgQCoh}(\dG/\dB)}(\msc{F},\sHom_{\Xq}(M,N)).
\]

\subsection{Enhancement for the derived category}
\label{sect:D_Enh}

We consider the central action
\[
\star:D(\dG/\dB)\times D(\Xq)\to D(\Xq)
\]
for the unbounded derived categories of quasi-coherent sheaves, and we have an adjunction
\[
\Hom_{D(\Xq)}(\msc{F}\star M,N)\cong \Hom_{D(\dG/\dB)}(\msc{F},\sRHom_{\Xq}(M,N))
\]
which one deduces abstractly, or from the adjunction at the cochain level described above.  Here the $\star$-action is derived by resolving $M$ by $K$-flat sheaves and $\sRHom_{\Xq}(M,N)=\sHom_{\Xq}(M,I_N)$ for a $K$-injective resolution $N\to I_N$ \cite[\href{https://stacks.math.columbia.edu/tag/079P}{Tag 079P}]{stacks}.  Evaluation
\[
ev:\sRHom_{\Xq}(M,N)\ot M\to N
\]
provides composition and tensor structure maps for derived sheaf-Hom $\sRHom_{\Xq}$ so that we obtain a monoidal category
\[
D^{\Enh}(\Xq)=(\text{ob}D(\Xq),\ \sRHom_{\Xq})
\]
which is enriched in the unbounded derived category of quasi-coherent sheaves on the flag variety.
\par

We have the derived global sections
\[
\mbb{H}^0(\dG/\dB,-)=\Hom_{D(\dG/\dB)}(\O_{\dG/\dB},-)
\]
and corresponding adjunction isomorphism
\begin{equation}\label{eq:1699}
\Hom_{D(\Xq)}\overset{\cong}\longrightarrow \mbb{H}^0(\dG/\dB,\sRHom_{\Xq}).
\end{equation}
Just as in the proof of Theorem \ref{thm:enhA} (Appendix \ref{sect:A}), one finds that the binatural isomorphism \eqref{eq:1699} realize $D^{\Enh}(\Xq)$ as an enhancement for the derved category $D(\Xq)$.

\begin{proposition}\label{prop:D_Enh}
The isomorphism \eqref{eq:1699} induces an isomorphism of monoidal categories
\[
\mbb{H}^0(\dG/\dB, D^{\operatorname{Enh}}(\Xq))\overset{\sim}\to D(\Xq).
\]
\end{proposition}

\begin{proof}
As in the proof of Theorem \ref{thm:enhA}, one sees that the adjunction map provides a morphism from the constant sheaf
\[
\Hom_{D(\Xq)}(M,N)\to \sRHom_{\Xq}(M,N)
\]
which recovers the standard evaluation maps for $\Hom_{D(\Xq)}$.  This is sufficient to deduce the result.
\end{proof}

\subsection{Coherent dg sheaves}
\label{sect:coherent_sheaves}

In order to gain a better handle on things, we might restrict our attention to objects in $D(\Xq)$ which satisfy certain finiteness conditions.  Abstractly, we consider the subcategory of perfect, or dualizable, objects.  In terms of our specific geometric presentations of these categories, we are interested in coherent dg sheaves.  Let us take a moment to describe this category clearly.
\par

We consider the derived category $D_{\coh}(\Xq)\subset D(\Xq)$ of coherent, $B_q$-equivariant dg sheaves over $\dG$.  These are, equivalently, complexes in $D(\Xq)$ with bounded coherent cohomology, or complexes in $D(\Xq)$ which are dualizable with respect to the product $\ot=\ot^{\operatorname{L}}$.

\subsection{Enhancements for the coherent derived categories}
\label{sect:perf_Enh}

When we consider the enhancements $D^{\Enh}(\Xq)$ provided above, we can be much more explicit about the evaluation and tensor maps when we restrict to the subcategories of coherent dg sheaves.  We let $D_{\coh}^{\Enh}(\Xq)$ denote the full (enriched, monoidal) subcategory consisting of such sheaves in $D^{\Enh}(\Xq)$.
\par

For coherent $M$ and bounded $N$ we can adopt any of the explicit models
\[
\begin{array}{rl}
\sRHom_{\Xq}(M,N)&=\sHom_{\Xq}(M,I_N)\ \text{or}\ \ \sHom_{\Xq}(P_M,I_N)\vspace{1mm}\\ &\hspace{1cm}\text{or}\ \ \sHom_{\Xq}(P_M,N)
\end{array}
\]
depending on our needs, where $P_M\to M$ and $N\to I_N$ are resolutions by relative projectives and injectives respectively.  The composition maps for these inner morphisms can then be obtained from composition at the dg level
\[
\circ:\sHom_{\Xq}(M,I_N)\ot\sHom_{\Xq}(P_L,M)\to \sHom_{\Xq}(P_L,I_N),
\]
as can the tensor structure maps
\[
tens:\sHom_{\Xq}(P_M,N)\ot\sHom_{\Xq}(P_{M'},N')\to \sHom_{\Xq}(P_M\ot P_{M'},N\ot N').
\]
Of course, the equivalence of Proposition \ref{prop:D_Enh} realizes $D_{\coh}^{\Enh}(\Xq)$ as an enhancement of $D_{\coh}(\Xq)$.

\subsection{Derived implications for the quantum Frobenius kernel}

As in the abelian setting, we restrict along the fully faithful embedding $\kappa^\ast:D(\FK{G}_q)\to D(\Xq)$ to obtain an enhancement
\[
D^{\Enh}(\FK{G}_q)=(\operatorname{obj}D(\FK{G}_q),\ \sRHom_{\FK{G}_q}),
\]
\[
\sRHom_{\FK{G}_q}:=\sRHom_{\Xq}(\kappa^\ast-,\kappa^\ast-),
\]
from the corresponding enhancement $D^{\Enh}(\Xq)$ for the half-quantum flag variety.  Furthermore, in the coherent setting, the analysis of Subsection \ref{sect:perf_Enh} restricts to provide explicit realizations of the composition and tensor structure maps for the derived sheaf-Homs $\sRHom_{\FK{G}_q}$.

\section{$\QCoh(\Xq)$ as a total space for the quantum Borels}
\label{sect:fibers_G/B}

We now have a basic understanding of the internal geometry for $\QCoh(\Xq)$, as well as the relationship between $\Xq$ and the quantum Frobenius kernel.  The point of this section is to examine some aspects of the external geometry for $\QCoh(\Xq)$.  In particular, we show that the fibers of the category $\QCoh(\Xq)$, considered as a flat family of tensor categories over the flag variety, recover the small quantum Borels of Section \ref{sect:borels}.  These fiber calculations are valid at both an abelian and derived level.  In this way $\QCoh(\Xq)$ serves as a total space for the family of Borels $\{\msc{B}_\lambda:\lambda:\Spec(K)\to \dG/\dB\}$ constructed in Section \ref{sect:borels}.

We furthermore show that the composite
\[
\FK{G}_q\overset{\kappa^\ast}\to \QCoh(\Xq)\to \QCoh(\Xq)|_\lambda\cong \msc{B}_\lambda
\]
recovers the restriction functor $\res_\lambda$.  Hence our understanding of $\kappa^\ast$ as a universal restriction functor for the Borels.
\par

We begin by recalling some basic information about categorical base change.

\subsection{Base change for cocomplete categories}
\label{sect:base_change}

Consider a presentable, and in particular cocomplete tensor category $\msc{C}$ (see Section \ref{sect:presentable}).  By a module category $\msc{M}$ over $\msc{C}$ we mean a presentable $k$-linear category with an associative action $\msc{C}\times \msc{M}\to \msc{M}$ (or $\msc{M}\times \msc{C}\to \msc{M}$) which commutes with arbitrary colimits in each factor.
\par

Fix a presentable tensor category $\msc{C}$, and presentable categories $\msc{M}$ and $\msc{N}$ with left and right actions of $\msc{C}$ respectively.  We consider, for any presentable $k$-linear category $\msc{D}$, the category
\[
\operatorname{Bil}_{\msc{C}}(\msc{N}\times \msc{M},\msc{D})
\]
of $\msc{C}$-bilinear maps $\msc{N}\times\msc{M}\to \msc{D}$ with natural transformations.  By a $\msc{C}$-bilinear map we mean a $k$-linear functor $F$ from the product $\msc{N}\times \msc{M}$ which commutes with arbitrary colimits in each factor, and comes equipped with a natural, associative, isomorphism $F(N,V\ot M)\cong F(N\ot V,M)$ for all $V$ in $\msc{C}$.  Natural transformations between bilinear functors are natural transformations of functors which respect the natural isomorphisms $F(-,-\ot-)\cong F(-\ot-,-)$.
\par

We similarly define the category of $\msc{C}$-linear functors $\operatorname{Fun}_{\msc{C}}(\msc{M},\msc{M}')$ for left (or right) module categories $\msc{M}$ and $\msc{M}'$.  These are functors which commute with arbitrary colimits and come equipped with associative natural isomorphisms $F(V\ot M)\cong V\ot F(M)$.  Morphisms in $\operatorname{Fun}_{\msc{C}}(\msc{M},\msc{M}')$ are natural transformations which respect the $\msc{C}$-linear structure.  When $\msc{C}=\Vect$ we take $\operatorname{Fun}_k=\operatorname{Fun}_{\Vect}$, and when $\msc{C}=\operatorname{QCoh}(X)$ for a scheme $X$ we take $\operatorname{Fun}_X=\operatorname{Fun}_{\operatorname{QCoh}(X)}$.  We note that $\operatorname{Fun}_k(\msc{M},\msc{M}')$ is just the category of cocontinuous $k$-linear functors between $k$-linear categories.
\par

Following \cite{kelly82,etingofnikshychostrik10,douglasspsnyder19}, a balanced tensor product $\msc{N}\ot_{\msc{C}}\msc{M}$ is a choice of presentable $k$-linear category, and a $\msc{C}$-bilinear map $F:\msc{N}\times\msc{M}\to \msc{N}\ot_{\msc{C}}\msc{M}$ for which restriction along $F$ provides an equivalence
\[
F^\ast:\operatorname{Fun}_k(\msc{N}\ot_{\msc{C}}\msc{M},\msc{D})\overset{\sim}\to \operatorname{Bil}_{\msc{C}}(\msc{N}\times \msc{M},\msc{D}),
\]
at arbitrary presentable $\msc{D}$.

\begin{remark}
In \cite{etingofnikshychostrik10,douglasspsnyder19} the authors are primarily concerned with finite categories, and so their functors are only assumed to commute with finite colimits.  In our cocomplete context it's appropriate to assume commutation with arbitrary set-indexed colimits.  (See \cite[Section 6.5]{kelly82} and \cite[Section 4.8]{lurieHA}.)
\end{remark}

By \cite[Section 6.5, Proposition 10.4]{kelly82} the product $\msc{N}\ot_k\msc{M}:=\msc{N}\ot_{\Vect}\msc{M}$ always exists (cf.\ \cite[Theorem 4]{franco13}).  The existence of colimits in the $2$-category of presentable linear categories implies that the product $\msc{N}\ot_{\msc{C}}\msc{M}$ exists as well \cite[Theorem 6.11]{bird84} \cite[Proposition 2.1.11]{chirvasitufreyd13}, as it is computed abstractly as the colimit of the diagram
\[
\xymatrix{
\msc{N}\ot_k\msc{C}\ot_k\msc{C}\ot_k\msc{M}\ar@<1ex>[r]\ar[r]\ar@<-1ex>[r]&
\msc{N}\ot_k\msc{C}\ot_k\msc{M}\ar@<.5ex>[r]\ar@<-.5ex>[r]&
\msc{N}\ot_k\msc{M}.
}
\]

We are particularly interested in the case where $\msc{C}=\operatorname{QCoh}(X)$ and $\msc{N}=\operatorname{QCoh}(Y)$ for schemes $X$ and $Y$, with action provided by pullback along a map $f:Y\to X$.  In this case the product
\[
\operatorname{QCoh}(Y)\ot_{\operatorname{QCoh}(X)}\msc{M}
\]
should be interpreted as a base change for $\msc{M}$.  The fiber of $\msc{M}$ at a given point $\lambda:\Spec(K)\to X$, where we have $\Vect(K)=\QCoh(\Spec(K))$, is defined as the base change
\[
\msc{M}|_\lambda:=\Vect(K)\ot_{\operatorname{QCoh}(X)}\msc{M}.
\]

We note that when $\msc{M}$ is a tensor category over $\msc{C}=\QCoh(X)$, i.e. a tensor category with a central tensor functor $\QCoh(X)\to \msc{M}$, the fiber $\msc{M}|_\lambda$ at a given $K$-point for $X$ inherits a unique $K$-linear tensor structure so that the universal (reduction) map $\msc{M}\to \msc{M}|_\lambda$ is a map of tensor categories.  Indeed, for any tensor category $\msc{A}$ over $\QCoh(X)$, the universal property of the product $\ot_{\QCoh(X)}$ and the tensor structures on $\msc{A}$ and $\msc{M}$ provide a natural tensor structure on the category $\msc{A}\ot_{\QCoh(X)}\msc{M}$, as outlined in \cite[Theorem 6.2]{greenough10}.

\subsection{The pullback functor $\iota^\ast_\lambda:\QCoh(\Xq)\to \msc{B}_\lambda$}
\label{sect:QCoh-linear}

Consider a geometric point $\lambda:\Spec(K)\to \dG/\dB$.  Pulling back along the $B_q$-equivariant map $\iota_\lambda:\dB_\lambda\to \dG$ provides a colimit preserving monoidal functor
\[
\iota_\lambda^\ast:\QCoh(\Xq)\to \msc{B}_\lambda.
\]
We let $\QCoh(\dG/\dB)$ act on $\msc{B}_\lambda$ via the fiber $\lambda^\ast:\QCoh(\dG/\dB)\to \Vect(K)$ and the $K$-linear structure on $\msc{B}_\lambda$.  We claim that under this action on $\msc{B}_\lambda$ the map $\iota^\ast_\lambda$ is $\QCoh(\dG/\dB)$-linear.  Indeed, for any sheaf $\msc{F}$ over $\dG/\dB$ we have
\[
\iota_\lambda^\ast(\msc{F}\star-)=\iota_\lambda^\ast(\pi^\ast \msc{F}\ot -)=\iota_\lambda^\ast\pi^\ast(\msc{F})\ot\iota_\lambda^\ast-
\]
and $\iota_\lambda^\ast\pi^\ast$ is naturally isomorphic to $\operatorname{unit}^\ast\lambda^\ast$, since we have the diagram
\[
\xymatrix{
\dB_\lambda\ar[rr]^{\iota_\lambda}\ar[d]_{\operatorname{unit}} & & \dG\ar[d]^\pi\\
\Spec(K)\ar[rr]_\lambda & & \dG/\dB.
}
\]
This identification of pullbacks therefore provides a natural isomorphism
\[
\iota_\lambda^\ast\pi^\ast(\msc{F})\ot\iota_\lambda^\ast-\cong \operatorname{unit}^\ast\lambda^\ast(\msc{F})\ot\iota^\ast_\lambda-=\lambda^\ast(\msc{F})\ot_k\iota^\ast_\lambda-.
\]
So in total we have the natural isomorphism
\[
\iota_\lambda^\ast(\msc{F}\star-)\cong \lambda^\ast(\msc{F})\ot_k\iota_\lambda^\ast(-)
\]
at all $\msc{F}$ in $\QCoh(\dG/\dB)$ which provides the pullback functor $\iota_\lambda^\ast:\QCoh(\Xq)\to \msc{B}_\lambda$ with the claimed $\QCoh(\dG/\dB)$-linear structure.
\par

One simply recalls the definition of the restriction maps $\res_\lambda:\FK{G}_q\to \msc{B}_\lambda$ to observe the following.

\begin{lemma}
The composition of universal restriction $\kappa^\ast:\FK{G}_q\to \QCoh(\Xq)$ with the pullback map $\iota_\lambda^\ast:\QCoh(\Xq)\to \msc{B}_\lambda$ recovers $\res_\lambda$,
\[
\res_\lambda\cong \iota^\ast_\lambda\circ \kappa^\ast.
\]
More precisely, the functor $\iota_\lambda^\ast:\QCoh(\Xq)\to \msc{B}_\lambda$ has the natural structure of a morphism between $\QCoh(\dG/\dB)\ot_k\FK{G}_q$-central tensor categories.
\end{lemma}

\subsection{Calculating the fibers}

\begin{theorem}\label{thm:calc_fib}
Consider any geometric point $\lambda:\Spec(K)\to \dG/\dB$.  The $\QCoh(\dG/\dB)$-linear map $\iota^\ast_\lambda:\QCoh(\Xq)\to \msc{B}_\lambda$ induces an equivalence of $K$-linear tensor categories
\[
\QCoh(\Xq)|_\lambda=\Vect(K)\ot_{\QCoh(\dG/\dB)}\QCoh(\Xq)\overset{\sim}\longrightarrow \msc{B}_\lambda.
\]
\end{theorem}

\begin{proof}
Since the pullback map $\iota^\ast_\lambda$ is already monoidal, it suffices to show that the induced map $fib_\lambda$ from the product is an equivalence of linear categories.  We take $\msc{M}=\QCoh(\Xq)$.  Suppose first that $\lambda$ is a closed point, i.e.\ that $K=k$.
\par

For any presentable $k$-linear category $\msc{D}$, a $\QCoh(\dG/\dB)$-bilinear map $F':\Vect\times \msc{M}\to \msc{D}$ is precisely the information of a map $F:\msc{M}\to \msc{D}$ equipped with a natural, associative, isomorphism
\[
F(\msc{F}\star M)\cong (\lambda^\ast \msc{F})\ot_k F(M)
\]
for all $\msc{F}$ in $\QCoh(\dG/\dB)$ and $M$ in $\msc{M}$.  So restricting along the inclusion $\msc{M}\to \Vect\times \msc{M}$, $M\mapsto (k,M)$, provides an equivalence of categories
\[
\operatorname{Bil}_{\QCoh(\dG/\dB)}(\Vect\times \msc{M},\msc{D})\cong\operatorname{Fun}_{\dG/\dB}(\msc{M},\msc{D})
\]
at arbitrary $\msc{D}$.
\par

We have the pushforward map $(\iota_\lambda)_\ast:\msc{B}_\lambda\to \msc{M}$ and the pullback map $\iota_\lambda^\ast:\msc{M}\to \msc{B}_\lambda$.  Restrictions along these maps provide functors
\begin{equation}\label{eq:1639}
\
I^\lambda:\operatorname{Fun}_{\dG/\dB}(\msc{M},\msc{D})\to \operatorname{Fun}_k(\msc{B}_\lambda,\msc{D})\ \ 
\text{and}\ \ I_\lambda:\operatorname{Fun}_k(\msc{B}_\lambda,\msc{D})\to \operatorname{Fun}_{\dG/\dB}(\msc{M},\msc{D}).
\end{equation}
We claim that $I^\lambda$ and $I_\lambda$ are mutually inverse equivalences.  The composite $I^\lambda\circ I_\lambda$ is isomorphic to the identity because $(\iota_\lambda)_\ast\circ (\iota_\lambda)^\ast:\msc{B}_\lambda\to \msc{B}_\lambda$ is isomorphic to the identity.  We consider the composition $I_\lambda\circ I^\lambda$, which is just pulling back functors along the endomorphism $(\iota_\lambda)_\ast\circ (\iota_\lambda)^\ast:\msc{M}\to \msc{M}$.
\par

Let us take $A_\lambda=K(\lambda)$, considered as an algebra object in $\QCoh(\dG/\dB)$.  Then pushing forward along $\lambda:\Spec(k)\to \dG/\dB$ identifies $\Vect$ with the category of arbitrary $A_\lambda$-modules in $\QCoh(\dG/\dB)$, and identifies $\msc{B}_\lambda$ with the category of $A_\lambda$-modules in $\msc{M}$.
\par

We consider the unit of the pullback-pushforward adjunction $M=\O_{\dG/\dB}\star M\to A_\lambda\star M$.  Via $\QCoh(\dG/\dB)$-linearity of $F$, we see that the unit map induces a natural isomorphism
\begin{equation}\label{eq:1649}
F(M)\overset{\cong}\to F(A_\lambda\star M)
\end{equation}
between $F$ and $F(A_\lambda\star-)$ in $\operatorname{Fun}_{\dG/\dB}(\msc{M},\msc{D})$.  But now, there is a natural isomorphism $A_\lambda\star -\cong (\iota_\lambda)_\ast\circ (\iota_\lambda)^\ast-$ via the identification between $\msc{B}_\lambda$ and the category of $A_\lambda$-modules in $\msc{M}$.  So the natural isomorphism \eqref{eq:1649} provides an isomorphism $id_{\operatorname{Fun}}\cong I_\lambda\circ I^\lambda$.  The maps \eqref{eq:1639} are therefore seen to be mutually inverse, as claimed.
\par

For an arbitrary geometric point $\lambda:\Spec(K)\to \dG/\dB$, we note that any bilinear functor $F:\Vect(K)\times \msc{M}\to \msc{D}$ factors through the base change $\msc{D}_K$, which is formally the non-full subcategory of objects $T$ in $\msc{D}$ equipped with a structural action $K\to \End_{\msc{D}}(T)$.  So by base change one can reduce to the argument for closed points.
\end{proof}

\begin{remark}
One can similarly calculate the fiber $\QCoh(S)\ot_{\QCoh(\dG/\dB)}\QCoh(\Xq)$ along any affine morphism $f:S\to \dG/\dB$ (cf.\ \cite[Theorem 3.3]{douglasspsnyder19}).
\end{remark}

Below we denote the equivalence of Theorem \ref{thm:calc_fib} by
\[
fib_\lambda:\QCoh(\Xq)|_\lambda\overset{\sim}\to \msc{B}_\lambda.
\]
The following result validates our interpretation of the embedding $\kappa^\ast:\FK{G}_q\to \QCoh(\Xq)$ as a universal restriction functor for the Borels.

\begin{proposition}\label{prop:univ_res}
The composition
\[
\FK{G}_q\overset{\kappa^\ast}\to \QCoh(\Xq)\overset{reduce}\longrightarrow \QCoh(\Xq)|_\lambda\overset{fib_\lambda}\to \msc{B}_\lambda
\]
is isomorphic to the restriction functor $\res_\lambda:\FK{G}_q\to \msc{B}_\lambda$, as a tensor functor.
\end{proposition}

\begin{proof}
The functor $fib_\lambda$, by definition, satisfies $fib_\lambda\circ reduce\cong \iota^\ast_\lambda$, and one has directly
\[
\begin{array}{l}
\iota^\ast_\lambda\kappa^\ast=\iota^\ast_\lambda\circ forget:\vspace{1mm}\\
\FK{G}_q=\QCoh(\dG)^{G_q}\to \QCoh(\dG)^{B_q}=\QCoh(\Xq)\to \QCoh(\dB_\lambda)^{B_q}=\msc{B}_\lambda.
\end{array}
\]
One consults Section \ref{sect:borels} to see that this map is precisely $\res_\lambda$.
\end{proof}

One can furthermore show that the sequence of Proposition \ref{prop:univ_res} is isomorphic to $\kappa^\ast$ as a \emph{central} tensor functor, though we leave the details to the interested reader.

\subsection{Taking fibers at the derived level}
\label{sect:derived_fib}

We record two versions of Theorem \ref{thm:calc_fib} which hold at a derived level, or more precisely at the level of presentable stable $\infty$-categories.  We consider the derived $\infty$-categories of quasi-coherent sheaves on the half-quantum and classical flag varieties
\[
\QCoh_{\dg}(X_q):=\mcl{D}(\QCoh(\Xq))\ \ \text{and}\ \ \QCoh_{\dg}(\dG/\dB):=\mcl{D}(\QCoh(\dG/\dB)),
\]
as constructed in \cite[Definition 1.3.5.8]{lurieHA}.  We also consider the derived $\infty$-category of sheaves for the small quantum Borel.  Here we recall our explicit geometric construction
\[
\msc{B}_\lambda=\QCoh(\dB_\lambda/B_q)\ \ \text{and take}\ \ \QCoh_{\dg}(\dB_\lambda/B_q):=\mcl{D}(\msc{B}_\lambda),
\]
for the sake of consistency.
\par

At the derived(= stable presentable $\infty$) level we have an action of $\QCoh_{\dg}(\dG/\dB)$ on $\QCoh_{\dg}(\Xq)$, and can again consider the fiber at a given point $\lambda:\Spec(K)\to \dG/\dB$,
\[
\QCoh_{\dg}(\Xq)|_\lambda:=\Vect_{\dg}(K)\ot_{\QCoh_{\dg}(\dG/\dB)}\QCoh_{\dg}(\Xq).
\]
In the above formula $\Vect_{\dg}(K)$ is the derived $\infty$-category of unbounded cochains over $K$, and the base change operation is defined via a universal property \cite[Lemma 4.8.4.3]{lurieHA}, as in the abelian setting.  We have the following derived analog of Theorem \ref{thm:calc_fib}.

\begin{theorem}\label{thm:derived_fibI}
Consider any geometric point $\lambda:\Spec(K)\to \dG/\dB$.  The $\QCoh_{\dg}(\dG/\dB)$-linear map $\operatorname{L}\iota^\ast_\lambda:\QCoh_{\dg}(\Xq)\to \QCoh_{\dg}(\dB_\lambda/B_q)$ induces an equivalence of stable $\infty$-categories
\[
\QCoh_{\dg}(\Xq)|_\lambda\overset{\sim}\longrightarrow \QCoh_{\dg}(\dB_\lambda/B_q).
\]
\end{theorem}

Note that such an equivalence is automatically exact, since it must preserve all limits and colimits.  A proof of Theorem \ref{thm:derived_fibI} is outlined below, following the statement of Theorem \ref{thm:derived_fibII}.
\par

We also have an ``Ind-finite" version of Theorem \ref{thm:derived_fibI}.  Here we consider the Ind-completions of the coherent subcategories $\Coh_{\dg}(Y)\subset \QCoh_{\dg}(Y)$, for $Y$ any of the (noncommutative) spaces considered above.  To be more specific, we take the full, stable, and essentially small $\infty$-subcategory $\Coh_{\dg}(Y)$ in $\QCoh_{\dg}(Y)$ whose objects are all dg sheaves which are bounded and coherent in each degree.  We then apply the Ind-construction \cite[Definition 5.3.5.1, Proposition 5.3.5.10]{lurie09} to obtain $\infty$-categories
\begin{equation}\label{eq:1922}
\IndCoh_{\dg}(\Xq),\ \ \IndCoh_{\dg}(\dG/\dB)=\QCoh_{\dg}(\dG/\dB),\ \ \IndCoh_{\dg}(\dB_\lambda/B_q)
\end{equation}
which are all presentable and stable \cite[Proposition 1.1.3.6]{lurieHA}.
\par

The identification of $\IndCoh_{\dg}(\dG/\dB)$ with $\QCoh_{\dg}(\dG/\dB)$ in \eqref{eq:1922} follows from the fact that all coherent dg sheaves are (already) compact in $\QCoh_{\dg}(\dG/\dB)$ and generate the whole category \cite[Proposition 5.3.5.11]{lurie09}.  Also, the $\infty$-category $\IndCoh_{\dg}(\dB_\lambda/B_q)$ is identified with the Ind-category of the dualizable objects in $\mcl{D}(\msc{B}_\lambda)$.
\par

In this Ind-finite, or Ind-coherent context we again have an action of $\QCoh_{\dg}(\dG/\dB)$ on $\IndCoh_{\dg}(\Xq)$ and can take the fibers along points in the flag variety.  We obtain an alternate version of Theorem \ref{thm:derived_fibI} which, in some instances, is more useful than its quasi-coherent cousin.

\begin{theorem}\label{thm:derived_fibII}
Consider any geometric point $\lambda:\Spec(K)\to \dG/\dB$.  The map $\operatorname{L}\iota^\ast_\lambda:\IndCoh_{\dg}(\Xq)\to \IndCoh_{\dg}(\dB_\lambda/B_q)$ induces an equivalence of stable $\infty$-categories
\[
\IndCoh_{\dg}(\Xq)|_\lambda\overset{\sim}\longrightarrow \IndCoh_{\dg}(\dB_\lambda/B_q).
\]
\end{theorem}

Providing full details for the proofs of these theorems would take us far outside of the (desired) scope of this text.  So we omit formal proofs.  However, we do record the main points here.
\par

For Theorem \ref{thm:derived_fibI}, for example, one can proceed as follows: Consider a flat resolution $\O_\lambda\overset{\sim}\to (\iota_\lambda)_\ast\O_{\dB_\lambda}$ of algebras in $\operatorname{Ch}(\QCoh(\dG_K)^{\dB_K})\subset \operatorname{Ch}(\QCoh(\Xq)_K)$.  We are free to assume additionally that $\O_\lambda$ is bounded and coherent in each degree.  Note that the descent of the aforementioned resolution to $\dG_K/\dB_K$ provides a coherent and flat resolution of the residue field $\overline{\O}_\lambda\overset{\sim}\to K(\lambda)$.
\par

We consider the categories of dg modules $\O_\lambda\text{-mod}_{\Xq}$ and $\overline{\O}_\lambda\text{-mod}_{\dG/\dB}$ in cochains over $\Xq$ and $\dG/\dB$, respectively, and take the corresponding derived $\infty$-categories.\footnote{Since one doesn't have immediate access to a nice model structure on these categories of dg modules, one should be slightly careful in defining their derived $\infty$-categories.  One can sidestep most model theoretic issues, however, by employing stable localization techniques from \cite[\S I.3]{nikolausscholze18} and \cite[\S 5.1, 5.2]{blumberggepnettabuada13}.}  We have obvious maps
\begin{equation}\label{eq:step1}
\mcl{D}(\O_\lambda\text{-mod}_{\Xq})\to \O_\lambda\text{-mod}_{\QCoh_{\dg}(\Xq)}\ \ \text{and}\ \ \mcl{D}(\overline{\O}_\lambda\text{-mod}_{\dG/\dB})\to \overline{\O}_\lambda\text{-mod}_{\QCoh_{\dg}(\dG/\dB)},
\end{equation}
which one shows are equivalences (cf.\ \cite[Theorems 4.1.1 \& 5.2.3]{hinich15} \cite[Corollary 4.2.2.16 \& Theorem 4.3.3.17]{lurieHA}).  On the right hand sides of the above equivalences we consider, specifically, modules for the given algebra objects \emph{in} their respective derived $\infty$-categories \cite[Definition 4.2.1.13]{lurieHA}.
\par

As one expects from the triangular setting \cite[\href{https://stacks.math.columbia.edu/tag/09S6}{Tag 09S6}]{stacks}, base change along the quasi-isomorphisms $\O_\lambda\overset{\sim}\to (\iota_\lambda)_\ast\O_{\dB_\lambda}$ and $\overline{\O}_\lambda\overset{\sim}\to K(\lambda)$ also provide equivalences
\begin{equation}\label{eq:step2}
(\iota_\lambda)_\ast\O_{\dB_\lambda}\ot^{\operatorname{L}}_{\O_\lambda}-:\mcl{D}(\O_\lambda\text{-mod}_{\Xq})\overset{\sim}\longrightarrow\mcl{D}((\iota_\lambda)_\ast\O_{\dB_\lambda}\text{-mod}_{\Xq})=\QCoh_{\dg}(\dB_\lambda/B_q)
\end{equation}
and
\[
K(\lambda)\ot^{\operatorname{L}}_{\overline{\O}_\lambda}-:\mcl{D}(\overline{\O}_\lambda\text{-mod}_{\dG/\dB})\overset{\sim}\longrightarrow \mcl{D}(K(\lambda)\text{-mod}_{\dG/\dB})=\Vect_{\dg}(K).
\]
So we have an equivalence $\Vect_{\dg}(K)\cong \overline{\O}_\lambda\text{-mod}_{\QCoh_{\dg}(\dG/\dB)}$, and the generic base change formula \cite[Theorem 4.8.4.6]{lurieHA} identifies the fiber in question as
\[
\QCoh_{\dg}(\Xq)|_\lambda=\overline{\O}_\lambda\text{-mod}_{\QCoh_{\dg}(\Xq)},
\]
i.e.\ as the category of $\overline{\O}_\lambda$-modules in $\QCoh_{\dg}(\Xq)$ under the derived $\star$-action of $\QCoh_{\dg}(\dG/\dB)$ on $\QCoh_{\dg}(\Xq)$.  However, since $\QCoh_{\dg}(\dG/\dB)$ acts via the tensor functor $\operatorname{L}\zeta^\ast:\QCoh_{\dg}(\dG/\dB)\to \QCoh_{\dg}(\Xq)$, and $\operatorname{L}\zeta^\ast\overline{\O}_\lambda=\O_\lambda$, we have a further identification
\begin{equation}\label{eq:step3}
\QCoh_{\dg}(\Xq)|_\lambda=\overline{\O}_\lambda\text{-mod}_{\QCoh_{\dg}(\Xq)}=\O_\lambda\text{-mod}_{\QCoh_{\dg}(\Xq)}.
\end{equation}
We consult the above formulas (\ref{eq:step1}, \ref{eq:step2}, \ref{eq:step3}) to obtain the proposed calculation of the fiber $\QCoh_{\dg}(\Xq)|_\lambda\overset{\sim}\to \QCoh_{\dg}(\dB_\lambda/B_q)$.
\par

As a final point, one can show that every equivalence above appears in a diagram over $\QCoh_{\dg}(\Xq)$, so that we have a large diagram
\[
\xymatrix{
 & \QCoh_{\dg}(\Xq)\ar@{..>}[r]^{\operatorname{univ}}\ar[ddr]^{\overline{\O}_\lambda\star-}\ar[d]|{\O_\lambda\ot-}\ar[ddl]_{\iota_\lambda^\ast} & \QCoh_{\dg}(\Xq)|_\lambda\ar@{..>}[dd]^=\\
 & \mcl{D}(\O_\lambda\text{-mod}_{\Xq})\ar[dr]_\sim\ar[dl]^\sim \\
\QCoh_{\dg}(\dB_\lambda/B_q) & & \overline{\O}_\lambda\text{-mod}_{\QCoh_{\dg}(\Xq)}\ar@{-->}[ll]_\sim^{\rm induced}. \\
}
\]
The above diagram tells us that the equivalence
\[
\QCoh_{\dg}(\Xq)|_\lambda\overset{\sim}\longrightarrow \QCoh_{\dg}(\dB_\lambda/B_q)
\]
obtained above is in fact induced by the derived pullback functor $\operatorname{L}\iota_\lambda^\ast$.
\par

The proof in the Ind-coherent setting is similar.  However, one has to take care to keep track of coherent objects at every step.  Dealing with the derived categories of dg modules, in particular, becomes subtle at that point.

\begin{remark}
In Section \ref{sect:fibers_sHom} we provide a calculation of fibers for the derived sheaf-Hom functor which occurs at a strictly triangular level, and directly reflects the fiber calculation of Theorem \ref{thm:derived_fibII}.
\end{remark}

\section{$\sRHom$ as a total space for morphisms over the Borels}
\label{sect:fibers_sHom}

In the previous section we calculated the fibers of the family of categories $\QCoh(\Xq)$ over $\dG/\dB$.  In particular, we saw that pulling back along the map $\iota_\lambda:\dB_\lambda\to \dG$ provides a calculation of the fiber
\[
\QCoh(\Xq)|_\lambda\overset{\sim}\to \msc{B}_\lambda
\]
at any geometric point for $\dG/\dB$, and that a similar calculation holds at the derived level.
\par

In this section we provide a parallel calculation of the fibers of the dg sheaves $\sRHom_{\dG/\dB}(M,N)$ in $D(\dG/\dB)$ via morphisms over $\msc{B}_\lambda$.  If one takes proper account of naturality, such a calculation of the fibers for $\sRHom_{\dG/\dB}$ equivalently calculates the fibers $\lambda^\ast D^{\operatorname{Enh}}(\Xq)$ of the enhanced derived category over $\dG/\dB$.  As in Section \ref{sect:fibers_G/B}, we employ $\QCoh(\dG/\dB)$-linearity of the pullback functor $\iota^\ast_\lambda:\QCoh(\Xq)\to \msc{B}_\lambda$ in an essential way in our analysis, and show that such linearity implies the desired calculation of the fibers for $\sRHom_{\dG/\dB}$.

\subsection{The natural map $\lambda^\ast\sHom_{\Xq}\to \Hom_{\msc{B}_\lambda}$}
\label{sect:1226}

Consider $M$ and $N$ in $\QCoh(\Xq)$.  As observed in Section \ref{sect:QCoh-linear}, the monoidal functor $\iota^\ast_\lambda:\QCoh(\Xq)\to \msc{B}_\lambda$ admits a natural $\QCoh(\dG/\dB)$-linear structure.  We can therefore apply $\iota^\ast_\lambda$ to evaluation
\[
ev:\sHom_{\Xq}(M,N)\star M\to N
\]
to obtain a map
\begin{equation}\label{eq:1705}
\iota^\ast_\lambda ev:\lambda^\ast\sHom_{\Xq}(M,N)\ot_k \iota_\lambda^\ast M\to\iota_\lambda^\ast N\ \ \text{in the category }\msc{B}_\lambda.
\end{equation}
By adjunction we then obtain a natural map
\begin{equation}\label{eq:1352}
\phi_{M,N}=\phi(\lambda)_{M,N}:\lambda^\ast\sHom_{\Xq}(M,N)\to \Hom_{\msc{B}_\lambda}(\iota^\ast_\lambda M,\iota^\ast_\lambda N)
\end{equation}
which is compatible with the evaluation, and hence compatible with composition and the tensor structure.  This is to say, the maps $\phi_{M,N}$ collectively provide a linear monoidal functor
\[
\phi(\lambda):\lambda^\ast\QCoh^{\Enh}(\Xq)\to \msc{B}_\lambda.
\]

\begin{proposition}\label{prop:1356}
The map $\phi_{M,N}$ of \eqref{eq:1352} is an isomorphism whenever $M$ is coherent and relatively projective.
\end{proposition}

We delay the proof to the end of the section, and focus instead on the implications of Proposition \ref{prop:1356} to our analysis of the derived inner-$\Hom$ functor.

\subsection{The derived map $\operatorname{L}\phi(\lambda):\operatorname{L}\lambda^\ast D^{\Enh}(\Xq)\to D^{\Enh}(\msc{B}_\lambda)$}

Let us begin with a basic lemma.

\begin{lemma}\label{lem:1399}
\begin{enumerate}
\item If $M$ is relatively projective in $\QCoh(\Xq)$ and $N$ is flat, then the sheaf $\sHom_{\Xq}(M,N)$ is flat over $\dG/\dB$.
\item For any sheaf $\msc{F}$ over $\dG/\dB$, complex $P$ of relatively projective sheaves over $\Xq$, and bounded complex $N$ of flat sheaves over $\Xq$, the natural map
\[
\msc{F}\ot^{\operatorname{L}}\sHom_{\Xq}(P,N)\to \msc{F}\ot\sHom_{\Xq}(P,N)
\]
is a quasi-isomorphism.
\item For any closed subscheme $i:Z\to \dG/\dB$, and $P$ and $N$ as in {\rm (2)}, the natural map
\[
\operatorname{L}i^\ast\sHom_{\Xq}(P,N)\to i^\ast\sHom_{\Xq}(P,N)
\]
is a quasi-isomorphism.
\end{enumerate}
\end{lemma}

\begin{proof}
Let us note, before beginning, that flatness of $N$ implies exactness of the operation $-\star N=\pi^\ast(-)\ot N$.  For the first point, consider an exact sequence $0\to \msc{F}'\to \msc{F}\to \msc{F}''\to 0$ of coherent sheaves and the corresponding possibly (non-)exact sequence
\[
0\to \msc{F}'\ot \sHom(M,N)\to \msc{F}\ot \sHom(M,N)\to \msc{F}''\ot \sHom(M,N)\to 0.
\]
By $\Coh(\dG/\dB)$-linearity of $\sHom(M,-)$, Lemma \ref{lem:sHom-linear}, the above sequence is isomorphic to the sequence
\[
0\to \sHom(M,\msc{F}'\star N)\to \sHom(M,\msc{F}\star N)\to \sHom(M,\msc{F}''\star N)\to 0.
\]
The second sequence is exact by flatness of $N$ and local projectivity of $M$.  So we see that $\sHom(M,N)$ is flat relative to the action of coherent sheaves on $\dG/\dB$.  This is sufficient to see that $\sHom(M,N)$ is flat in $\QCoh(\dG/\dB)$.
\par

For the second point, resolve $\msc{F}$ by a finite complex of flat sheaves $\msc{F}'\to \msc{F}$.  Via a spectral sequence argument, using flatness of $\sHom(P,N)$ in each degree, one sees that the induced map
\[
\msc{F}\ot^{\operatorname{L}}\sHom(P,N)=\msc{F}'\ot\sHom(P,N)\to \msc{F}\ot\sHom(P,N)
\]
is a quasi-isomorphism.  Point (3) follows from point (2) and the identification $\operatorname{L}i^\ast(-)=i_\ast\O_Z\ot^{\operatorname{L}}-$.
\end{proof}

We now consider the derived category $D(\Xq)$.  We have the derived pullback
\[
\operatorname{L}\iota^\ast_\lambda:D(\Xq)\to D(\msc{B}_\lambda)
\]
which still annihilates the $D(\dG/\dB)$-action.  So, as in Section \ref{sect:1226}, we get an induced map on inner-Homs
\[
\operatorname{L}\phi_{M,N}:\operatorname{L}\lambda^\ast\sRHom_{\Xq}(M,N)\to \RHom_{\msc{B}_\lambda}(\operatorname{L}\iota^\ast_\lambda M,\operatorname{L}\iota^\ast_\lambda N)
\]
which is compatible with composition and the tensor structure.  We therefore obtain a monoidal functor
\begin{equation}\label{eq:1470}
\operatorname{L}\phi(\lambda):\operatorname{L}\lambda^\ast D^{\Enh}(\Xq)\to D^{\Enh}(\msc{B}_\lambda),
\end{equation}
where $D^{\Enh}(\msc{B}_\lambda)$ is the linear enhancement of $D(\msc{B}_\lambda)$ implied by the action of $D(\Vect)$.
\par

Consider coherent $M$ and bounded $N$ in $D(\Xq)$.  (By coherent we mean that $M$ is in $D_{\coh}(\Xq)$.)   If we express $\sRHom_{\Xq}$ by resolving the first coordinate by relatively projective sheaves, which are necessarily flat, and we replace $N$ with a bounded complex of flat sheaves if necessary, the map $\operatorname{L}\phi_{M,N}$ is simply the fiber map
\[
\operatorname{L}\phi_{M,N}=\phi_{M,N}:\lambda^\ast\sHom_{\Xq}(P_M,N)\to \Hom_{\msc{B}_\lambda}(\iota_\lambda^\ast P_M,\iota_\lambda^\ast N)
\]
defined at equation \eqref{eq:1352}, via Lemma \ref{lem:1399} (3).  By Proposition \ref{prop:1356}, the map $\operatorname{L}\phi_{M,N}$ is then seen to be an isomorphism whenever $M$ is coherent and $N$ is bounded.

\begin{theorem}\label{thm:fiber_l}
The map
\[
\operatorname{L}\phi_{M,N}:\operatorname{L}\lambda^\ast\sRHom_{\Xq}(M,N)\to \RHom_{\msc{B}_\lambda}(\operatorname{L}\iota^\ast_\lambda M,\operatorname{L}\iota^\ast_\lambda N)
\]
is a quasi-isomorphism whenever $M$ is coherent and $N$ is bounded.  Consequently, the monoidal functor
\[
\operatorname{L}\phi(\lambda):\operatorname{L}\lambda^\ast D^{\Enh}(\Xq)\to D^{\Enh}(\msc{B}_\lambda)
\]
is fully faithful when restricted to the full subcategory of coherent dg sheaves.
\end{theorem}

\begin{proof}
We have already argued above that $\operatorname{L}\phi_{M,N}$ is an isomorphism at such $M$ and $N$.  Fully faithfulness of the restriction of $\operatorname{L}\phi(\lambda)$ to $D^{\Enh}_{coh}(\Xq)$ follows.
\end{proof}

\begin{remark}
From the perspective of Theorem \ref{thm:derived_fibII}, the finiteness conditions appearing in statement of Theorem \ref{thm:fiber_l} are relatively intuitive.  Namely, the fiber map $\IndCoh_{\dg}(\Xq)|_\lambda\overset{\sim}\to \IndCoh_{\dg}(\dB_\lambda/B_q)$ is basically given by base change $K(\lambda)\star-$, and the natural map
\[
K(\lambda)\ot\sRHom_{\Xq}(M,N)\to \sRHom_{\Xq}(M,K(\lambda)\star N)=\lambda_\ast\RHom_{\msc{B}_\lambda}(M|_\lambda, N|_\lambda)
\]
should be an isomorphism whenever $M$ is coherent, and hence compact in the Ind-category.  The restrictions on $N$ in Theorem \ref{thm:derived_fibII} then appear because $\QCoh_{\dg}(\Xq)$ and $\IndCoh_{\dg}(\Xq)$ differ in general, but intersect at the bounded derived $\infty$-category of quasi-coherent sheaves.  In particular, the functor $\IndCoh_{\dg}(\Xq)\to\QCoh_{\dg}(\Xq)$ admits a fully faithful section over the subcategory of bounded quasi-coherent dg sheaves.
\end{remark}

\subsection{Proof of Proposition \ref{prop:1356}}

\begin{proof}[Proof of Proposition \ref{prop:1356}]
Recall the explicit expression of $\sHom_{\Xq}$ in terms of morphisms in $\QCoh(Y_q)=\QCoh(\dG)^{\uqB}$, provided by Proposition \ref{prop:shom_exp}.  After pulling back $\pi^\ast:\QCoh(\dG/\dB)\overset{\cong}\to \QCoh(\dG)^{\dB}$ we, equivalently, have a morphism
\begin{equation}\label{eq:1723}
(\O(\dB_\lambda)\ot_{\O(\dG)}\Hom_{Y_q}(M,N))^{\dB}\to \Hom_{\msc{B}_\lambda}(\iota^\ast_\lambda M,\iota^\ast_\lambda N)
\end{equation}
which is compatible with evaluation, and we are claiming that this map is an isomorphism.  By Lemma \ref{lem:sHom-linear}, or rather the proof of Lemma \ref{lem:sHom-linear}, and local projectivity of $M$ the map
\[
\O(\dB_\lambda)\ot_{\O(\dG)}\Hom_{Y_q}(M,N)\to \Hom_{Y_q}(M,\iota^\ast_\lambda N),\ \ f\ot \xi\mapsto f\cdot \operatorname{red}_\lambda \xi
\]
is an isomorphism, where $\operatorname{red}_\lambda:N\to \iota^\ast_\lambda N$ is the reduction map.  (Here we are viewing sheaves on $\dB_\lambda$ as sheaves on $\dG$ via pushforward.)  Furthermore $\Hom_{Y_q}(M,\iota^\ast_\lambda N)=\Hom_{Y_q}(\iota^\ast_\lambda M,\iota^\ast_\lambda N)$ and when we take invariants we have
\[
\Hom_{Y_q}(\iota^\ast_\lambda M,\iota^\ast_\lambda N)^{\dB}=\Hom_{\Coh(\Xq)}(\iota^\ast_\lambda M,\iota^\ast_\lambda N).
\]
We therefore have an isomorphism
\begin{equation}\label{eq:1731}
(\O(\dB_\lambda)\ot_{\O(\dG)}\Hom_{Y_q}(M,N))^{\dB}\overset{\cong}\to \Hom_{\Xq}(\iota^\ast_\lambda M,\iota^\ast_\lambda N)
\end{equation}
given by $f\ot \xi\mapsto f\cdot \operatorname{red}_\lambda(\xi)$, and under this isomorphism the reduction of evaluation for $\Hom_{Y_q}(M,N)$ appears as the expected evaluation map
\[
\Hom_{\Xq}(\iota^\ast_\lambda M,\iota^\ast_\lambda N)\ot_k \iota^\ast_\lambda M\to \iota^\ast_\lambda N,\ \ \xi\ot m\mapsto \xi(m).
\]
\par

Under the identification \eqref{eq:1731} the map \eqref{eq:1723} now appears as
\begin{equation}\label{eq:1742}
\Hom_{\Xq}(\iota^\ast_\lambda M,\iota^\ast_\lambda N)\to \Hom_{\msc{B}_\lambda}(\iota_\lambda^\ast M,\iota^\ast_\lambda N),\ \ \xi\mapsto \xi.
\end{equation}
One simply observes that these morphism spaces are literally equal, since $\msc{B}_\lambda=\Coh(\dB_\lambda)^{B_q}$ and $\Coh(\Xq)=\Coh(\dG)^{B_q}$, and notes the map \eqref{eq:1742} is the identity to see that \eqref{eq:1742} is an isomorphism.  It follows that our original map \eqref{eq:1723} is an isomorphism.
\end{proof}

\section{Accessing the Springer resolution via $\sRHom$}
\label{sect:RGq}

We now turn away from our discussion of the triad, consisting of the half-quantum flag variety, small quantum group, and small quantum Borel, and observe a connection between the Springer resolution and sheaves on $\Xq$.
\par

We show that the sheafy endomorphism algebra
\[
\msc{A}_{\tN}:=\sRHom_{\Xq}(\1,\1)
\]
for the unit in $\QCoh(\Xq)$, which is formally an algebra in the derived category $D(\dG/\dB)$, has cohomology equal to the structure sheaf for the Springer resolution.  This \emph{suggests} a further enhancement $D^{\operatorname{Enh}^2}(\Xq)$ in the derived category of (dg) sheaves over the Springer resolution.  We elaborate on the latter point in the subsequent, and final, section of the text.
\par

We begin by recalling the necessary background concerning the Springer resolution.  We then provide the calculation of cohomology
\[
H^\ast(\msc{A}_{\tN})=p_\ast\O_{\tN}
\]
at sufficiently large odd order parameter $q$, which we deduced from results of Ginzburg and Kumar \cite{ginzburgkumar93}.

\subsection{The Springer resolution}

Let $\mfk g$ be the Lie algebra for $\dG$. Recall that the Springer resolution $\tN$ is the affine bundle $\tN=\dG\times_{\dB}\mfk{n}$ over the flag variety $\dG/\dB$, where $\mathfrak{n}$ is the (positive) nilpotent subalgebra $\mfk{n}\subset \mfk{g}$.  Equivalently, the Springer resolution is obtained as the relative spectrum \cite[\S 9.1]{EGAI} of the descent of the $\dB$-equivariant algebra $\O_{\dG}\ot \Sym(\mfk{n}^\ast)$ over $\dG$,
\begin{equation}\label{eq:1932}
\tN=\dG\times_{\dB}\mfk{n}=\Spec_{\dG/\dB}\left(\text{descent\ of }\O_{\dG}\ot_k \Sym(\mfk{n}^\ast)\right)=\Spec_{\dG/\dB}(\Sym(\msc{E})).
\end{equation}
In the above formula $\msc{E}$ is the equivariant vector bundle on $\dG/\dB$ associated to the $\dB$-representation $\mfk{n}^\ast$.
\par

From this construction of $\tN$ as the relative spectrum of a sheaf of algebras on the flag variety, we see that pushing forward along the (affine) structure map $p:\tN\to \dG/\dB$ provides an identification $p_\ast\O_{\tN}=\Sym(\msc{E})$ and also an equivalence of monoidal categories
\[
p_\ast:\QCoh(\tN)\overset{\sim}\longrightarrow \QCoh(p_\ast\O_{\tN})
\]
\cite[Th\'{e}or\`{e}me 9.2.1]{EGAI}.  To be clear, the latter category is the category of modules over the commutative algebra object $p_\ast\O_{\tN}$ in $\QCoh(\dG/\dB)$, and the monoidal product is as expected $\ot_{p_\ast\O_{\tN}}$.
\par

The Springer resolution $\tN$ can alternatively be identified with the moduli space of choices of a Borel in $\mfk{g}$, and a nilpotent element in the given Borel,
\begin{equation}\label{eq:moduli}
\mcl{M}oduli=\{(\mfk{b}_\lambda,x):\mfk{b}_\lambda\subset\mfk{g}\ \text{a Borel}\ x\in \mfk{b}_\lambda\ \text{is nilpotent}\}
\end{equation}
\cite[\S 3.2]{chrissginzburg09}.  This moduli space sits in the product $\mcl{M}oduli\subset \dG/\dB\times \mcl{N}$, where $\mcl{N}$ is the nilpotent cone in $\mfk{g}$, and we have an explicit isomorphism between $\tN$ and $\mcl{M}oduli$ given by the $\dG$-actions on the two factors in this product,
\[
\tN=\dG\times_{\dB} \mfk{n}\overset{\cong}\to \mcl{M}oduli,\ \ (g,x)\mapsto (\operatorname{Ad}_g(\mfk{b}),\operatorname{Ad}_g(x)).
\]
In the above formula $\mfk{b}$ is the positive Borel in $\mfk{g}$.  We identify $\tN$ with this moduli space when convenient, via the above isomorphism.

We consider $\tN$ as a conical variety over $\dG/\dB$ by taking the generating bundle $\msc{E}$ to be in (cohomological) degree $2$.  In terms of the moduli description given above, this conical structure corresponds to a $\mbb{G}_m$-action defined by the squared scaling $c\cdot (\mfk{b}_\lambda, x)=(\mfk{b}_\lambda, c^2\cdot x)$.

\begin{remark}\label{rem:dg}
For certain applications $\tN$ might be viewed, more fundamentally, as a dg scheme over $\dG/\dB$ which has generating bundle in cohomological degree $2$ and vanishing differential.  Compare with \cite{arkhipovbezrukavnikovginzburg04}.
\end{remark}

\subsection{The moment map}

Given our identification of the Springer resolution $\tN$ with the moduli of pairs \eqref{eq:moduli}, we have two projections $\mfk{b}_\lambda\leftarrow (\mfk{b}_\lambda,x)\to x$ which define maps $p:\tN\to \dG/\dB$ and $\mu:\tN\to \mcl{N}$.  The map $p$ is affine, and simply recovers the structural map $\tN=\Spec_{\dG/\dB}(\Sym(\msc{E}))\to \dG/\dB$.  The map $\mu$ provides an identification of the affinization of $\tN$ with the nilpotent cone,
\[
\bar{\mu}:\tN_{\rm aff}\overset{\cong}\longrightarrow \mcl{N}.
\]
The map $\mu$ is called the \emph{moment map}.  The moment map is a proper birational equivalence, and so realizes the Springer resolution as a resolution of singularities for the nilpotent cone \cite[Theorem 10.3.8]{hottatanisaki07}.

\subsection{A calculation of cohomology}
\label{sect:GK}

The following is deduced from results of Ginzburg and Kumar \cite{ginzburgkumar93}.

\begin{theorem}\label{thm:enh_GK}
Suppose that $q$ is of odd order $\ord(q)>h$, or that $G$ is of type $A_1$.  There is a canonical identification
\[
H^\ast(\msc{A}_{\tN})=p_\ast\O_{\tN},
\]
as sheaves of graded algebras over $\dG/\dB$.
\end{theorem}

\begin{proof}
Recall that $\QCoh(\dG/u)=\QCoh(\dG)^{\uqB}$, by definition.  We have
\[
\sRHom_{\dG/u}(\1,\1)=\O_{\dG}\ot_k \RHom_{\uqB}(k,k),
\]
where $\RHom_{\uqB}(k,k)$ is given its natural $\dB$-action as inner morphisms for the $\Rep \dB$-action on $\Rep B_q$.  Hence, by the calculation of $\Ext_{\uqB}(k,k)$ provided in \cite[Lemma 2.6]{ginzburgkumar93}, we have
\[
H^\ast(\sRHom_{\dG/u}(\1,\1))=\O_{\dG}\ot_k \Ext_{\uqB}(k,k)=\O_{\dG}\ot_k \Sym(\mfk{n}^\ast).
\]
One therefore applies Proposition \ref{prop:shom_exp} to obtain
\[
H^\ast(\msc{A}_{\mcl{N}})=H^\ast(\sRHom_{\Xq}(\1,\1))=\text{descent of }\O_{\dG}\ot_k \Sym(\mfk{n}^\ast)=p_\ast\O_{\tN}.
\]
\par

The proof in the type $A_1$ setting, at even order $q$, is completely similar.  Here one calculates the cohomology directly $\Ext_{\uqB}(k,k)=k[x]$, where $x$ is of degree $2$ and the $\dB$-action is determined by the associated character for the torus.  We have $kx=\Ext_{\uqB}^2(k,k)=(kE^l)^\ast\cong \mfk{n}^\ast$, giving again $\Ext_{\uqB}(k,k)=\Sym(\mfk{n}^\ast)$.
\end{proof}

As noted in the original work \cite{ginzburgkumar93}, the higher global section $H^{>0}(\dG/\dB,p_\ast\O_{\tN})$ vanish so that Theorem \ref{thm:enh_GK} implies a computation of extensions for $\QCoh(\Xq)$.  In the following statement $\O(\mcl{N})$ is considered as a cohomologically graded algebra with generators in degree $2$.

\begin{corollary}[\cite{ginzburgkumar93}]\label{cor:GK}
Suppose $q$ is of odd order $\ord(q)>h$, or that $G$ is of type $A_1$.  There is an identification of graded algebras $\Ext_{\Xq}(\1,\1)=\O(\mcl{N})$.
\end{corollary}

\begin{remark}
The calculation $H^\ast(\msc{A}_{\tN})=p_\ast\O_{\tN}$ is expected to hold at arbitrary $G$ and arbitrary (large order) $q$.  In particular, whether $q$ is of odd order or even order shouldn't matter.  In order to determine this cohomology one need only calculate the specific $\dB$-representation $\mfk{e}=\Ext^2_{\uqB}(k,k)$ at even order $q$, as in general the cohomology $H^\ast(\msc{A}_{\tN})$ is identified with the pushforward of the structure sheaf on the affine bundle $\dG\times_{\dB}\mfk{e}^\ast$.  Such a calculation has not appeared in the literature however.
\end{remark}

\section{Formality conjectures and geometric representation theory}
\label{sect:GRT}

As a final point in the paper, we conjecture that the derived endomorphism algebra $\msc{A}_{\tN}=\RHom_{\Xq}(\1,\1)$ is formal, and that this formality lifts to a categorical level.  In particular, we conjecture that the derived ($\infty$-)category of sheaves on the half-quantum flag variety forms a sheaf of tensor categories over the Springer resolution.  We provide related conjectures, which claim that certain non-monoidal equivalences of Arkhipov-Bezrukavnikov-Ginzburg and Bezrukavnikov-Lachowska \cite{arkhipovbezrukavnikovginzburg04,bezrukavnikovlachowska07} are recoverable in a manner which is manifestly monoidal.
\par

We first record our conjectures, then provide a rationalle for these claims, based on the findings of this text and preexisting results from \cite{arkhipovbezrukavnikovginzburg04,bezrukavnikovlachowska07}.  In this section we exchange our geometric interpretation for quantum group representations $\FK{G}_q$ for the standard algebraic construction $\Rep(\uqG)$ (see Theorem \ref{thm:ag}).

\subsection{Main conjectures}

As in Section \ref{sect:derived_fib}, we let $\QCoh_{\dg}(Y)$ denote the derived $\infty$-category of sheaves on a (possibly noncommutative) space $Y$.  For an algebraic group, or quantum group $T$, we take
\begin{equation}\label{eq:2317}
\Rep_{\dg}(T):=\Ind \mcl{D}_{fin}(T),
\end{equation}
where $\mcl{D}_{fin}(T)$ is the derived $\infty$-category which sits over the usual derived category of (bounded) finite-dimensional dg representations for $T$.  The $\Ind$-construction is as in \cite[Definition 5.3.5.1]{lurie09}.
\par

We also consider the Ind-category $\IndCoh_{\dg}(Y)$, where $\Coh_{\dg}(Y)$ is the derived $\infty$-category of coherent dg sheaves.  When $Y$ is smooth dg scheme, or smooth and sufficiently tame stack, we have
\[
\IndCoh_{\dg}(Y)=\QCoh_{\dg}(Y)
\]
\cite[Proposition 5.3.5.11]{lurie09}.

\begin{conjecture}[Strong Formality]\label{conj:sfc}
There is a central, fully faithful, $\QCoh_{\dg}(\dG/\dB)$-linear tensor functor
\[
\eta^\ast:\QCoh_{\dg}(\tN)\to \IndCoh_{\dg}(\Xq).
\]
The functor $\eta^\ast$ is an equivalence onto the localizing $\QCoh_{\dg}(\dG/\dB)$-submodule category generated by the unit in $\IndCoh_{\dg}(\Xq)$.  Similarly, there is a fully faithful tensor functor $\bar{\eta}^\ast:\QCoh_{\dg}(\mcl{N})\to \Rep_{\dg}(\uqG)$ which is an equivalence onto the localizing subcategory generated by the unit in $\Rep_{\dg}(\uqG)$.
\end{conjecture}

We note that the claim for the quantum group would follow from the corresponding claim for $\Xq$, via fully faithfulness of the embedding $\kappa^\ast:\Rep_{\dg}(\uqG)\to \IndCoh_{\dg}(\Xq)$ (Theorem \ref{thm:Kempf}) and the fact that pullback identifies $\QCoh_{\dg}(\mcl{N})$ with the localizing subcategory generated by the unit in $\QCoh_{\dg}(\tN)$.  So the existence of the functor $\eta^\ast$ for $\Xq$ is the essential point here.
\par

At a global level, one can approach Conjecture \ref{conj:sfc} by first recognizing a non-monoidal embedding
\[
\QCoh_{\dg}(\mfk{n}/\dB)\to \Rep_{\dg}(B_q),
\]
which one more-or-less extracts from \cite{arkhipovbezrukavnikovginzburg04}, then by ``proceeding towards" the Springer resolution and $\Xq$ via a sequence of categorical base change opperations.  In this way one can reduce Conjecture \ref{conj:sfc} to an analysis of the nilpotent subalgebra $\mfk{n}$ and its relations to the quantum Borel.  At this point, however, we leave any additional details on this line of inquiry to a later text.

We note that such a functor $\eta^\ast$ would realize sheaves on the half-quantum flag variety as a tensorial \emph{correspondence} between the small quantum group and the Springer resolution
\begin{equation}\label{eq:2164}
\xymatrix{
 & \IndCoh_{\dg}(\Xq)\\
\QCoh_{\dg}(\tN)\ar[ur]^{\eta^\ast} & & \Rep_{\dg}(u(G_q))\ar[ul]_{\kappa^\ast},
}
\end{equation}
where by a correspondence we mean a pair of (ideally mutually centralizing) monoidal functors to a common target.  Such a pair \eqref{eq:2164} might also be thought of as a morphism $\QCoh_{\dg}(\tN)\to \Rep_{\dg}(u(G_q))$ in a $4$-category of braided tensor $\infty$-categories over $\Vect_{\dg}$ \cite[Example 1.14]{johnsonscheimbauer17}.
\par

As with any correspondence, our Pavlovian response is to consider the associated push-pull functor.  We note that although this correspondence is manifestly monoidal in nature, the associated push-pull functors will not be monoidal, since pushforward functors are generally non-monoidal.  In the following conjecture we consider the principal block $\operatorname{PrinBlock}_{\dg}(u(G_q))$ in $\Rep_{\dg}(\uqG)$.

\begin{conjecture}\label{conj:grt}
The push-pull functor $\eta_\ast\kappa^\ast$ restricts to an equivalence
\[
\eta_\ast\kappa^\ast:\operatorname{PrinBlock}_{\dg}(u(G_q))\overset{\sim}\to \QCoh_{\dg}(\tN).
\]
Furthermore, this equivalence recovers the equivalence of Arkhipov-Bezrukavnikov-Ginzburg and Bezrukavnikov-Lachowska \cite{arkhipovbezrukavnikovginzburg04,bezrukavnikovlachowska07}.
\end{conjecture}

The point here is that the non-monoidal analyses of the texts \cite{arkhipovbezrukavnikovginzburg04,bezrukavnikovlachowska07} should be recoverable in a manner which clearly takes the monoidal structures on $\QCoh_{\dg}(\tN)$ and $\Rep_{\dg}(\uqG)$ into account.  We close the paper with some elaborations on Conjecture \ref{conj:grt}.

\subsection{Additional comments for Conjecture \ref{conj:grt}}

The rationalle behind this second conjecture is rather simple minded.  We forgo a direct comparison with Bezrukavnikov-Lachowska \cite{bezrukavnikovlachowska07}, and deal with the claim that the given map is an equivalence via Arkhipov-Bezrukavnikov-Ginzburg \cite{arkhipovbezrukavnikovginzburg04}.
\par

Let us suppose that Conjecture \ref{conj:sfc} is in fact valid.\footnote{One only needs a non-monoidal embedding $\QCoh_{\dg}(\tN)\to \QCoh_{\dg}(\Xq)$ here.  Having forgone monoidality, such a nice functor almost certainly exists.}  Via a calculus of equivariantization and de-equivariantization, and the fact that the map $\eta_\ast\kappa^\ast$ is (or rather, will be) $\dG$-equivariant, we have a corresponding functor $F:\operatorname{PrinBlock}_{\dg}(G_q)\to \QCoh_{\dg}(\tN)^{\dG}$ which fits into a diagram
\[
\xymatrix{
\operatorname{PrinBlock}_{\dg}(G_q)\ar[rr]^F\ar[d]& & \QCoh_{\dg}(\tN)^{\dG}\ar[d]\\
\operatorname{PrinBlock}_{\dg}(u(G_q))\ar[rr]^{\eta_\ast\kappa^\ast} & & \QCoh_{\dg}(\tN)
}
\]
\cite[Theorem 5.8]{arkhipovgaitsgory03}.  This calculus furthermore implies that $F$ is an equivalence if and only if $\eta_\ast\kappa^\ast$ is an equivalence (cf.\ \cite[Proposition 8.6]{negron21}), and via faithfulness of the forgetful functor $\QCoh_{\dg}(\tN)^{\dG}\to \QCoh_{\dg}(\tN)$ the functor $F$ is completely determined by its composite to $\QCoh_{\dg}(\tN)$.  The above diagram identifies this composite as
\[
\operatorname{PrinBlock}_{\dg}(G_q)\to \QCoh_{\dg}(\tN),\ \ V\mapsto \eta_\ast\kappa^\ast(V),
\]
and
\[
\eta_\ast\kappa^\ast(V)=\sRHom_{\Xq}(\1,E_V)=\text{descent of }\O_{\dG}\ot_k\RHom_{u(B_q)}(k,V)^\sim.
\]
Here we understand $\RHom_{u(B_q)}(k,V)$ as a $\dB$-equivariant $\O(\mfk{n})$-module (see \cite[Section 1.4]{arkhipovbezrukavnikovginzburg04}).  So we identify $F$ as the functor
\[
F:V\mapsto \text{descent of }\O_{\dG}\ot_k\RHom_{u(B_q)}(k,V),
\]
where $\dG$-equivariance is deduced via $\dG$-equivariance of $\O_{\dG}$.  Now, upon considering the discussions around \cite[Equation 1.4.1]{arkhipovbezrukavnikovginzburg04}, this appears to be precisely the equivalence from \cite{arkhipovbezrukavnikovginzburg04}.

\appendix

\section{Inner-Hom nonsense}
\label{sect:A}

\subsection{Proof outline for Proposition \ref{prop:enriched}}

\begin{proof}[Outline for Proposition \ref{prop:enriched}]
If we take $\msc{F}=\sHom_{\Xq}(M,N)$, composition is alternatively defined by the $\QCoh(\dG/\dB)$-linearity $\msc{F}\star \sHom(L,M)\to \sHom(L,\msc{F}\star M)$ composed with evaluation $ev:\msc{F}\star M\to N$ in the second coordinate.  One uses these two descriptions to deduce associativity of composition.
\par

Associativity for the monoidal structure follows from the fact that both maps
\[
\sHom(M_1,N_2)\ot \sHom(M_2,N_2)\ot \sHom(M_3,N_2)\rightrightarrows \sHom(M_1\ot M_2\ot M_2, N_1\ot N_2\ot N_3)
\]
are adjoint to the map
\[
\begin{array}{l}
\sHom(M_1,N_2)\ot \sHom(M_2,N_2)\ot \sHom(M_3,N_2)\star (M_1\ot M_2\ot M_3)\vspace{2mm}\\
\hspace{1cm}\overset{symm}\longrightarrow (\sHom(M_1,N_2)\star M_1)\ot (\sHom(M_2,N_2)\star M_2)\ot (\sHom(M_3,N_2)\star M_3)\vspace{2mm}\\
\hspace{2cm}\overset{comp}\longrightarrow N_1\ot N_2\ot N_3.
\end{array}
\]
Compatibilities with composition appears as an equality
\begin{equation}\label{eq:1268}
(g_1\ot g_2)\circ (f_1\ot f_2)=(g_1\circ f_1)\ot (g_2\circ f_2):\msc{G}_1\ot \msc{F}_1\ot \msc{G}_2\ot\msc{F}_2\to \sHom(L_1\ot L_2,N_1\ot N_2)
\end{equation}
where $f_i$ and $g_i$ are ``generalized sections", i.e.\ maps
\[
f_i:\msc{F}_i\to \sHom(L_i,M_i),\ \ g_i:\msc{G}_i\to \sHom(M_i,N_i).
\]
To equate these two sections \eqref{eq:1268} one must equate the corresponding morphisms
\[
(\msc{G}_1\ot \msc{F}_1\ot \msc{G}_2\ot \msc{F}_2)\star(L_1\ot L_2)\to N_1\ot N_2,
\]
which involve various applications of the half-braiding for $\QCoh(\dG/\dB)$ acting on $\QCoh(\Xq)$.  One represents these two morphisms via string diagrams and observes the desired equality via naturality of the half-braiding for the $\QCoh(\dG/\dB)$-action.
\end{proof}

\subsection{Proof of Theorem \ref{thm:enhA}}

\begin{lemma}
The adjunction isomorphism
\[
\Hom_{\Xq}(M,N)\to \Hom_{\dG/\dB}(\O_{\dG/\dB},\sHom(M,N))=\Gamma(\dG/\dB,\sHom(M,N))
\]
is precisely the global sections of the natural map \eqref{eq:978}.
\end{lemma}

\begin{proof}
The adjunction map sends a morphism $f:M\to N$ to the unique map $\O_{\dG/\dB}\to \sHom(M,N)$ for which the composite
\[
M=\O_{\dG/\dB}\star M\to \sHom(M,N)\star M\overset{ev}\to N
\]
is the morphism $f$.  Let us call this section $f:\O_{\dG/\dB}\to \sHom(M,N)$.  By considering the fact that we have the surjective map of sheaves
\[
\oplus_{f\in \Hom(M,N)}\O_{\dG/\dB}\star M\overset{\oplus f\ot id}\longrightarrow \Hom(M,N)\ot_k M
\]
we see that the above property implies that the uniquely associated map
\[
\Hom(M,N)\ot_k\O_{\dG/\dB}\to \sHom(M,N)
\]
whose global sections are the adjunction isomorphism fits into the diagram \eqref{eq:982}, and is therefore equal to the morphism \eqref{eq:978}.
\end{proof}

We now prove our theorem.

\begin{proof}[Proof of Theorem \ref{thm:enhA}]
Compatibility with evaluation \eqref{eq:982} implies that the restrictions of the composition and tensor structure on $\sHom_{\Xq}$ along the inclusion
\[
a:\Hom_{\Xq}(M,N)\subset \Hom_{\Xq}(M,N)\ot_k\O_{\Xq}\to \sHom_{\Xq}(M,N)
\]
provided by adjunction recovers the composition and tensor structure maps for $\Hom_{\Xq}(M,N)$.  (Here $\Hom_{\Xq}(M,N)$ denotes the constant sheaf.)  For composition for example we understand, via \eqref{eq:982} and the manner in which composition and evaluation are related for $\Hom_{\Xq}$, that the map
\begin{equation}\label{eq:1022}
\Hom(M,N)\ot_k\Hom(L,M)\overset{\circ}\longrightarrow \Hom(L,N)\overset{a}\to \sHom(L,N)
\end{equation}
is the unique one so that the composite
\[
\Hom(M,N)\ot_k\Hom(L,M)\ot_k L\to \sHom(L,N)\ot L\overset{ev}\to N
\]
is just the squared $k$-linear evaluation map
\[
\Hom(M,N)\ot_k\Hom(L,M)\ot_k L\to \Hom(M,N)\ot_k M\to N
\]
But by \eqref{eq:982} this second map is equal to the composite
\[
\Hom(M,N)\ot_k\Hom(L,M)\ot_k L\overset{a\ot a\ot id}\longrightarrow \sHom(M,N)\ot \sHom(L,M)\ot L\overset{ev^2}\to N
\]
Hence \eqref{eq:1022} is equal to the map
\[
\Hom(M,N)\ot_k\Hom(L,M)\overset{a\ot a}\to \sHom(M,N)\ot \sHom(L,M)\overset{\circ}\to \sHom(L,N),
\]
which just says that restricting along the adjunction map $a$ recovers composition for $\Hom_{\Xq}$ via the global sections of composition for $\sHom_{\Xq}$.  Compatibility with evaluation also implies that the aforementioned map between $\Hom$ spaces respects the $\dG$-action.
\end{proof}

\bibliographystyle{abbrv}

\end{document}